\documentclass[10pt]{amsart}
\usepackage{mathrsfs}
\usepackage{}
\usepackage{txfonts}
\usepackage{amsfonts}
\usepackage{amsfonts}
\usepackage{amssymb,latexsym}
\usepackage{enumerate}
\newtheorem{theorem}{Theorem}
\newcommand{\R}{\mathbb R}

\newtheorem{remark}[theorem]{Remark}
\newtheorem{proposition}[theorem]{Proposition}
\newtheorem{lemma}[theorem]{Lemma}
\newtheorem{definition}[theorem]{Definition}
\newtheorem{corollary}[theorem]{Corollary}

\numberwithin{equation}{section}

\def\a {\alpha}
\def\b{\beta}

\def\n{\nabla}
\begin{document}
\setlength{\baselineskip}{1.2\baselineskip}
\title  [spacetime quasiconcave solutions]
{ON THE microscopic spacetime convexity principle For FULLY NONLINEAR PARABOLIC EQUATIONS II: spacetime quasiconcave solutions }
\author{Chuanqiang Chen}
\address{Wu Wen-Tsun Key Laboratory of Mathematics\\
         University of Science and Technology of China\\
         Hefei, 230026, Anhui Province, CHINA}
\email{cqchen@mail.ustc.edu.cn}
\thanks{Research of the author was supported by the National Natural Science Foundation of China (No. 11301497).}
\maketitle

\begin{abstract}
In \cite{CMS}, Chen-Ma-Salani established the strict convexity of spacetime level sets of solutions to heat equation in convex rings, using the constant rank theorem and a deformation method. In this paper, we generalize the constant rank theorem in \cite{CMS} to fully nonlinear parabolic equations, that is,
establish the corresponding microscopic spacetime convexity principles for spacetime level sets. In fact, the results hold for fully nonlinear parabolic equations under a general structural condition, including the $p$-Laplacian parabolic equations ($p >1$) and some mean curvature type parabolic equations.
\end{abstract}


\section{Introduction}

This paper is devoted to the microscopic spacetime convexity principle for the second
fundamental forms of the spatial and spacetime level sets of the solutions to fully nonlinear
parabolic equations. In this paper, we consider the spatial convexity and the spacetime convexity
of the level sets of the spacetime quasiconcave solutions to the heat
equation. A continuous function $u(x,t)$ on $\Omega \times (0,T]$ is called {\it spacetime quasiconcave}
if the spacetime  superlevel sets $\{(x,t) \in \Omega \times (0, T)| u(x,t) \geq
c\}$ are convex for each constant $c$.

Spacetime convexity is a basic geometric property of the solutions
of parabolic equations.In \cite{Bo82,Bo96,Bo00}, Borell used the
Brownian motion to study certain spacetime convexities of the
solutions of diffusion equations and the level sets of the solution to a heat equations with
Schr$\ddot{o}$dinger potential. Ishige-Salani introduced some notions of parabolic quasiconcavity in \cite{IS10, IS11} and
parabolic power concavity in \cite{IS14}, which are some kinds of spacetime convexity. In \cite{IS10, IS11, IS14}, they studied the corresponding parabolic
boundary value problems using the convex envelope method, which is a macroscopic method. At the same time, Hu-Ma \cite{HM} established a constant rank theorem for the space-time Hessian of space-time convex solutions to the heat equation, which is the microscopic method. Chen-Hu \cite{CH13} generalized the microscopic spacetime convexity principle to fully nonlinear parabolic equations. Recently, Chen-Ma-Salani \cite{CMS} established the strict convexity of spacetime level sets of solutions to the heat equation in convex rings, using the constant rank theorem and a deformation process. In this paper, we generalize the constant rank theorem in \cite{CMS} to fully nonlinear parabolic equations. And the results hold for fully nonlinear parabolic equations under a general structural condition, including the $p$-Laplacian parabolic equations ($p >1$) and some mean curvature type parabolic equations. As in \cite{CMS}, this approach can be used to establish some spacetime convexity of the solutions of some parabolic equations in convex rings, by combining a deformation process.

The convexity of the level sets of the solutions to elliptic partial
differential equations has been studied extensively. For instance,
Ahlfors \cite{AH} contains the well-known result that level curves
of Green function on simply connected convex domain in the plane are
 the convex Jordan curves. In 1956, Shiffman \cite{Sh56} studied the
 minimal annulus in $\mathbb{R}^3$ whose boundary consists of two closed
 convex curves in parallel planes $P_1, P_2$. He proved that the intersection
 of the surface with any parallel plane $P$, between $P_1$ and $P_2$, is a
 convex Jordan curve. In 1957, Gabriel \cite{Ga57} proved that the level sets
 of the Green function on a 3-dimensional bounded convex domain are strictly
 convex. In 1977, Lewis \cite{Le77} extended Gabriel's result to $p$-harmonic
 functions in higher dimensions. Caffarelli-Spruck \cite{CS82} generalized the
 Lewis \cite{Le77} results to a class of semilinear elliptic partial differential equations.
 Motivated by the result of Caffarelli-Friedman \cite{CF85}, Korevaar \cite{Ko90}
 gave a new proof on the results of Gabriel and Lewis by applying the deformation
 process and the constant rank theorem of the second fundamental form of the convex level sets of $p$-harmonic function. A survey of this subject is given by Kawohl \cite{Ka85}. For more recent related extensions, please see the papers by Bianchini-Longinetti-Salani \cite{BLS},  Bian-Guan \cite{BG09}, Xu \cite{Xu08} and Bian-Guan-Ma-Xu \cite{BGMX}.

For the convexity of spacetime level sets, Ishige-Salani \cite{IS10, IS11} introduced some notions of parabolic quasiconcavity, and studied the corresponding parabolic boundary value problems using the convex envelope method. Recently, Chen-Ma-Salani made a break through for the heat equation in \cite{CMS} using the constant rank theorem method.

There is also an extensive literature on the curvature estimates of the level sets of the solutions to elliptic partial differential equations, see \cite{OS}, \cite{Lo83}, \cite{Lo87}, \cite{MOZ09}, \cite{CMY}, \cite{GX} and references therein.

Let us introduce some notations.
\begin{definition}
For each $\theta \in \mathbb S^{n-1}$, denote $\theta^\bot$ the
linear subspace in $\mathbb{R}^n$ which is orthogonal to $\theta$.
Define $\mathcal {S}_n^-(\theta)$ to be the class of $ n\times n$
symmetric real matrices which are negative definite on
$\theta^\bot$. Denote $\mathcal {S}_n^{0-}(\theta)$ the subclass of
$\mathcal {S}_n^-(\theta)$ of matrices that have $\theta$ as
eigenvector with corresponding null eigenvalue. For any $ b \in
\mathbb{R}^n$ with $s = \left\langle {b, \theta} \right\rangle > 0$,
define
\begin{equation*}
\mathscr{B}_\theta ^ - (\Upsilon) = \left\{ {B \in {\mathcal {S}^{n
+ 1}} | B = \left( {\begin{array}{*{20}{c}}
{\widetilde B}&{{b^T}}\\
b&\chi
\end{array}} \right)\quad with \quad \widetilde B \in \mathcal {S}_n^{0-}(\theta)\cap \Upsilon, \chi \in \mathbb{R}} \right\},
\end{equation*}
where $\mathcal {S}^{n+1}$ denote the space of real symmetric $n+1 \times n+1$ matrices,
\end{definition}

Denote $J = (I_n |0)$ the $n \times (n + 1)$ matrix, where $I_n$ is the $n \times n$ identity matrix and
0 is the null vector in $\mathbb{R}^n$.

In this paper, we consider the spacetime quasiconcave solution to fully nonlinear parabolic equation as
follows,
\begin{equation}\label{1.1}
\frac{{\partial u}} {{\partial t}} = F(\n ^2u, \n u, u, x, t),
\quad~\text{in}~ \Omega  \times (0,T],
\end{equation}
where $\Omega$ is a domain in $\R^n$, $\nabla u = (u_{x_1}, \cdots, u_{x_n})$ and $\nabla ^2 u = (u_{x_ix_j})_{n \times n}$. Let
$\Lambda\subset \mathcal {S}^n$ be an open set,  $\mathbb{S}^{n-1}$ a unit sphere and $F=F(r, p, u, x, t)$  a
$C^{2,1}$ function in $\Lambda \times \R^n \times \R \times \Omega \times [0,T] $.
We will assume that $F$ satisfies the following conditions:
\begin{equation}\label{1.2}
\Big( F^{\alpha\beta} \Big) :=\left(\frac{\partial F}{\partial
r_{\alpha\beta}}(\n ^2u, \n u, u, x, t)\right) >0, \quad \forall \; (x, t) \in  \Omega \times [0,T],
\end{equation}
and for each $(\theta,u)\in \mathbb S^{n-1}\times \mathbb R$
fixed,
\begin{equation}\label{1.3}
F(s^{-1} JB^{ - 1} J^T, s^{-1} \theta, u, x, t) \text{ is locally concave in  } (B, x, t).
\end{equation}
In fact, we always assume
\begin{equation}\label{1.4}
|\nabla u| > 0 \text{ and } u_t >0  \text{ in } \Omega \times (0, T].
\end{equation}

Now we state our theorems.
\begin{theorem}\label{th1.2}
Suppose $u \in C^{3,1}(\Omega \times (0,T])$ is a spacetime quasiconcave to fully nonlinear parabolic equation \eqref{1.1}, and $F$ satisfies conditions \eqref{1.2}-\eqref{1.4}.
Then the second fundamental form of spacetime level sets $\hat \Sigma^c=\{(x,t) \in \Omega \times (0, T)| u(x,t) =
c\}$ has the same constant rank in $\Omega$ for each fixed $t \in (0, T]$. Moreover, let $l(t)$ be the minimal
rank of the second fundamental form in $\Omega$, then $l(s)
\leqslant l(t)$ for all $0< s \leqslant t \leqslant T$.
\end{theorem}

For the study of the spacetime level sets of fully nonlinear equation, we should consider the spatial level sets first.
Suppose $u$ is the spacetime quasiconcave solution to fully nonlinear parabolic equation \eqref{1.1}, then $u$ is also spatial quasiconcave, that is the spatial level sets $\Sigma^{c}=\{ x \in \Omega| u(x,t) = c \}$ are all convex. And we get the following constant rank theorem for the second fundamental form of the spatial level sets.
\begin{theorem}\label{th1.3}
Suppose $u \in C^{3,1}(\Omega \times (0,T])$ is a spacetime quasiconcave to fully nonlinear parabolic equation \eqref{1.1}, and $F$ satisfies conditions \eqref{1.2}-\eqref{1.4}.
Then the second fundamental form of spatial level sets $\Sigma^{c} = \{x \in \Omega | u(x,t) =
c\}$ has the same constant rankin $\Omega$ for each fixed $t \in (0, T]$. Moreover, let $l(t)$ be the minimal
rank of the second fundamental form in $\Omega$, then $l(s)
\leqslant l(t)$ for all $0< s \leqslant t \leqslant T$.
\end{theorem}

As it is well known, one needs to choose a
suitable coordinate system to simplify the calculations in the proof of the constant rank theorem. In \cite{CMS}, the proof for the heat equation is based on a coordinate system such that the spatial second fundamental form $a(x,t)$ (see \eqref{2.3}) is diagonalized at each point. In this paper, we generalize the constant rank theorem in \cite{CMS} to fully nonlinear parabolic equations, and we give a more technical proof under the coordinate system such that the spacetime second fundamental form $\hat a(x,t)$ (see \eqref{2.7}) is diagonalized at each point. As in \cite{HM, CH13} and \cite{C}, the key difficulties of two calculations are the same, and the processes of the two proofs are also the same. So the corresponding proof holds for fully nonlinear equations based on the coordinate system such that the spatial second fundamental form $a(x,t)$ is diagonalized at each point, and the calculations must be more complicate than \cite{CMS}.

\begin{remark}
In fact, in the proof of Theorem \ref{th1.3}, we just need a weaker structural condition as follows instead of \eqref{1.3},
\begin{equation}\label{1.5}
F(s^{-1} JB^{ - 1} J^T, s^{-1} \theta, u, x, t) \text{ is locally concave in  } (B, x),  \quad
\text{ for fixed } (\theta,u)\in \mathbb S^{n-1}\times \mathbb R.
\end{equation}
But the condition that $u$ is spacetime quasiconcave is necessary, and it is the main difference between Theorem \ref{th1.3} and the result in \cite{CS}.
That is, if $u$ is spacetime quasiconcave, the constant rank theorem for spatial level sets holds for the parabolic equations with \eqref{1.5}.
Otherwise, if $u$ is spatial quasiconcave, the constant rank theorem for spatial level sets holds for the parabolic equations with a totally different structural condition in \cite{CS}.
\end{remark}

In fact, the $p$-Laplacian operator and the mean curvature operators, that is, the parabolic equations
\begin{align}
\label{1.12} & u_t =\texttt{ div }(|\nabla u|^{p-2} \nabla u), \quad p > 1,  \\
\label{1.13} & u_t = \texttt{div}( \frac{\nabla u}{ \sqrt{1+|\nabla u|^2}} ), \\
\label{1.14}&  u_t =(1+|\nabla u|^2)^{\frac{3}{2}} \texttt{div}( \frac{\nabla u}{ \sqrt{1+|\nabla u|^2}} ),
\end{align}
do not satisfy the structure condition \eqref{1.3} or \eqref{1.5}. But we have the following theorem.

\begin{theorem}\label{th1.7}
Theorem \ref{th1.2} and Theorem \ref{th1.3} holds for the spacetime quasiconcave solutions to the parabolic equations \eqref{1.12}, \eqref{1.13} and \eqref{1.14}.
\end{theorem}

\begin{remark}
Theorem \ref{th1.3} can be looked as a parabolic version for Theorem 1.1 in \cite{BGMX}.
\end{remark}

\begin{remark}
The microscopic spacetime convexity principle can be used to establish the spacetime convexity of the solutions of some parabolic equations in convex rings, by combining the deformation process. For example, Chen-Ma-Salani \cite{CMS} consider the heat equation in convex rings, and get the strict convexity of the  spacetime level sets, with some compatible conditions on the initial data and the convex rings.
\end{remark}

The rest of the paper is organized as follows. Section
2 contains some preliminaries. In Section 3, we prove the constant rank theorem of the spatial second fundamental from, that is Theorem \ref{th1.3} and  the corresponding part of Theorem \ref{th1.7}. For the constant rank theorem of the spacetime second fundamental form is proved in Section 4, including Theorem \ref{th1.2} and the corresponding part of Theorem \ref{th1.7}.

\section{Preliminaries}
\setcounter{equation}{0} \setcounter{theorem}{0}

In this section, we will give some preliminaries.

First, we introduce the definitions of spatial quasiconcave and spacetime quasiconcave.

\begin{definition}\label{def2.1}
A continuous function $u(x,t)$ on $\Omega \times (0,T]$ is called {\it
spatial quasiconcave} if its superlevel sets $\{x \in \Omega| u(x,t)\ge
c\}$ are convex for each constant $c$ and any fixed $t \in (0,T]$.
And $u(x,t)$ is called {\it spacetime quasiconcave} if its superlevel sets $\{(x, t) \in \Omega \times (0,T) | u(x,t)\ge
c\}$ are convex for each constant $c$.
\end{definition}

In the following, we always assume $\nabla u=(u_1, \cdots, u_n)$ is the spatial gradient of
$u$ and $D u=(u_1, \cdots, u_n, u_t)$ the spacetime gradient.

\subsection{Spatial level sets and the spatial second fundamental form }

Suppose $u(x,t) \in C^{2}(\Omega \times (0, T])$, and $u_n \ne 0$
for any fixed $(x,t)\in \Omega \times (0, T]$. It follows that the
upward inner normal direction of the spatial level sets $\Sigma^c = \{x \in
\Omega |u(x,t) = c\}$ is
\begin{eqnarray}\label{2.1}
\vec{\nu}= \frac{|u_n|}{|\nabla u|u_n}(u_1, u_2, \cdots, u_{n-1}, u_n),
\end{eqnarray}
where $\nabla u=(u_1, u_2, \cdots, u_{n-1}, u_n)$ is the spatial gradient of
$u$.

The second fundamental form $II$ of the spatial level sets of
function $u$ with respect to the upward normal direction \eqref{2.1}
is
\begin{equation}
b_{ij} =  - \frac{|u_n|(u_n^2 u_{ij} + u_{nn}u_{i}u_j - u_n u_j u_{in}
- u_n u_i u_{jn})}{|\nabla u| u_n^3}, \quad  1 \leq i,j \leq n-1. \notag
\end{equation}
Set
\begin{equation}\label{2.2}
h_{ij} =u_n^2 u_{ij} + u_{nn}u_{i}u_j - u_nu_ju_{in} - u_n u_i u_{jn}, \quad  1 \leq i,j \leq n-1,
\end{equation}
then we may write
\begin{equation}
b_{ij} = - \frac{|u_n|h_{ij}}{|\nabla u|u_n^3}. \notag
\end{equation}
Note that if $\Sigma ^{c} = \{x \in \Omega |u(x,t) = c\}$ is
locally convex, then the second fundamental form of $\Sigma^{c}$
is semipositive definite with respect to the upward normal direction
\eqref{2.1}. Let $a(x,t)=(a_{ij}(x,t))$ be the symmetric Weingarten
tensor of $\Sigma^{c}= \{x \in \Omega |u(x,t) = c \}$, then $a$ is
semipositive definite. As computed in \cite{BGMX}, if $u_n \ne 0$,
and the Weingarten tensor is
\begin{equation}\label{2.3}
a_{ij} =-\frac{|u_n|}{|\nabla u|{u_n}^3}A_{ij}, \quad  1 \leq i,j \leq n-1,
\end{equation}
where
\begin{equation}\label{2.4}
A_{ij} = h_{ij}
-\frac{u_iu_lh_{jl}}{W(1+W)u_n^2} -\frac{u_ju_lh_{il}}{W(1+W)u_n^2}
+ \frac{u_iu_ju_ku_l h_{kl}}{W^2(1+W)^2u_n^4}, \quad W=\frac{|\nabla u|}{|u_n|}.
\end{equation}
With the above notations, at the
point $(x,t)$ where $u_n(x,t)=|\nabla u(x,t)|>0,\, u_i(x,t)=0$,
$i=1, \cdots, n-1$, $a_{ij,k}$ is commutative, that is, they satisfy
the Codazzi property $a_{ij,k}=a_{ik,j}, \; \forall i,j,k \le n-1$.

\subsection{Spacetime level sets and the spacetime second fundamental form }

Suppose $u(x,t) \in C^{2}(\Omega \times (0, T])$, and $u_t \ne 0$
for any fixed $(x,t)\in \Omega \times (0, T]$. It follows that the
upward inner normal direction of the spatial level sets $\hat \Sigma^c =\{(x,t) \in
\Omega \times (0, T)|u(x,t) = c\}$ is
\begin{eqnarray}\label{2.5}
\vec{\hat{\nu}}= \frac{|u_t|}{|D u|u_t}(u_1, u_2, \cdots, u_{n-1}, u_n, u_t),
\end{eqnarray}
where $D u=(u_1, u_2, \cdots, u_{n-1}, u_n, u_t)$ is the spacetime gradient of
$u$.

The second fundamental form $II$ of the spacetime level sets of
function $u$ with respect to the upward normal direction \eqref{2.5}
is
\begin{equation}
\hat b_{\a \b} =  - \frac{|u_t|(u_t^2 u_{\a \b} + u_{tt}u_{\a}u_{\b} - u_t u_{\b} u_{\a t}
- u_t u_{\a} u_{\b t})}{|D u| u_t^3}, \quad  1 \leq \a, \b \leq n. \notag
\end{equation}
Set
\begin{equation}\label{2.6}
\hat h_{\a \b} =u_t^2 u_{\a \b} + u_{tt}u_{\a}u_{\b} - u_t u_{\b} u_{\a t}
- u_t u_{\a} u_{\b t}, \quad  1 \leq \a, \b \leq n,
\end{equation}
then we may write
\begin{equation}
\hat b_{\a \b} = - \frac{|u_t|\hat h_{\a \b}}{|D u|u_t^3}. \notag
\end{equation}
Note that if $\hat \Sigma ^{c} = \{(x,t) \in \Omega \times (0,  T)|u(x,t) = c\}$ is
locally convex, then the second fundamental form of $\Sigma^{c,t}$
is semipositive definite with respect to the upward normal direction
\eqref{2.5}. Let $\hat a(x,t)=(\hat a_{ij}(x,t))$ be the symmetric Weingarten
tensor of $\hat \Sigma^{c}= \{(x,t) \in \Omega \times (0, T)|u(x,t) = c\}$, then $\hat a$ is
semipositive definite. If $u_t \ne 0$,
and the Weingarten tensor is
\begin{equation}\label{2.7}
\hat a_{\a \b} =-\frac{|u_t|}{|D u|{u_t}^3} \hat A_{\a \b}, \quad  1 \leq \a, \b \leq n,
\end{equation}
where
\begin{equation}\label{2.8}
\hat A_{\a \b} = \hat h_{\a \b}
-\frac{u_\a u_\gamma \hat h_{\b \gamma}}{\hat W(1+\hat W)u_t^2} -\frac{u_\b u_\gamma \hat h_{\a \gamma}}{\hat W(1+ \hat W)u_t^2}
+ \frac{u_\a u_\b u_\gamma u_\eta \hat h_{\gamma \eta}}{\hat W^2(1+\hat W)^2u_t^4}, \quad \hat W = \frac{|D u|}{|u_t|}.
\end{equation}
With the above notations, at the point $(x,t)$ where $u_t(x,t)>0$, $u_n(x,t)=|\nabla u(x,t)|>0,\, u_i(x,t)=0$,
$i=1, \cdots, n-1$, we get
\begin{align}\label{2.9}
1-\frac{u_n^2}{\hat W(1+\hat W)u_t^2} \equiv \frac{\hat W u_t^2+\hat W^2 u_t^2-u_n^2}{\hat W(1+\hat W)u_t^2}
= \frac{1}{\hat W}.
\end{align}
So
\begin{align}
\label{2.10}\hat A_{\alpha \beta} =& \hat h_{\alpha \beta}= u_t^2 u_{\alpha \beta},  \qquad 1 \leq \alpha, \beta \leq n-1;\\
\label{2.11}\hat A_{\alpha n} =& \hat h_{\alpha n}  -\frac{u_n^2\hat h_{\alpha n}}{\hat W(1+\hat W)u_t^2} =\frac{1}{\hat W}\hat h_{\alpha n}
=\frac{1}{\hat W}[ u_t^2 u_{\alpha n}-u_t u_n u_{\alpha t}], \qquad 1 \leq \alpha \leq n-1;  \\
\label{2.12}\hat A_{nn}=&  \hat h_{nn}-2\frac{u_n^2 \hat h_{nn}}{\hat W(1+\hat W)u_t^2}
+ \frac{u_n^4 \hat h_{nn}}{\hat W^2(1+\hat W)^2u_t^4}
= \frac{1}{\hat W^2} \hat h_{nn} \notag \\
= & \frac{1}{\hat W^2} [u_t^2 u_{n n}+ u_n^2 u_{tt} - 2 u_t u_n u_{n t}  ].
\end{align}

Also, at any point $(x,t)$, we can translate the spacetime coordinate systems. When we choose the coordinates $y=(y_1, \cdots, y_n, y_{n+1})$ as a new spacetime coordinates, such that $u_{y_{n+1}} >0$, the Weingarten tensor is
\begin{equation}\label{2.13}
\bar a_{\a \b} =-\frac{|u_{y_{n+1}}|}{|D u|{u_{y_{n+1}}}^3} \bar A_{\a \b}, \quad  1 \leq \a, \b \leq n,
\end{equation}
where
\begin{equation}\label{2.14}
\bar A_{\a \b} = \bar h_{\a \b}
-\frac{u_{y_\a} u_{y_\gamma} \bar h_{\b \gamma}}{\bar W(1+\bar W)u_{y_{n+1}}^2} -\frac{u_{y_\b} u_{y_\gamma} \bar h_{\a \gamma}}{\bar W(1+ \bar W)u_{y_{n+1}}^2}
+ \frac{u_{y_\a} u_{y_\b} u_{y_\gamma} u_{y_\eta} \bar h_{\gamma \eta}}{\bar W^2(1+\bar W)^2 u_{y_{n+1}}^4}, \quad \bar W = \frac{|D u|}{|u_{y_{n+1}}|},
\end{equation}
\begin{equation}\label{2.15}
\bar h_{\a \b} =u_{y_{n+1}}^2 u_{y_\a y_\b} + u_{y_{n+1} y_{n+1}}u_{y_\a} u_{y_\b} - u_{y_{n+1}} u_{y_\b} u_{y_\a y_{n+1}}
- u_{y_{n+1}} u_{y_\a} u_{y_\b y_{n+1}}, \quad  1 \leq \a, \b \leq n.
\end{equation}
With the above notations, at the point $(x,t)$ with the new coordinates $y$ such that $u_{y_i}=0$ for any $1 \leq i \leq n$ and $u_{y_{n+1}}=|Du| >0$,
we get
\begin{equation}\label{2.16}
\bar A_{\a \b} = \bar h_{\a \b}=u_{y_{n+1}}^2 u_{y_\a y_\b} ,\quad  1 \leq \a, \b \leq n,
\end{equation}

\subsection{Elementary symmetric functions}

In this subsection, we recall the definition and some basic properties
of elementary symmetric functions, which could be found in
\cite{L96}.

\begin{definition}
For any $k = 1, 2,\cdots, n,$ we set
\begin{align}
\sigma_k(\lambda) = \sum _{1 \le i_1 < i_2 <\cdots<i_k\leq n}\lambda_{i_1}\lambda_{i_2}\cdots\lambda_{i_k},
 \qquad \text {for any} \quad\lambda=(\lambda_1,\cdots,\lambda_n)\in {\Bbb R}^n.
\end{align}
 We also set $\sigma_0=1$ and $\sigma_k =0$ for $k>n$.
\end{definition}

We denote by $\sigma _k (\lambda \left| i \right.)$ the symmetric
function with $\lambda_i = 0$ and $\sigma _k (\lambda \left| ij
\right.)$ the symmetric function with $\lambda_i =\lambda_j = 0$.

We need the following standard formulas of elementary symmetric
functions.
\begin{proposition}\label{prop2.3}
Let $\lambda=(\lambda_1,\dots,\lambda_n)\in\mathbb{R}^n$ and $k
=0, 1, \cdots, n$, then
\begin{align*}
&\sigma_k(\lambda)=\sigma_k(\lambda|i)+\lambda_i\sigma_{k-1}(\lambda|i), \quad \forall \,1\leq i\leq n,\\
&\sum_{i=1}^n \lambda_i\sigma_{k-1}(\lambda|i)=k\sigma_{k}(\lambda),\\
&\sum_{i=1}^n \sigma_{k}(\lambda|i)=(n-k)\sigma_{k}(\lambda).
\end{align*}
\end{proposition}

The definition can be extended to symmetric matrices by letting
$\sigma_k(W) = \sigma_k(\lambda(W))$, where $ \lambda(W)= (\lambda
_1(W),\lambda _2 (W), \cdots ,\lambda _{n}(W))$ are the eigenvalues
of the symmetric matrix $W$. We also denote by $\sigma _k (W \left|
i \right.)$ the symmetric function with $W$ deleting the $i$-row and
$i$-column and $\sigma _k (W \left| ij \right.)$ the symmetric
function with $W$ deleting the $i,j$-rows and $i,j$-columns. Then
we have the following identities.

\begin{proposition}\label{prop2.4}
Suppose $W=(W_{ij})$ is diagonal, and $m$ is a positive integer,
then
\begin{align}
\frac{{\partial \sigma _m (W)}} {{\partial W_{ij} }} = \begin{cases}
\sigma _{m - 1} (W\left| i \right.), &\text{if } i = j, \\
0, &\text{if } i \ne j.
\end{cases}
\end{align}
and
\begin{align}
\frac{{\partial ^2 \sigma _m (W)}} {{\partial W_{ij} \partial W_{kl}
}} =\begin{cases}
\sigma _{m - 2} (W\left| {ik} \right.), &\text{if } i = j,k = l,i \ne k,\\
- \sigma _{m - 2} (W\left| {ik} \right.), &\text{if } i = l,j = k,i \ne j,\\
0, &\text{otherwise }.
\end{cases}
\end{align}
\end{proposition}

To study the rank of the spacetime second fundamental form $\hat a$, we
need the following simple lemma.
\begin{lemma} \label{lem2.5}
Suppose $ \hat a \geq 0$, and $l=Rank \{\hat a(x_0,t_0) \}\leqslant n-1$, and
$\bigg( \hat a_{\a \b} (x_0,t_0) \bigg)_{n-1 \times n-1} $ is diagonal with $ \hat a_{11} \geq \hat a_{22} \geq \cdots
\geq \hat a_{n-1n-1} $, then at $(x_0,t_0)$, there is a positive constant $
C_0$ such that

CASE 1:
\begin{eqnarray*}
&&\hat a_{11}  \geq \cdots \geq \hat a_{l-1l-1} \geq
C_0 , \quad \hat a_{ll} = \cdots = \hat a_{n-1n-1} =0 , \\
&&\hat a_{nn} -\sum\limits_{i = 1}^{l-1} {\frac{{\hat a_{in} ^2 }} {{\hat a_{ii}
}}} \geq C_0 ,  \quad \hat a_{in} = 0, \quad l \leqslant i \leqslant n-1.
\end{eqnarray*}

CASE 2:
\begin{eqnarray*}
&&\hat a_{11}  \geq \cdots \geq \hat a_{ll} \geq C_0 , \quad \hat a_{l+1l+1} =
\cdots = \hat a_{n-1n-1} =0, \\
&& \hat a_{nn}  = \sum\limits_{i = 1}^{l}{\frac{{\hat a_{in} ^2 }} {{\hat a_{ii} }}}
,  \quad \hat a_{in} = 0, \quad  l+1 \leqslant i \leqslant n-1.
\end{eqnarray*}
\end{lemma}
PROOF. Set $M=\bigg( \hat a_{\a \b}(x_0,t_0) \bigg)_{n-1 \times n-1}=diag(\hat a_{11}, \hat a_{22}, \cdots, \hat a_{n-1n-1}) \geq 0$ and we can assume $Rank \{ M \}= k$ at $(x_0, t_0)$,
then we can obtain $k=l-1$ or $k =l$. Otherwise, if $k <l-1$, we know
\begin{align*}
\hat a_{l-1 l-1}= \cdots=\hat a_{n-1n-1} =0 \quad \text{ at }(x_0, t_0),
\end{align*}
and from $\hat a (x_0,t_0) \geq 0$, we get
\begin{align*}
\hat a_{l-1 n}= \cdots=\hat a_{n-1n} =0  \quad \text{ at }(x_0, t_0).
\end{align*}
So $Rank \{ \hat a \} \leq l-1$, contradiction. If $k>l$, we have
\begin{align*}
l= Rank \{\hat a \} \geq Rank \{ M \} =k \geq l+1  \quad \text{ at }(x_0, t_0).
\end{align*}
This is impossible.

For $k =l-1$, we know at $(x_0, t_0)$
\begin{align*}
\hat a_{11}  \geq \cdots \geq \hat a_{l-1l-1} > 0  , \quad \hat a_{ll} = \cdots = \hat a_{n-1n-1} =0 ,
\end{align*}
and due to $\hat a (x_0,t_0) \geq 0$, we get
\begin{align*}
\hat a_{ln}= \cdots = \hat a_{n-1 n} = 0.
\end{align*}
Since $Rank \{ \hat a \}=l$, then $\sigma_l (\hat a) >0$. Direct computation yields
\begin{align*}
\sigma_l (\hat a) = \hat a_{nn}\sigma_{l-1}(M)
-\sum_{i=1}^{l-1}\hat a_{ni}\hat a_{in}\sigma_{l-2}(M|i) = \sigma_{l-1}(M) [\hat a_{nn} - \sum_{i=1}^{l-1}\frac{\hat a_{in}^2}{\hat a_{ii}} ] >0,
\end{align*}
so we have
\begin{align*}
\hat a_{nn} - \sum_{i=1}^{l-1}\frac{\hat a_{in}^2}{\hat a_{ii}} >0.
\end{align*}
This is CASE 1.

For $k =l$, we know at $(x_0, t_0)$
\begin{align*}
\hat a_{11}  \geq \cdots \geq \hat a_{ll} > 0  , \quad \hat a_{l+1 l+1} = \cdots = \hat a_{n-1n-1} =0 ,
\end{align*}
and due to $\hat a (x_0,t_0) \geq 0$, we get
\begin{align*}
\hat a_{l+1 n}= \cdots = \hat a_{n-1 n} = 0.
\end{align*}
Since $Rank \{ \hat a \}=l$, then $\sigma_{l+1} (\hat a) = 0$. Direct computation yields
\begin{align*}
\sigma_{l+1} (\hat a)= \hat a_{nn}\sigma_{l}(M)
-\sum_{i=1}^{l}\hat a_{ni}\hat a_{in}\sigma_{l-1}(M|i) = \sigma_{l}(M) [\hat a_{nn} - \sum_{i=1}^{l}\frac{\hat a_{in}^2}{\hat a_{ii}} ] =0,
\end{align*}
so we have
\begin{align*}
\hat a_{nn} - \sum_{i=1}^{l}\frac{\hat a_{in}^2}{\hat a_{ii}} =0.
\end{align*}
This is CASE 2. \qed

\subsection{Structural conditions \eqref{1.3} and \eqref{1.5}}

Now we discuss the structural conditions \eqref{1.3}and \eqref{1.5}.

First, we introduce the following set to study the matrix $B^-$.
\begin{definition}
For each $\theta \in \mathbb S^{n-1}$, define $\mathscr{A}_\theta ^ - $ as follows
\begin{equation}
\mathscr{A}_\theta ^ - (\Upsilon) = \left\{ {A \in {\mathcal {S}^{n
+ 1}} | A = \left( {\begin{array}{*{20}{c}}
{\widetilde A}&{{ \mu \theta^T}}\\
\mu \theta & 0
\end{array}} \right)\quad with \quad \widetilde A \in \mathcal {S}_n^{-}(\theta)\cap \Upsilon, \mu >0} \right\}.
\end{equation}
\end{definition}

Properties of  $\mathscr{A}_\theta ^ - $ ,  $\mathscr{B}_\theta ^ - $  and their relationship have been studied in \cite{BLS}.
In particular, if $ \theta = (0, \cdots , 0, 1)$,
\begin{equation}
B = \left( {\begin{array}{*{20}{c}}
{}&{}&{}&0& \times \\
{}&{{a^{ij}}}&{}& \vdots & \vdots \\
{}&{}&{}&0& \times \\
0& \cdots &0&0&s\\
 \times & \cdots & \times &s&\chi
\end{array}} \right) \in \mathscr{B}_\theta ^ - (\Upsilon),
\end{equation}
then
\begin{equation}
A = B^{-1}=\left( {\begin{array}{*{20}{c}}
{}&{}&{}& \times &0\\
{}&{{a_{ij}}}&{}& \vdots & \vdots \\
{}&{}&{}& \times &0\\
 \times & \cdots & \times & \times &\mu \\
0& \cdots &0&\mu &0
\end{array}} \right).
\end{equation}
where the $(n-1)\times (n-1)$ matrix $(a_{ij})$ is negative definite and can be assumed diagonal,
$(a^{ij})$ is the inverse matrix of $(a_{ij})$, $s = B_{n+1,n} = \frac{1}{\mu} > 0$. The values at the positions
denoted by "$\times$" which are not important in the calculations.

For any given $V = ((X_{\alpha\beta} ),Y,(Z_i), D) \in
\mathcal{S}^n \times \mathbb{R} \times \mathbb{R}^{n}\times
\mathbb{R} $, we define a quadratic form

\begin{eqnarray}\label{2.23}
Q^*(V,V) &=& F^{\alpha\beta,\gamma \eta}X_{\alpha\beta}X_{\gamma
\eta}+2F^{\alpha\beta,u_l}
\theta_l X_{\alpha\beta} Y + 2 F^{\alpha\beta,x_k}
X_{\alpha\beta}Z_k  + 2 F^{\alpha\beta,t}
X_{\alpha\beta} D \notag \\
&&+F^{u_k, u_l}\theta_k \theta_l Y^2+ 2 F^{u_k , x_l}
\theta_k Y Z_l  + 2 F^{ u_k ,t} \theta_k Y D +  F^{x_k , x_l}
Z_k Z_l  \nonumber\\
&& + 2 F^{ x_k ,t} Z_k D+  F^{t ,t} D^2 +2 s F^{u_k} \theta_k Y^2 +6 s F^{\alpha\beta}X_{\alpha\beta}Y-6 sF^{\alpha\beta}A_{\alpha\beta}Y^2\nonumber\\
&&+2 s \sum_{i \in T}\frac{F^{\alpha\beta}}{A_{ii}}[X_{i
\alpha}-2A_{i \alpha}Y][X_{i \beta}-2A_{i \beta}Y],
\end{eqnarray}
where the derivative functions of $F$ are evaluated at $( s ^{-1}\widetilde A , s ^{-1} \theta,u,x,t)
$ and $T:= \{ 1, 2, \cdots, n-1\}$.

Through direct calculations, we can get

\begin{lemma} \label{lem2.7}
$F$ satisfies the condition \eqref{1.3} if and only if for each $p \in \mathbb{R}^{n}$
\begin{align} \label{2.24}
Q^*(V,V) \leq 0, \quad \forall \quad
V = ((X_{\alpha\beta} ),Y,(Z_i), D) \in
\mathcal{S}^n \times \mathbb{R} \times \mathbb{R}^{n}\times
\mathbb{R},
\end{align}
where the derivative functions of $F$ are evaluated at $( s ^{-1}\widetilde A, s ^{-1} \theta,u,x,t) $, and $Q^*$ is defined in
\eqref{2.23}.
\end{lemma}

The proof of Lemma \ref{lem2.7} is similar to the discussion in \cite{BG09}, and we omit it.

\begin{remark} \label{rem2.8}
$F$ satisfies the condition \eqref{1.5} if and only if for each fixed $p \in \mathbb{R}^{n}$, and for any $ \tilde V = ((X_{\alpha\beta} ),Y,(Z_i), 0) \in
\mathcal{S}^n \times \mathbb{R} \times \mathbb{R}^{n}\times
\mathbb{R}$
\begin{align} \label{2.25}
Q^*(\tilde V, \tilde V)=& F^{\alpha\beta,\gamma \eta}X_{\alpha\beta}X_{\gamma
\eta}+2F^{\alpha\beta,u_l}
\theta_l X_{\alpha\beta} Y + 2 F^{\alpha\beta,x_k}
X_{\alpha\beta}Z_k +F^{u_k, u_l}\theta_k \theta_l Y^2 \notag \\
&+ 2 F^{u_k , x_l}\theta_k Y Z_l  +  F^{x_k , x_l}
Z_k Z_l +2 s F^{u_k} \theta_k Y^2 +6 s F^{\alpha\beta}X_{\alpha\beta}Y-6 sF^{\alpha\beta}A_{\alpha\beta}Y^2\\
&+2 s \sum_{i \in T}\frac{F^{\alpha\beta}}{A_{ii}}[X_{i
\alpha}-2A_{i \alpha}Y][X_{i \beta}-2A_{i \beta}Y] \notag \\
\leq& 0,\nonumber
\end{align}
where the derivative functions of $F$ are evaluated at $( s ^{-1}\widetilde A, s ^{-1} \theta,u,x,t) $, and $Q^*$ is defined in
\eqref{2.23}. Obviously, the condition \eqref{1.5} is weaker than the condition \eqref{1.3}.
\end{remark}

\subsection{An auxiliary lemma}

Similarly to the Lemma 2.5 in Bian-Guan\cite{BG09}, we have

\begin{lemma}\label{lem2.9}
Suppose $W(x)=(W_{ij}(x))_{N \times N} \geq 0$ for every $x \in \Omega \subset
\mathbb{R}^{n}$, and $W_{ij} (x) \in C^{1,1} (\Omega )$, then for
every $\mathcal {O} \subset \subset \Omega$, there exists a positive
constant $C$ depending only on the $dist\{\mathcal {O},
\partial \Omega \}$ and $\left\| W \right\|_{C^{1,1} (\Omega )}$
such that
\begin{equation}\label{2.26}
\left| {\nabla W_{ij} } \right| \leqslant C(W_{ii} W_{jj} )^{\frac{1}{4}},
\end{equation}
for every $x \in \mathcal {O}$ and $1 \leq i, j \leq N$.
\end{lemma}

\textbf{Proof:} The same arguments as in the proof of Lemma 2.5 in
\cite{BG09} carry through with a small modification since $W$ is a
general matrix instead of a Hessian of a convex function.

It's known that for any nonnegative $C^{1,1}$ function $h$, $|\nabla
h(x)| \leq Ch ^{\frac{1}{2}} (x)$ for all $x \in \mathcal {O}$,
where $C$ depends only on $||h||_{C^{1,1}(\Omega)} $ and
$dist\{\mathcal {O}, \partial \Omega\}$ (see \cite{Tr71}).

Since $W (x) \geq 0$, so we choose $h(x) =W_{ii}(x) \geq 0$. Then we can get from the above argument
\begin{equation*}
\left| {\nabla W_{ii} } \right| \leqslant C_1 (W_{ii})^{\frac{1}{2}} = C_1 (W_{ii} W_{ii} )^{\frac{1}{4}},
\end{equation*}
so \eqref{2.26} holds for $i=j$.

Similarly, for $i\ne j$, we choose $ h=\sqrt{W_{ii}W_{jj}} \geq 0 $, then we get
 \begin{equation}\label{2.27}
\left| {\nabla \sqrt{W_{ii}W_{jj}} } \right| \leqslant C_2 (\sqrt{W_{ii}W_{jj}} )^{\frac{1}{2}} = C_2 (W_{ii} W_{jj} )^{\frac{1}{4}}.
\end{equation}
And for $h=\sqrt{W_{ii}W_{jj}} - W_{ij}$,  we have
 \begin{equation}\label{2.28}
\left| {\nabla (\sqrt{W_{ii}W_{jj}} - W_{ij} )} \right| \leqslant C_3 (\sqrt{W_{ii}W_{jj}} - W_{ij} )^{\frac{1}{2}} \leq C_3 (W_{ii} W_{jj} )^{\frac{1}{4}}.
\end{equation}
So from \eqref{2.27} and \eqref{2.28}, we get
\begin{align*}
\left| {\nabla W_{ij}} \right|  =& \left| {\nabla \sqrt{W_{ii}W_{jj}} } - {\nabla (\sqrt{W_{ii}W_{jj}} - W_{ij} )} \right|  \\
\leq& \left| {\nabla \sqrt{W_{ii}W_{jj}} } \right| + \left| {\nabla (\sqrt{W_{ii}W_{jj}} - W_{ij} )} \right|  \\
\leq& (C_2 +C_3) (W_{ii} W_{jj} )^{\frac{1}{4}}.
\end{align*}
So \eqref{2.26} holds for $i \ne j$. \qed

\begin{remark}\label{rem2.10}
If $W(x,t)=(W_{ij}(x,t))_{N \times N}\geq 0$ for every $(x,t) \in \Omega \times (0, T]$, and $W_{ij} (x,t) \in C^{1,1} (\Omega  \times (0, T])$, then for
every $\mathcal {O} \times (t_0 - \delta, t_0] \subset \subset \Omega \times (0, T]$ with $t_0 < T$, there exists a positive
constant $C$ depending only on the $dist(\mathcal {O} \times (t_0 - \delta, t_0] ,
\partial (\Omega \times (0, T]))$, $t_0$, $\delta$ and $\left\| W \right\|_{C^{1,1} (\Omega \times (0, T] )}$
such that
\begin{equation}\label{2.29}
\left| {D W_{ij} } \right| \leqslant C(W_{ii} W_{jj} )^{\frac{1}{4}},
\end{equation}
for every $(x, t) \in \mathcal {O} \times (t_0 - \delta, t_0]$ and $1 \leq i, j \leq N$. Notice that $D W_{ij} = (\nabla_x W_{ij}, \partial_t W_{ij} )$.
In fact, if $t_0 = T$, it only holds
\begin{equation}\label{2.30}
\left| {\nabla_x W_{ij} } \right| \leqslant C(W_{ii} W_{jj} )^{\frac{1}{4}}.
\end{equation}
for every $(x, t) \in \mathcal {O} \times (t_0 - \delta, t_0]$ and $1 \leq i, j \leq N$.
\end{remark}

\section{Constant rank theorem of the spatial second fundamental form}
\setcounter{equation}{0} \setcounter{theorem}{0}

In this section, we consider the spatial level sets $\Sigma^c = \{ x \in \Omega| u(x,t) = c \}$. Since $u$ is the spacetime quasiconcave solution to fully nonlinear parabolic equation \eqref{1.1}, $u$ is also spatial quasiconcave, that is the spatial level sets $\Sigma^c$ are all convex for $t \in (0, T]$, that is the spatial second fundamental form $a \geq 0$.  We will establish the constant rank theorem
for the spatial second fundamental form $a$ under the structural condition \eqref{1.5} as follows.
\begin{theorem}\label{th3.1}
Suppose $u \in C^{3,1}(\Omega \times (0,T])$ is a spacetime quasiconcave to fully nonlinear parabolic equation \eqref{1.1}, and $F$ satisfies conditions \eqref{1.2}, \eqref{1.4} and \eqref{1.5}.
Then the second fundamental form of spatial level sets $\Sigma^c = \{x \in \Omega | u(x,t) =
c\}$ has the same constant rank in $\Omega$ for any fixed $t \in
(0, T]$. Moreover, let $l(t)$ be the minimal
rank of the second fundamental form in $\Omega$, then $l(s)
\leqslant l(t)$ for all $0< s \leqslant t \leqslant T$.
\end{theorem}
From the discussion in Section 2, the structural condition \eqref{1.5} is weaker than the structural condition \eqref{1.3},
then Theorem \ref{th1.3} holds directly from Theorem \ref{th3.1}.

In the following of this section, we will prove Theorem \ref{th3.1}, and discuss some constant rank properties of the spatial second fundamental form $a$.
And we will prove the constant rank theorem of the spatial fundamental form of the spacetime quasiconcave solutions to the parabolic equations \eqref{1.12}-\eqref{1.14}.

\subsection{Proof of Theorem \ref{th3.1}}

Suppose $a(x,t)$ attains
minimal rank $l$ at some point $(x_0,t_0) \in \Omega \times (0, T]$. We may
assume $l\leqslant n-2$, otherwise there is nothing to prove. And
we assume $u \in C^{4}(\Omega\times (0,T])$ and $u_n (x_0, t_0)>0$. So there is a neighborhood $\mathcal {O}\times
(t_0-\delta, t_0]$ of $(x_0, t_0)$, such that there are $l$
"good" eigenvalues of $(a_{ij})$ which are bounded below by a
positive constant, and the other $n-1-l$ "bad" eigenvalues of
$(a_{ij})$ are very small. Denote $G$ be the index set of these
"good" eigenvalues and $B$ be the index set of "bad" eigenvalues.
And for any fixed point  $(x,t) \in \mathcal {O}\times (t_0-\delta,
t_0]$, we may express $(a_{ij})$ in a form of
\eqref{2.3}, by choosing $e_1,\cdots, e_{n-1},e_n$ such that
\begin{equation}\label{3.1}
 |\n u(x,t)|=u_n(x,t)>0\ \mbox{ and }\
\Big( u_{ij} \Big)_{1\leq i,j \leq n-1} \mbox{ is diagonal at }\ (x,t).
\end{equation}

Without loss of generality we assume $ u_{11} \leq u_{22}\leq \cdots
\leq u_{n-1n-1} $. So, at $(x,t) \in \mathcal {O}\times (t_0-\delta,
t_0]$, from \eqref{2.2}-\eqref{2.4}, we have the matrix
$\Big( a_{ij} \Big)_{1\leq i,j \leq n-1}$ is also diagonal, and $a_{11} \geq a_{22} \geq \cdots \geq a_{n-1,
n-1}$. There is a positive constant $C>0$ depending only on
$\|u\|_{C^{2}}$ and $\mathcal {O}\times (t_0-\delta, t_0]$,
such that $a_{11} \ge a_{22} \ge \cdots \ge a_{ll} > C$ for all $(x,t)
\in \mathcal {O}\times (t_0-\delta, t_0]$. For convenience we
denote $ G = \{ 1, \cdots ,l \} $ and $ B = \{
l+1, \cdots, n-1 \} $ be the "good" and "bad" sets of indices
respectively. If there is no confusion, we also denote
\begin{eqnarray}\label{3.2}
G=\{a_{11}, \cdots, a_{ll}\}, \quad B=\{a_{l+1l+1}, \cdots, a_{n-1n-1}\}.
\end{eqnarray}
Note that for any
$\delta>0$, we may choose $\mathcal {O}\times (t_0-\delta,
t_0]$ small enough such that $a_{jj} <\delta$ for all $j \in
B$ and $(x,t) \in \mathcal {O}\times (t_0-\delta, t_0]$.

For each $c$, let $a=(a_{ij})$ be the symmetric Weingarten tensor of
$\Sigma^{c}$. Set
\begin{eqnarray}\label{3.3}
 p(a)=\sigma_{l+1}(a_{ij}), \quad
q(a) &=& \left\{
 \begin{array}{llr}
\frac{\sigma_{l+2}(a_{ij})}{\sigma_{l+1}(a_{ij})},  &\mbox{if} \; \sigma_{l+1}(a_{ij})>0, &\\
0,& \mbox{otherwise}. &
\end{array}
 \right.
\end{eqnarray}
Since we are dealing with general fully nonlinear equation (\ref{1.1}), as in
the case for the convexity of solutions in \cite{BG09}, there are
technical difficulties to deal with $p(a)$ alone. A key idea in
\cite{BG09} is the introduction of function $q$ as in \eqref{3.3}
and explore some crucial concavity properties of $q$. We consider
function
\begin{equation} \label{3.4}
 \phi(x,t)= p(a)+q(a),
\end{equation}
where $p$ and $q$ as in \eqref{3.3}. We will prove the differential inequality
\begin{equation} \label{3.5}
\sum_{\alpha,\beta=1}^{n}F^{\alpha\beta}\phi_{\alpha\beta}(x,t)-\phi_t\leq
C(\phi+|\nabla \phi|),~~\forall ~(x,t)\in \mathcal {O}
\times(t_0-\delta,t_0],
\end{equation}
where $C$ is a positive constant  independent of $\phi$. Combining with the conditions
\begin{align*}
&\phi \geq 0,  \text{ in } \mathcal {O} \times(t_0-\delta,t_0],  \\
&\phi(x_0,t_0) =0,
\end{align*}
we can get by the strong maximum principle
\begin{align*}
\phi \equiv 0,  \text{ in } \mathcal {O} \times(t_0-\delta,t_0].
\end{align*}
Hence
\begin{align*}
\sigma_{l+1} (a) \equiv 0,  \text{ in } \mathcal {O} \times(t_0-\delta,t_0].
\end{align*}
By the method of continuity, Theorem \ref{th3.1} holds. In the following, we prove the differential inequality \eqref{3.5}.

To get around $\sigma_{l+1}(a) =0$ in $q(a)$, for $\varepsilon>0$ sufficiently small, we instead
consider
\begin{equation}\label{3.6}
\phi_\varepsilon (a)= \phi(a_\varepsilon),
\end{equation}
where $a_\varepsilon=a+\varepsilon I.$ We will also denote
$G_\varepsilon=\{a_{ii}+\varepsilon, i\in G\},$
$B_\varepsilon=\{a_{ii}+\varepsilon, i\in B\}.$

To simplify the notations, we will drop subindex $\varepsilon$
 with the understanding that all the estimates will be independent
 of $\varepsilon.$ In this setting, if we pick $\mathcal {O}\times (t_0-\delta,
t_0]$ small enough, there is $C>0$ independent of $\varepsilon$ such that
\begin{equation}\label{3.7}
\phi(a(x,t))\geq C\varepsilon, \quad \sigma_1(B)\geq C\varepsilon,
~\quad ~\mbox{ for all}\ (x,t) \in \mathcal {O}\times (t_0-\delta,
t_0].
\end{equation}

In the following, we denote
\begin{equation*}
\mathcal H_{\phi} = \sum_{i,j\in B}|\nabla a_{ij}|+\phi.
\end{equation*}
We will use notion $h=O(\mathcal H_{\phi})$ if
$|h(x,t)| \le C\mathcal H_{\phi}$ for $(x,t) \in \mathcal {O}\times (t_0-\delta, t_0]$ with positive constant $C$ under
control.

For any fixed point $(x,t) \in \mathcal
{O}\times (t_0-\delta, t_0]$, we choose a coordinate system as
in (\ref{3.1}) so that $|\n u|=  u_n >0$ and the matrix
$(a_{ij}(x,t))$ is diagonal for $1 \le i,j\le n-1$ and semipositive definite.
From the definition of $\phi$, we get
\begin{align*}
\phi \geq \sigma_l(G) \sum_{i \in B} a_{ii} \geq 0,
\end{align*}
so
\begin{equation}\label{3.8}
a_{ii} =O (\phi )= O(\mathcal H_{\phi}), ~~\forall i\in B.
\end{equation}
And from \eqref{2.2} - \eqref{2.4}, we get
\begin{eqnarray*}
a_{ii}=-\frac{h_{ii}}{u^3_n}=-\frac{u_{ii}}{u_n},
\end{eqnarray*}
so
\begin{equation}\label{3.9}
h_{ii} = O(\mathcal H_{\phi}), u_{ii} = O(\mathcal H_{\phi}), ~~\forall i\in B.
\end{equation}
From the definition of $a_{ij}$, and $u_k = 0$ for $k = 1, \cdots
, n - 1$, we can get
\begin{eqnarray}\label{3.10}
a_{ij,\a} &=&\Big(-\frac{|u_n|}{|\nabla u|{u_n}^3} \Big)_\a h_{ij} +\Big(-\frac{|u_n|}{|\nabla u|{u_n}^3} \Big) h_{ij, \a} \notag \\
&=&3 {u_n}^{-4} u_n^2 u_{ij}- {u_n}^{-3}[u_n^2 u_{ij\a} +2 u_n u_{n\a} u_{ij}-u_{i\a}u_n u_{jn}-u_{j\a}u_n u_{in}] \notag \\
&=&- {u_n}^{-2}[u_n u_{ij\a} -  u_{n\a} u_{ij}-u_{i\a}u_{jn}-u_{j\a} u_{in}],
\end{eqnarray}
so for $i,j \in B$, we get
\begin{eqnarray}
\label{3.11}& u_{ij\alpha } = O(\mathcal H_{\phi}),  ~~\forall \a <n, \\
\label{3.12}&u_n u_{ijn} = 2u_{in}u_{jn}+O(\mathcal H_{\phi}) .
\end{eqnarray}
In fact, from \eqref{2.6}-\eqref{2.8},
\begin{eqnarray*}
\hat a_{jj}=-\frac{|u_t|}{|D u|{u_t}^3} \hat h_{jj}=-\frac{u_{jj}}{|D u|} = O(\mathcal H_{\phi}) , \quad \forall j \in B,
\end{eqnarray*}
and from the spacetime convexity, we can get
\begin{eqnarray*}
\hat a_{jn}^2=\Big[-\frac{|u_t|}{|D u|{u_t}^3} \frac{1}{\hat W} \hat h_{jn} \Big]^2 \leq \hat a_{jj} \hat a_{nn} = O(\mathcal H_{\phi}) , \quad \forall j \in B,
\end{eqnarray*}
so it yields
\begin{eqnarray}\label{3.13}
\hat h_{jn}^2= O(\mathcal H_{\phi}) ,\forall j \in B.
\end{eqnarray}

Following the proof of Lemma 2.1 in \cite{CS}, we can get
\begin{align}\label{3.14}
\phi_t =& \sum_{ij=1}^{n-1}\frac{\partial \phi}{\partial a_{ij}} a_{ij,t}  \notag  \\
=& -u_n^{-3}\sum_{j \in B} \left[ \sigma_l(G) +\frac{{\sigma}^2_1(B|j)-{\sigma}_2(B|j)}
{{\sigma}^2_1(B)}\right][u_n^2u_{jjt}-2u_nu_{jn}u_{jt} ]+O(\mathcal
H_{\phi})
\end{align}
and
\begin{eqnarray}\label{3.15}
&&\sum_{\alpha,\beta=1}^{n}F^{\alpha\beta}\phi_{\alpha\beta} = \sum_{\alpha,\beta=1}^{n}F^{\alpha\beta} \Big[\sum_{ij=1}^{n-1}\frac{\partial \phi}{\partial a_{ij}} a_{ij,\a \beta} + \sum_{ijkl=1}^{n-1}\frac{\partial^2 \phi}{\partial a_{ij}\partial a_{kl}} a_{ij,\a}a_{kl,\b} \Big] \notag \\
&=&u_n^{-3}\sum_{j\in B}\left[ \sigma_l(G) +
\frac{{\sigma}^2_1(B|j)-{\sigma}_2(B|j)} {{\sigma}^2_1(B)}\right]
[-u_n^2\sum_{\alpha,\beta=1}^n F^{\alpha\beta}u_{\alpha\beta j
j}+6u_nu_{nj}\sum_{\alpha,\beta=1}^nF^{\alpha\beta}u_{\alpha\beta j}
\notag \\
&& \qquad\qquad\qquad\qquad\qquad\qquad\qquad\qquad\quad
-6u_{nj}^2\sum_{\alpha,\beta=1}^{n}F^{\alpha\beta}u_{\alpha\beta
}]  \notag  \\
&&+2u_n^{-3}\sum_{j\in B,i\in G}\left[ \sigma_l(G) +
\frac{{\sigma}^2_1(B|j)-{\sigma}_2(B|j)}
{{\sigma}^2_1(B)}\right]\sum_{\alpha,\beta=1}^{n}
F^{\alpha\beta}\frac{1}{u_{ii}}[u_nu_{ij\alpha}-2u_{i\a}u_{jn}][u_nu_{ij\b}-2u_{i\b}u_{jn}] \notag \\
&&-\frac{1}{{\sigma}^3_1(B)}\sum_{\alpha,\beta=1}^n\sum_{i\in
B}F^{\alpha\beta}[{\sigma}_1(B)a_{ii,\alpha}-a_{ii}\sum_{j\in
B}a_{jj,\alpha}][{\sigma}_1(B)a_{ii,\beta}-a_{ii}\sum_{j\in
B}a_{jj,\beta}]\\
&& -\frac{1}{{\sigma}_1(B)}\sum_{\alpha,\beta=1}^n\sum_{i\neq j\in
B}F^{\alpha\beta}a_{ij,\alpha}a_{ij,\beta}+O(\mathcal H_{\phi}).
\notag
\end{eqnarray}
Hence
\begin{eqnarray}\label{3.16}
&&\sum_{\alpha,\beta=1}^nF^{\alpha\beta}\phi_{\alpha\beta}-\phi_t\nonumber\\
&=& -u_n^{-3}\sum_{j\in
B}\left[\sigma_l(G)+\frac{{\sigma}^2_1(B|j)-{\sigma}_2(B|j)}{{\sigma}^2_1(B)}\right]\Big[
u_n^2(\sum_{\alpha,\beta=1}^nF^{\alpha\beta}u_{jj\alpha\beta}-u_{jjt})+ 2u_nu_{jn}u_{jt}
\nonumber\\
&&\qquad \qquad \qquad \qquad \qquad \qquad\qquad \qquad  -6u_nu_{jn}\sum_{\alpha,\beta=1}^n
F^{\alpha\beta}u_{j\alpha\beta} +6u_{jn}^2\sum_{\alpha,\beta=1}^n
F^{\alpha\beta}u_{\alpha\beta}\Big]\nonumber\\
&&+2u_n^{-3}\sum_{j\in B,i\in G}\left[ \sigma_l(G) +
\frac{{\sigma}^2_1(B|j)-{\sigma}_2(B|j)}
{{\sigma}^2_1(B)}\right]\sum_{\alpha,\beta=1}^{n}
F^{\alpha\beta}\frac{1}{u_{ii}}[u_nu_{ij\alpha}-2u_{i\a}u_{jn}][u_nu_{ij\b}-2u_{i\b}u_{jn}] \notag \\
&&-\frac{1}{{\sigma}^3_1(B)}\sum_{\alpha,\beta=1}^n\sum_{i\in
B}F^{\alpha\beta}[{\sigma}_1(B)a_{ii,\alpha}-a_{ii}\sum_{j\in
B}a_{jj,\alpha}][{\sigma}_1(B)a_{ii,\beta}-a_{ii}\sum_{j\in
B}a_{jj,\beta}]\\
&&-\frac{1}{{\sigma}_1(B)}\sum_{\alpha,\beta=1}^n\sum_{i\neq j,
i,j\in B}F^{\alpha\beta}a_{ij,\alpha}a_{ij,\beta}+O(\mathcal H_{\phi}).\nonumber
\end{eqnarray}

From $ \hat h_{jn} = u_t^2 u_{jn} - u_n u_t u_{jt}$, we have
\begin{align}\label{3.17}
2u_nu_{jn}u_{jt} =& - \frac{1}{u_t} [( \frac{\hat h_{jn}}{u_t})^2 -( u_t u_{jn})^2 - (u_n u_{jt})^2  ] \notag \\
=& O(\mathcal H_{\phi}) + \frac{1}{u_t} [( u_t u_{jn})^2 + (u_n u_{jt})^2  ],
\end{align}
where the second "=" holds from \eqref{3.13}.

For each $j\in B$, differentiating equation (\ref{1.1}) in $e_j$ direction at
$x$,
\begin{align}\label{3.18}
u_{jjt} =& \sum\limits_{kl = 1}^n {F^{kl} } u_{kljj}  + \sum\limits_{i = 1}^n {F^{u_i } } u_{jji}  + F^u u_{jj} \notag  \\
&+ \sum\limits_{klpq = 1}^n {F^{kl,pq} u_{klj} u_{pqj} }  + 2\sum\limits_{ikl = 1}^n {F^{kl,u_i } } u_{klj} u_{ij}  + 2\sum\limits_{kl = 1}^n {F^{kl,u} } u_{klj} u_j  \notag \\
&+ 2\sum\limits_{kl = 1}^n {F^{kl,x_j } } u_{klj}  + \sum\limits_{ik = 1}^n {F^{u_i ,u_k } } u_{ij} u_{kj}  + 2\sum\limits_{i = 1}^n {F^{u_i ,u} } u_{ij} u_j  + 2\sum\limits_{i = 1}^n {F^{u_i ,x_j } } u_{ij}  \notag \\
&+ F^{u,u} u_j ^2  + 2F^{u,x_j } u_j  + F^{x_j ,x_j } \notag \\
=& \sum\limits_{kl = 1}^n {F^{kl} } u_{kljj}  + 2 \frac{F^{u_n} } {u_n} u_{jn}^2  + \sum\limits_{klpq = 1}^n {F^{kl,pq} u_{klj} u_{pqj} }  + 2\sum\limits_{kl = 1}^n {F^{kl,u_n } } u_{klj} u_{jn}   \notag \\
&+ 2\sum\limits_{kl = 1}^n {F^{kl,x_j } } u_{klj}  +  {F^{u_n,u_n } } u_{jn}^2   + 2{F^{u_n ,x_j } } u_{jn}  + F^{x_j ,x_j }+O(\mathcal H_{\phi}).
\end{align}
Set
\begin{align}\label{3.19}
Q_j=& \sum\limits_{klpq = 1}^n {F^{kl,pq} u_{klj} u_{pqj} u_n^2 }  + 2\sum\limits_{kl = 1}^n {F^{kl,u_n } } u_{klj} u_{jn}u_n^2 + 2\sum\limits_{kl = 1}^n {F^{kl,x_j } } u_{klj} u_n^2 +  {F^{u_n,u_n } } u_{jn}^2 u_n^2   \notag  \\
&  + 2{F^{u_n ,x_j } } u_{jn} u_n^2  + F^{x_j ,x_j }u_n^2 + 2 F^{u_n}u_n u_{jn}^2 +6u_nu_{jn}\sum_{\alpha,\beta=1}^n
F^{\alpha\beta}u_{j\alpha\beta} - 6u_{jn}^2\sum_{\alpha,\beta=1}^n F^{\alpha\beta}u_{\alpha\beta}  \\
&+2\sum_{i\in G}\sum_{\alpha,\beta=1}^{n}
F^{\alpha\beta}\frac{1}{u_{ii}}[u_nu_{ij\alpha}-2u_{i\a}u_{jn}][u_nu_{ij\b}-2u_{i\b}u_{jn}],\notag
\end{align}
and denote
\begin{align*}
&s=\frac{1}{u_{n}} = \frac{1}{|\nabla u|},  A_{ij} = s u_{ij} = \frac{u_{ij}}{u_{n}},\theta=(0, 0, \cdots, 0,1);  \\
&X_{\alpha\beta}=2 u_{\a \b}u_{j n}, \quad \a \in B \text{ or } \b \in B;  \\
&X_{\alpha\beta}=u_{\alpha \beta j} u_{n}, \quad \text{ otherwise };  \\
&Y=u_{jn} u_{n};  \\
&Z_i= \delta_{ij}, \quad i =1,2, \cdots, n;  \\
&\tilde V = ((X_{\alpha\beta} ),Y,(Z_i), 0)\in \mathcal{S}^n \times \mathbb{R} \times \mathbb{R}^{n}\times
\mathbb{R} ;
\end{align*}
then we can get
\begin{align*}
X_{i \alpha} -2 A_{i \a} Y =0,  \quad i \in B.
\end{align*}
So it yields
\begin{align} \label{3.20}
Q_j= Q^*(\tilde V, \tilde V),
\end{align}
where $Q^*(\tilde V, \tilde V)$ is defined in \eqref{2.25}.

From \eqref{3.16} - \eqref{3.18}), it yields
\begin{eqnarray}\label{3.21}
&&F^{\alpha\beta}\phi_{\alpha\beta}- \phi_t\nonumber\\
&=& u_n^{-3}\sum_{j\in
B}\left[\sigma_l(G)+\frac{{\sigma}^2_1(B|j)-{\sigma}_2(B|j)}{{\sigma}^2_1(B)}\right] \Big(  Q_j -2 u_n u_{jn} u_{jt} \Big) \nonumber\\
&&-\frac{1}{{\sigma}^3_1(B)}\sum_{\alpha,\beta=1}^n\sum_{i\in
B}F^{\alpha\beta}[{\sigma}_1(B)a_{ii,\alpha}-a_{ii}\sum_{j\in
B}a_{jj,\alpha}][{\sigma}_1(B)a_{ii,\beta}-a_{ii}\sum_{j\in
B}a_{jj,\beta}]\nonumber\\
&&-\frac{1}{{\sigma}_1(B)}\sum_{\alpha,\beta=1}^n\sum_{i\neq j,
i,j\in B}F^{\alpha\beta}a_{ij,\alpha}a_{ij,\beta}+O(\mathcal H_{\phi}) \notag \\
&\leq& u_n^{-3}\sum_{j\in
B}\left[\sigma_l(G)+\frac{{\sigma}^2_1(B|j)-{\sigma}_2(B|j)}{{\sigma}^2_1(B)}\right] Q_j\nonumber\\
&&-\frac{1}{{\sigma}^3_1(B)}\sum_{\alpha,\beta=1}^n\sum_{i\in
B}F^{\alpha\beta}[{\sigma}_1(B)a_{ii,\alpha}-a_{ii}\sum_{j\in
B}a_{jj,\alpha}][{\sigma}_1(B)a_{ii,\beta}-a_{ii}\sum_{j\in
B}a_{jj,\beta}]\nonumber\\
&&-\frac{1}{{\sigma}_1(B)}\sum_{\alpha,\beta=1}^n\sum_{i\neq j,
i,j\in B}F^{\alpha\beta}a_{ij,\alpha}a_{ij,\beta}+O(\mathcal H_{\phi}).
\end{eqnarray}

From the structural condition \eqref{1.5} (i.e. Remark \ref{rem2.8}), it implies
\begin{equation*}
Q^*(\tilde V, \tilde V)\le 0.
\end{equation*}
so for $j \in B$, we get
\begin{align}\label{3.22}
Q_j= Q^*(\tilde V, \tilde V) \leq 0.
\end{align}
Condition \eqref{1.2} implies
\begin{equation}\label{3.23}
 (F^{\alpha\beta}) \ge \delta_0 I_n, \; \quad \mbox{for some $\delta_0>0$, and
$\forall x \in \mathcal {O}$.}
\end{equation}
Set
\[
V_{i\alpha}={\sigma}_1(B)a_{ii,\alpha}-a_{ii}\sum_{j\in B}a_{jj,\alpha}.
\]
Combining (\ref{3.21}), (\ref{3.22}) and (\ref{3.23}),
\begin{eqnarray}\label{3.24}
F^{\alpha\beta}\phi_{\alpha\beta}\le  C(\phi+\sum_{i,j\in B}|\nabla
a_{ij}|)-\delta_0[\frac{\sum_{i\neq j\in
B, \alpha=1}^na^2_{ij,\alpha}}{\sigma_1(B)}+\frac{\sum_{i\in
B,\alpha=1}^nV_{i\alpha}^2}{{\sigma}^3_1(B)}].
\end{eqnarray}
By Lemma~3.3 in \cite{BG09}, for each $M\ge 1$, for any $M\ge
|\gamma_i|\ge \frac{1}{M}$, there is a constant $C$ depending only
on $n$ and $M$ such that, $\forall \alpha$,
\begin{equation}\label{3.25}
\sum_{i,j\in B}|a_{ij,\alpha}|\le C(1+\frac{1}{\delta_0^2})
(\sigma_1(B)+|\sum_{i\in B} \gamma_ia_{ii,\alpha}|)
+\frac{\delta_0}{2}[\frac{\sum_{i\neq j \in B}|a_{ij,\alpha}|^{2}}
{\sigma_{1}(B)}+
\frac{\sum_{i\in B}V_{i\alpha}^{2}}{\sigma_{1}^{3}(B)}].
\end{equation}
Taking $\gamma_i=\sigma_l(G)+\frac{{\sigma}^2_1(B|i)-{\sigma}_2(B|i)}{{\sigma}^2_1(B)}$ for each $i\in B$,
the Newton-MacLaurine inequality implies
\[
\sigma_l(G)+1 \ge
\sigma_l(G)+\frac{{\sigma}^2_1(B|j)-{\sigma}_2(B|j)}{{\sigma}^2_1(B)}
\ge \sigma_l(G), \quad \forall j\in B.
\]
and
\begin{equation}\label{3.26}
\phi_{\alpha}=\sum_{ij=1 }^{n-1} \frac{\partial \phi}{\partial a_{ij}} a_{ij, \a} = \sum_{i\in B} \gamma_ia_{ii,\alpha} + O(\phi).
\end{equation}
Therefore we conclude from \eqref{3.25} and (\ref{3.26}) that
$\sum_{i,j \in B} |\nabla a_{ij}|$  can be controlled by the rest
terms on the right hand side in (\ref{3.24}) and $\phi + |\nabla
\phi|$. So \eqref{3.5} holds, and the proof of Theorem \ref{th3.1} is complete.
\qed

\subsection{Constant rank properties of $a$}

In the proof of Theorem \ref{th3.1}, we can get for any $(x,t)\in \mathcal {O}
\times(t_0-\delta,t_0]$ with the suitable coordinate \eqref{3.1},
\begin{eqnarray}\label{3.27}
&&\sum_{\alpha,\beta=1}^nF^{\alpha\beta}\phi_{\alpha\beta}-\phi_t\nonumber\\
&=& u_n^{-3}\sum_{j\in B,i\in G}\left[ \sigma_l(G) +
\frac{{\sigma}^2_1(B|j)-{\sigma}_2(B|j)}
{{\sigma}^2_1(B)}\right]\bigg(Q_j - \frac{1}{u_t} [( u_t u_{jn})^2 + (u_n u_{jt})^2  ]\bigg) \notag \\
&&-\frac{1}{{\sigma}^3_1(B)}\sum_{\alpha,\beta=1}^n\sum_{i\in
B}F^{\alpha\beta}[{\sigma}_1(B)a_{ii,\alpha}-a_{ii}\sum_{j\in
B}a_{jj,\alpha}][{\sigma}_1(B)a_{ii,\beta}-a_{ii}\sum_{j\in
B}a_{jj,\beta}]\\
&&-\frac{1}{{\sigma}_1(B)}\sum_{\alpha,\beta=1}^n\sum_{i\neq j,
i,j\in B}F^{\alpha\beta}a_{ij,\alpha}a_{ij,\beta}+O(\mathcal H_{\phi})\nonumber \\
&\leq& C( \phi + |\nabla \phi |),\nonumber
\end{eqnarray}
and by the strong maximum principle,
\begin{align}\label{3.28}
\phi =0  \quad \text{ for } (x,t) \in \mathcal {O} \times (t_0-\delta,t_0].
\end{align}
So it must have  for any $(x,t)\in \mathcal {O}
\times(t_0-\delta,t_0]$ with the suitable coordinate \eqref{3.1}
\begin{align}\label{3.29}
a_{jj} = 0,\quad  \text{ for } j \in B.
\end{align}

In fact, we can get more information from the differential inequality, and the constant rank properties is as follows

\begin{corollary} \label{cor3.2}
For any $(x,t)\in \mathcal {O} \times(t_0-\delta,t_0)$ with the suitable coordinate \eqref{3.1}
\begin{align}
\label{3.30}&u_{jn}=u_{jt} =0, |D u_{j}|=0, \quad  \text{ for } j \in B;  \\
\label{3.31}&\sum\limits_{kl = 1}^n {F^{kl} } u_{kljj} - u_{jjt}= 2\sum_{i\in G}\sum_{\alpha,\beta=1}^{n}
F^{\alpha\beta}\frac{u_n^2 u_{ij\alpha} u_{ij\b} }{u_{ii}},  \quad  \text{ for } j \in B; \\
\label{3.32}&|D u_{ij}|=0, \quad \text{ for }  i \in B, j= 1, 2, \cdots, n-1.
\end{align}

\end{corollary}

PROOF. For $(x,t)\in \mathcal {O} \times(t_0-\delta,t_0)$ with the suitable coordinate \eqref{3.1}, we have from \eqref{2.3} and
\eqref{3.29}
\begin{align*}
u_{jj}=0, \text{ for } j \in B.
\end{align*}
and from \eqref{3.28} and \eqref{3.27}, we get for $j \in B$,
\begin{align*}
&Q_{j}=0,  \\
&u_t u_{jn}=u_n u_{jt} =0.
\end{align*}
So
\begin{align} \label{3.33}
u_{jn}=u_{jt} =0,  \quad  \text{ for } j \in B,
\end{align}
then
\begin{align}\label{3.34}
 |D u_{j}|=0, \quad  \text{ for } j \in B.
\end{align}
From the definition of $Q_j$, and \eqref{3.28}, \eqref{3.32}, \eqref{3.33}, we get
\begin{eqnarray*}
0= Q_j &=&\sum\limits_{klmn = 1}^n {F^{kl,mn} u_{klj} u_{mnj} u_{n}^2 }   + 2 \sum\limits_{kl = 1}^n {F^{kl,u_n } } u_{klj} u_{jn} u_{n}^2 + 2 \sum\limits_{kl = 1}^n {F^{kl,x_j } } u_{klj}  u_{n}^2 \notag \\
&&+  {F^{u_n,u_n } } u_{jn}^2  u_{n}^2   + 2{F^{u_n ,x_j } } u_{jn} u_{n}^2   + F^{x_j ,x_j } u_{n}^2  +2\sum_{i\in G}\sum_{\alpha,\beta=1}^{n}
F^{\alpha\beta}\frac{u_n^2 u_{ij\alpha} u_{ij\b} }{u_{ii}} \notag \\
&=& u_{jjt}-\sum\limits_{kl = 1}^n {F^{kl} } u_{kljj}  + 2\sum_{i\in G}\sum_{\alpha,\beta=1}^{n}
F^{\alpha\beta}\frac{u_n^2 u_{ij\alpha} u_{ij\b} }{u_{ii}},  \quad  \text{ for } j \in B.
\end{eqnarray*}
Also, we can get from \eqref{3.28}, and Lemma \ref{lem2.9} (i.e. Remark \ref{2.10})
\begin{align*}
|D a_{ij}|=0, \quad  \text{ for } i \in B, j =1, 2, \cdots, n-1,
\end{align*}
then from \eqref{3.10} and \eqref{3.32}
\begin{align*}
|D u_{ij}|=0, \quad  \text{ for } i \in B, j =1, 2, \cdots, n-1.
\end{align*}
So the proof is complete.
\qed

\subsection{Constant rank theorem of the spatial fundamental form for the equation \eqref{1.12}}

In this subsection, we consider the $p$-Laplacian parabolic equation, that is
\begin{align}\label{3.35}
 u_t  = \texttt{div}(|\nabla u|^{p - 2} \nabla u) = L_{\alpha \beta} (\nabla u)u_{\alpha \beta}, \text{ in } \Omega \times (0, T],
 \end{align}
where
\begin{align}\label{3.36}
 L_{\alpha \beta} (\nabla u) = |\nabla u|^{p - 2} \delta _{\alpha \beta}  + (p - 2)|\nabla u|^{p - 4} u_\alpha  u_\beta, \quad 1 \leq \alpha, \beta \leq n.
\end{align}
It is easy to know the equation \eqref{3.35} is parabolic when $p>1$ and $|\nabla u | >0$ in $\Omega \times [0, T]$.
We will establish the constant rank theorem
for the spatial second fundamental form $a$ as follows.
\begin{theorem}\label{th3.3}
Suppose $u \in C^{3,1}(\Omega \times (0,T])$ is a spacetime quasiconcave to the parabolic equation \eqref{3.35} and satisfies \eqref{1.4}.
Then the second fundamental form of spatial level sets $\Sigma^c = \{x \in \Omega | u(x,t) =
c\}$ has the same constant rank in $\Omega$ for any fixed $ t \in (0, T]$. Moreover, let $l(t)$ be the minimal
rank of the second fundamental form in $\Omega$, then $l(s)
\leqslant l(t)$ for all $0< s \leqslant t \leqslant T$.
\end{theorem}

PROOF. The proof is similar to the the proof of Theorem \ref{th3.1}, with some modifications.

Suppose $a(x,t)$ attains
minimal rank $l$ at some point $(x_0,t_0) \in \Omega \times (0, T]$. We may
assume $l\leqslant n-2$, otherwise there is nothing to prove. So there is a small neighborhood $\mathcal {O}\times
(t_0-\delta, t_0]$ of $(x_0, t_0)$, such that there are $l$
"good" eigenvalues of $(a_{ij})$ which are bounded below by a
positive constant, and the other $n-1-l$ "bad" eigenvalues of
$(a_{ij})$ are very small. Denote $G$ be the index set of these
"good" eigenvalues and $B$ be the index set of "bad" eigenvalues.
We will prove the differential inequality
\begin{equation} \label{3.37}
\sum_{\alpha,\beta=1}^{n}L_{\alpha\beta}\phi_{\alpha\beta}(x,t)-\phi_t\leq
C(\phi+|\nabla \phi|),~~\forall ~(x,t)\in \mathcal {O}
\times(t_0-\delta,t_0],
\end{equation}
where $\phi$ is defined in \eqref{3.4} and $C$ is a positive constant  independent of $\phi$.
Then by the strong maximum principle and the method of continuity, Theorem \ref{th3.3} holds.

For any fixed point  $(x,t) \in \mathcal {O}\times (t_0-\delta,
t_0]$, we may express $(a_{ij})$ in a form of
\eqref{2.3}, by choosing $e_1,\cdots, e_{n-1},e_n$ such that
\begin{equation}\label{3.38}
 |\n u(x,t)|=u_n(x,t)>0\ \mbox{and} \
\Big(u_{ij} \Big)_{1 \leq i,j \leq n-1} \mbox{is diagonal at}\ (x,t).
\end{equation}
Following the proof of Theorem \ref{th3.1}, we get from \eqref{3.21}
\begin{eqnarray}\label{3.39}
&&L_{\alpha\beta}\phi_{\alpha\beta}- \phi_t\nonumber\\
&=& u_n^{-3}\sum_{j\in
B}\left[\sigma_l(G)+\frac{{\sigma}^2_1(B|j)-{\sigma}_2(B|j)}{{\sigma}^2_1(B)}\right] \Big(  Q_j -2 u_n u_{jn} u_{jt} \Big) \nonumber\\
&&-\frac{1}{{\sigma}^3_1(B)}\sum_{\alpha,\beta=1}^n\sum_{i\in
B} L_{\alpha\beta}[{\sigma}_1(B)a_{ii,\alpha}-a_{ii}\sum_{j\in
B}a_{jj,\alpha}][{\sigma}_1(B)a_{ii,\beta}-a_{ii}\sum_{j\in
B}a_{jj,\beta}]\nonumber\\
&&-\frac{1}{{\sigma}_1(B)}\sum_{\alpha,\beta=1}^n\sum_{i\neq j,
i,j\in B}L_{\alpha\beta}a_{ij,\alpha}a_{ij,\beta}+O(\mathcal H_{\phi}),
\end{eqnarray}
where
\begin{align*}
Q_j=&  2\sum\limits_{kl = 1}^n \frac{{\partial L_{kl} }}{{\partial u_n }} u_{klj} u_{jn}u_n^2 +  \sum\limits_{kl = 1}^n \frac{{\partial ^2 L_{kl} }}{{\partial u_n ^2 }} u_{kl} u_{jn}^2 u_n^2    \\
&  + 2 \sum\limits_{kl = 1}^n \frac{{\partial L_{kl} }}{{\partial u_n }} u_{kl}u_n u_{jn}^2 +6u_nu_{jn}\sum_{kl=1}^n
L_{kl}u_{jkl} - 6u_{jn}^2\sum_{kl=1}^n L_{kl}u_{kl}  \\
&+2\sum_{i\in G}\sum_{\alpha,\beta=1}^{n}
\frac{1}{u_{ii}}L_{\alpha\beta}[u_nu_{ij\alpha}-2u_{i\a}u_{jn}][u_nu_{ij\b}-2u_{i\b}u_{jn}].
\end{align*}

Under the coordinate \eqref{3.38}, we get
\begin{align}\label{a1}
L_{kl}  = 0 , k \ne l ; \quad L_{kk}  = u_n ^{p - 2} , k < n ;  \quad L_{nn}  = (p - 1)u_n ^{p - 2};
\end{align}

\begin{align}
 \label{a2}&\frac{{\partial L_{kl} }}{{\partial u_n }} = 0, \quad k \ne l; \\
 \label{a3}&\frac{{\partial L_{kk} }}{{\partial u_n }} = (p - 2)u_n ^{p - 3}  = (p - 2)\frac{{L_{kk} }}{{u_n }}, \quad k < n;  \\
 \label{a4}&\frac{{\partial L_{nn} }}{{\partial u_n }} = (p - 1)(p - 2)u_n ^{p - 3}  = (p - 2)\frac{{L_{nn} }}{{u_n }}.
\end{align}
and
\begin{align}
\label{a5}&\frac{{\partial ^2 L_{kl} }}{{\partial u_n ^2 }} = 0, k \ne l; \\
\label{a6}&\frac{{\partial ^2 L_{kk} }}{{\partial u_n ^2 }} = (p - 2)(p - 3)u_n ^{p - 4}  = (p - 2)(p - 3)\frac{{L_{kk} }}{{u_n ^2 }}, k < n; \\
\label{a7}&\frac{{\partial ^2 L_{nn} }}{{\partial u_n ^2 }} = (p - 1)(p - 2)(p - 3)u_n ^{p - 4}  = (p - 2)(p - 3)\frac{{L_{nn} }}{{u_n ^2 }}.
\end{align}

From the equation \eqref{3.35}, we know
\begin{align*}
 u_t  = L_{kk} u_{kk},
\end{align*}
and for $j \in B$
\begin{align}
u_{tj}  =& L_{kk} u_{kkj}  + \frac{{\partial L_{kl} }}{{\partial u_p }}u_{pj} u_{kl}= L_{kk} u_{kkj}  + \frac{{\partial L_{kl} }}{{\partial u_j }}u_{jj} u_{kl}+ \frac{{\partial L_{kl} }}{{\partial u_n }}u_{nj} u_{kl} \notag\\
=& L_{kk} u_{kkj}  +O(\mathcal H_{\phi})+ \frac{{\partial L_{kk} }}{{\partial u_n }}u_{nj} u_{kk}  = L_{kk} u_{kkj}  + (p - 2)\frac{{L_{kk} }}{{u_n }}u_{nj} u_{kk} +O(\mathcal H_{\phi})\notag \\
=& L_{kk} u_{kkj}  + (p - 2)\frac{{u_t }}{{u_n }}u_{nj}+O(\mathcal H_{\phi}).
\end{align}
And from \eqref{3.13}, we get
\begin{eqnarray}
\hat h_{jn}^2= O(\mathcal H_{\phi}) ,\forall j \in B.
\end{eqnarray}

So
\begin{align*}
Q_j=&  2\sum\limits_{k = 1}^n (p - 2)\frac{{L_{kk} }}{{u_n }} u_{kkj} u_{jn}u_n^2 +  \sum\limits_{k = 1}^n (p - 2)(p - 3)\frac{{L_{kk} }}{{u_n ^2 }} u_{kk} u_{jn}^2 u_n^2    \\
&  + 2 \sum\limits_{k = 1}^n (p - 2)\frac{{L_{kk} }}{{u_n }} u_{kk}u_n u_{jn}^2 +6u_nu_{jn}\sum_{k=1}^n
L_{kk}u_{kkj} - 6u_{jn}^2 u_t  \\
&+2\sum_{i\in G}\sum_{\alpha=1}^{n}
\frac{1}{u_{ii}}L_{\alpha\a}[u_nu_{ij\alpha}-2u_{i\a}u_{jn}]^2  \\
=&  2(p - 2)[u_{tj}  - (p - 2)\frac{{u_t }}{{u_n }}u_{nj}] u_{jn}u_n +  (p - 2)(p - 3) u_{t} u_{jn}^2   \\
&  + 2 (p - 2)u_t u_{jn}^2 + 6 u_n u_{jn}[u_{tj}  - (p - 2)\frac{{u_t }}{{u_n }}u_{nj}]- 6u_{jn}^2 u_t  \\
&+2\sum_{i\in G}\sum_{\alpha=1}^{n}
\frac{1}{u_{ii}}L_{\alpha\a}[u_nu_{ij\alpha}-2u_{i\a}u_{jn}]^2 +O(\mathcal H_{\phi})\\
=&  ( 2p+2) u_n u_{jn}u_{tj} - (p^2+p)u_t u_{jn}^2   \\
&+2\sum_{i\in G}\sum_{\alpha=1}^{n}
\frac{1}{u_{ii}}L_{\alpha\a}[u_nu_{ij\alpha}-2u_{i\a}u_{jn}]^2 +O(\mathcal H_{\phi}).
\end{align*}
Hence
\begin{align*}
Q_j -2u_n u_{jn}u_{tj} =& 2p u_n u_{jn}u_{tj} - (p^2+p)u_t u_{jn}^2   \\
&+2\sum_{i\in G}\sum_{\alpha=1}^{n}
\frac{1}{u_{ii}}L_{\alpha\a}[u_nu_{ij\alpha}-2u_{i\a}u_{jn}]^2  +O(\mathcal H_{\phi}) \\
=& 2p  u_{jn}[ u_t u_{jn} - \frac{\hat h_{jn}}{u_t}] - (p^2+p)u_t u_{jn}^2   \\
&+2\sum_{i\in G}\sum_{\alpha=1}^{n}
\frac{1}{u_{ii}}L_{\alpha\a}[u_nu_{ij\alpha}-2u_{i\a}u_{jn}]^2 +O(\mathcal H_{\phi}) \\
=& \frac{p}{u_t} [\frac{1}{p-1} \frac{\hat h_{jn}^2}{u_t^2}-(p-1)(u_t u_{jn}+ \frac{1}{p-1} \frac{\hat h_{jn}}{u_t})^2]   \\
&+2\sum_{i\in G}\sum_{\alpha=1}^{n}
\frac{1}{u_{ii}}L_{\alpha\a}[u_nu_{ij\alpha}-2u_{i\a}u_{jn}]^2 +O(\mathcal H_{\phi})\\
\leq& \frac{p}{u_t} \cdot \frac{1}{p-1} \frac{\hat h_{jn}^2}{u_t^2} +O(\mathcal H_{\phi})=O(\mathcal H_{\phi}).
\end{align*}
So we can get
\begin{eqnarray}\label{3.42}
L_{\alpha\beta}\phi_{\alpha\beta}- \phi_t &\le&  C(\phi+\sum_{i,j\in B}|\nabla
a_{ij}|) \nonumber\\
&&-\frac{1}{{\sigma}^3_1(B)}\sum_{\alpha,\beta=1}^n\sum_{i\in
B} L_{\alpha\beta}[{\sigma}_1(B)a_{ii,\alpha}-a_{ii}\sum_{j\in
B}a_{jj,\alpha}][{\sigma}_1(B)a_{ii,\beta}-a_{ii}\sum_{j\in
B}a_{jj,\beta}]\nonumber\\
&&-\frac{1}{{\sigma}_1(B)}\sum_{\alpha,\beta=1}^n\sum_{i\neq j,
i,j\in B}L_{\alpha\beta}a_{ij,\alpha}a_{ij,\beta},
\end{eqnarray}
Following the proof of Theorem \ref{th3.1}, we get \eqref{3.37}.

\qed

\begin{remark}\label{rem3.4}
The constant rank properties ( that is, Corollary \ref{cor3.2}) still holds for the equation \eqref{1.12}.
\end{remark}

\subsection{Constant rank theorem of the spatial fundamental form for the equation \eqref{1.13}}

In this subsection, we consider the mean curvature parabolic equation, that is
\begin{align}\label{3.43}
 u_t  = \texttt{div}(\frac{{\nabla u}}{{\sqrt {1 + |\nabla u|^2 } }}) = m_{\alpha \beta} (\nabla u)u_{\alpha \beta}, \text{ in } \Omega \times (0, T],
 \end{align}
where
\begin{align}\label{3.44}
m_{\alpha \beta} (\nabla u) = (1 + |\nabla u|^2 )^{ - \frac{1}{2}} \delta _{\alpha \beta}  - (1 + |\nabla u|^2 )^{ - \frac{3}{2}} u_\alpha u_\beta.
\end{align}
We will establish the constant rank theorem
for the spatial second fundamental form $a$ as follows.
\begin{theorem}\label{th3.5}
Suppose $u \in C^{3,1}(\Omega \times (0,T])$ is a spacetime quasiconcave to the parabolic equation \eqref{3.43} and satisfies \eqref{1.4}.
Then the second fundamental form of spatial level sets $\Sigma^c = \{x \in \Omega | u(x,t) =
c\}$ has the same constant rank in $\Omega$ for any fixed $t \in (0, T]$. Moreover, let $l(t)$ be the minimal
rank of the second fundamental form in $\Omega$, then $l(s)
\leqslant l(t)$ for all $0< s \leqslant t \leqslant T$.
\end{theorem}

PROOF. The proof is similar to the the proof of Theorem \ref{th3.1} and Theorem \ref{th3.3}, with some modifications.

Suppose $a(x,t)$ attains minimal rank $l$ at some point $(x_0,t_0) \in \Omega \times (0, T]$. We may
assume $l\leqslant n-2$, otherwise there is nothing to prove. So there is a small neighborhood $\mathcal {O}\times
(t_0-\delta, t_0]$ of $(x_0, t_0)$, such that there are $l$
"good" eigenvalues of $(a_{ij})$ which are bounded below by a
positive constant, and the other $n-1-l$ "bad" eigenvalues of
$(a_{ij})$ are very small. Denote $G$ be the index set of these
"good" eigenvalues and $B$ be the index set of "bad" eigenvalues.
We will prove the differential inequality
\begin{equation} \label{3.45}
\sum_{\alpha,\beta=1}^{n}m_{\alpha\beta}\phi_{\alpha\beta}(x,t)-\phi_t\leq
C(\phi+|\nabla \phi|),~~\forall ~(x,t)\in \mathcal {O}
\times(t_0-\delta,t_0],
\end{equation}
where $\phi$ is defined in \eqref{3.4} and $C$ is a positive constant  independent of $\phi$.
Then by the strong maximum principle and the method of continuity, Theorem \ref{th3.5} holds.

For any fixed point  $(x,t) \in \mathcal {O}\times (t_0-\delta,
t_0]$, we may express $(a_{ij})$ in a form of
\eqref{2.3}, by choosing $e_1,\cdots, e_{n-1},e_n$ such that
\begin{equation}\label{3.46}
 |\n u(x,t)|=u_n(x,t)>0\ \mbox{and}\
\Big(u_{ij} \Big)_{1 \leq i,j \leq n-1} \mbox{is diagonal at}\ (x,t).
\end{equation}
Following the proof of Theorem \ref{th3.1}, we get from \eqref{3.21}
\begin{eqnarray}\label{3.47}
&&m_{\alpha\beta}\phi_{\alpha\beta}- \phi_t\nonumber\\
&=& u_n^{-3}\sum_{j\in
B}\left[\sigma_l(G)+\frac{{\sigma}^2_1(B|j)-{\sigma}_2(B|j)}{{\sigma}^2_1(B)}\right] \Big(  Q_j -2 u_n u_{jn} u_{jt} \Big) \nonumber\\
&&-\frac{1}{{\sigma}^3_1(B)}\sum_{\alpha,\beta=1}^n\sum_{i\in
B} m_{\alpha\beta}[{\sigma}_1(B)a_{ii,\alpha}-a_{ii}\sum_{j\in
B}a_{jj,\alpha}][{\sigma}_1(B)a_{ii,\beta}-a_{ii}\sum_{j\in
B}a_{jj,\beta}]\nonumber\\
&&-\frac{1}{{\sigma}_1(B)}\sum_{\alpha,\beta=1}^n\sum_{i\neq j,
i,j\in B}m_{\alpha\beta}a_{ij,\alpha}a_{ij,\beta}+O(\mathcal H_{\phi}),
\end{eqnarray}
where
\begin{align*}
Q_j=&  2\sum\limits_{kl = 1}^n \frac{{\partial m_{kl} }}{{\partial u_n }} u_{klj} u_{jn}u_n^2 +  \sum\limits_{kl = 1}^n \frac{{\partial ^2 m_{kl} }}{{\partial u_n ^2 }} u_{kl} u_{jn}^2 u_n^2    \\
&  + 2 \sum\limits_{kl = 1}^n \frac{{\partial m_{kl} }}{{\partial u_n }} u_{kl}u_n u_{jn}^2 +6u_nu_{jn}\sum_{kl=1}^n
m_{kl}u_{jkl} - 6u_{jn}^2\sum_{kl=1}^n m_{kl}u_{kl}  \\
&+2\sum_{i\in G}\sum_{\alpha,\beta=1}^{n}
\frac{1}{u_{ii}}m_{\alpha\beta}[u_nu_{ij\alpha}-2u_{i\a}u_{jn}][u_nu_{ij\b}-2u_{i\b}u_{jn}].
\end{align*}

Under the coordinate \eqref{3.46}, we get
\begin{align*}
m_{kl}  = 0, ~~k \ne l; \quad m_{kk}  = (1 + u_n ^2 )^{ - \frac{1}{2}},~~k < n; \quad m_{nn}  = (1 + u_n ^2 )^{ - \frac{3}{2}};
\end{align*}

\begin{align*}
&\frac{{\partial m_{kl} }}{{\partial u_n }} = 0, \quad k \ne l; \\
&\frac{{\partial m_{kk} }}{{\partial u_n }} =  - (1 + u_n ^2 )^{ - \frac{3}{2}} u_n, \quad k < n; \\
&\frac{{\partial m_{nn} }}{{\partial u_n }} =  - 3(1 + u_n ^2 )^{ - \frac{5}{2}} u_n;
\end{align*}
and
\begin{align*}
&\frac{{\partial ^2 m_{kl} }}{{\partial u_n ^2 }} = 0, \quad k \ne l; \\
&\frac{{\partial ^2 m_{kk} }}{{\partial u_n ^2 }} =  - (1 + u_n ^2 )^{ - \frac{3}{2}}  + 3(1 + u_n ^2 )^{ - \frac{5}{2}} u_n ^2, \quad k < n; \\
&\frac{{\partial ^2 m_{nn} }}{{\partial u_n ^2 }} =  - 3(1 + u_n ^2 )^{ - \frac{5}{2}}  + 15(1 + u_n ^2 )^{ - \frac{7}{2}} u_n ^2.
\end{align*}

From \eqref{3.9}, and \eqref{3.11}
\begin{equation}\label{3.48}
u_{kk} = O(\mathcal H_{\phi}), \quad u_{kkj} = O(\mathcal H_{\phi}), ~~\forall k \in B, j \in B,
\end{equation}
From \eqref{3.13}, we get
\begin{eqnarray}\label{3.49}
\hat h_{jn}^2= O(\mathcal H_{\phi}) ,\forall j \in B.
\end{eqnarray}

From the equation \eqref{3.43}, we know
\begin{align*}
 u_t  = \sum\limits_{k = 1}^n {m_{kk} u_{kk} }  = (1 + u_n ^2 )^{ - \frac{1}{2}} \sum\limits_{k = 1}^{n - 1} {u_{kk} }  + (1 + u_n ^2 )^{ - \frac{3}{2}} u_{nn},
\end{align*}
so we get
\begin{align*}
  u_t  - (1 + u_n ^2 )^{ - \frac{3}{2}} u_{nn}=&  (1 + u_n ^2 )^{ - \frac{1}{2}} \sum\limits_{k = 1}^{n - 1} {u_{kk} } \\
  =& (1 + u_n ^2 )^{ - \frac{1}{2}} \sum\limits_{k \in G} {u_{kk} }+ O(\mathcal H_{\phi}),
\end{align*}
and since $u_{kk} \leq 0 $ for $k <n$, it yields
\begin{align*}
 (1 + u_n ^2 )^{ - \frac{3}{2}} u_{nn}  \ge u_t.
\end{align*}
Hence we can get
\begin{align*}
\sum\limits_{k l= 1}^n {\frac{{\partial m_{kl} }}{{\partial u_n }} u_{kl} }  =&  \sum\limits_{k= 1}^n {\frac{{\partial m_{kk} }}{{\partial u_n }} u_{kk} }    = \sum\limits_{k= 1}^{n-1} {\frac{{\partial m_{kk} }}{{\partial u_n }} u_{kk} }+ \frac{{\partial m_{nn} }}{{\partial u_n }} u_{nn} \\
=& - (1 + u_n ^2 )^{ - \frac{3}{2}} u_n\sum\limits_{k= 1}^{n-1} { u_{kk} } - 3(1 + u_n ^2 )^{ - \frac{5}{2}} u_n  u_{nn} \\
=& - (1 + u_n ^2 )^{ - 1} u_n [ u_t  - (1 + u_n ^2 )^{ - \frac{3}{2}} u_{nn} ]- 3(1 + u_n ^2 )^{ - \frac{5}{2}} u_n  u_{nn}\\
=& - (1 + u_n ^2 )^{ - 1} u_t u_n - 2(1 + u_n ^2 )^{ - \frac{5}{2}} u_n  u_{nn}.
\end{align*}
and
\begin{align*}
\sum\limits_{kl = 1}^n \frac{{\partial ^2 m_{kl} }}{{\partial u_n ^2 }} u_{kl}=&\sum\limits_{k = 1}^n \frac{{\partial ^2 m_{kk} }}{{\partial u_n ^2 }} u_{kk}=\sum\limits_{k = 1}^{n-1} \frac{{\partial ^2 m_{kk} }}{{\partial u_n ^2 }} u_{kk} + \frac{{\partial ^2 m_{nn} }}{{\partial u_n ^2 }} u_{nn} \\
=& \Big[ - (1 + u_n ^2 )^{ - \frac{3}{2}}  + 3(1 + u_n ^2 )^{ - \frac{5}{2}} u_n ^2 \Big] \sum\limits_{k = 1}^{n-1} u_{kk} + \Big[- 3(1 + u_n ^2 )^{ - \frac{5}{2}}  + 15(1 + u_n ^2 )^{ - \frac{7}{2}} u_n ^2\Big] u_{nn} \\
=& \Big[ - (1 + u_n ^2 )^{ - 1}  + 3(1 + u_n ^2 )^{ - 2} u_n ^2 \Big][ u_t  - (1 + u_n ^2 )^{ - \frac{3}{2}} u_{nn}  ]+ \Big[- 3(1 + u_n ^2 )^{ - \frac{5}{2}}  + 15(1 + u_n ^2 )^{ - \frac{7}{2}} u_n ^2\Big] u_{nn} \\
=&  - (1 + u_n ^2 )^{ - 1}u_t  + 3(1 + u_n ^2 )^{ - 2} u_n ^2 u_t + \Big[- 2(1 + u_n ^2 )^{ - \frac{5}{2}}  + 12(1 + u_n ^2 )^{ - \frac{7}{2}} u_n ^2\Big] u_{nn} \\
\end{align*}
For $j \in B$, differentiating the equation once in $x_j$, we get
\begin{align*}
u_{tj}  =& \sum\limits_{k = 1}^n {m_{kk} u_{kkj} }  + \sum\limits_{k = 1}^n {\frac{{\partial m_{kk} }}{{\partial u_n }}u_{nj} u_{kk} }+ O(\mathcal H_{\phi}),
\end{align*}
so
\begin{align*}
\sum\limits_{kl = 1}^n {m_{kl} u_{klj} }   =&  \sum\limits_{k = 1}^n {m_{kk} u_{kkj} }  = u_{tj} - \sum\limits_{k = 1}^n {\frac{{\partial m_{kl} }}{{\partial u_n }}u_{nj} u_{kl} } -\sum\limits_{k = 1}^n {\frac{{\partial m_{kl} }}{{\partial u_j }}u_{jj} u_{kl} } \\
=& u_{tj} - \sum\limits_{k = 1}^n {\frac{{\partial m_{kk} }}{{\partial u_n }}u_{nj} u_{kk} } + O(\mathcal H_{\phi}) \\
=& u_{tj}+ (1 + u_n ^2 )^{ - \frac{3}{2}} u_n u_{nj} \sum\limits_{k = 1}^{n - 1} {u_{kk} }  + 3(1 + u_n ^2 )^{ - \frac{5}{2}} u_n u_{nj} u_{nn}+ O(\mathcal H_{\phi})\\
=& u_{tj}+ (1 + u_n ^2 )^{ - 1} u_n u_{nj} [u_t  - (1 + u_n ^2 )^{ - \frac{3}{2}} u_{nn} ]  + 3(1 + u_n ^2 )^{ - \frac{5}{2}} u_n u_{nj} u_{nn}+ O(\mathcal H_{\phi})\\
=& u_{tj}+ (1 + u_n ^2 )^{ - 1} u_t u_n u_{nj}   + 2(1 + u_n ^2 )^{ - \frac{5}{2}} u_n u_{nj} u_{nn}+ O(\mathcal H_{\phi}),
\end{align*}
and from \eqref{3.49},
\begin{align*}
 (1 + u_n ^2 )^{ - \frac{3}{2}} u_{nnj}  =& u_{tj}  - (1 + u_n ^2 )^{ - \frac{1}{2}} \sum\limits_{k = 1}^{n - 1} {u_{kkj} }  + (1 + u_n ^2 )^{ - 1} u_t u_n u_{nj}  + 2(1 + u_n ^2 )^{ - \frac{5}{2}} u_n u_{nj} u_{nn}+ O(\mathcal H_{\phi}) \\
 =& u_{tj}  - (1 + u_n ^2 )^{ - \frac{1}{2}} \sum\limits_{k \in G}{u_{kkj} }  + (1 + u_n ^2 )^{ - 1} u_t u_n u_{nj}  + 2(1 + u_n ^2 )^{ - \frac{5}{2}} u_n u_{nj} u_{nn}+ O(\mathcal H_{\phi}).
\end{align*}
Hence
\begin{align*}
\sum\limits_{k l= 1}^n {\frac{{\partial m_{kl} }}{{\partial u_n }} u_{klj} }  =&  \sum\limits_{k= 1}^n {\frac{{\partial m_{kk} }}{{\partial u_n }} u_{kkj} }  = \sum\limits_{k= 1}^{n-1} {\frac{{\partial m_{kk} }}{{\partial u_n }} u_{kkj} }+ \frac{{\partial m_{nn} }}{{\partial u_n }} u_{nnj} \\
=& - (1 + u_n ^2 )^{ - \frac{3}{2}} u_n\sum\limits_{k= 1}^{n-1} { u_{kkj} } - 3(1 + u_n ^2 )^{ - \frac{5}{2}} u_n  u_{nnj} \\
=&  - (1 + u_n ^2 )^{ - \frac{3}{2}} u_n\sum\limits_{k\in G}  { u_{kkj} } + O(\mathcal H_{\phi}) \\
& - 3(1 + u_n ^2 )^{ - 1} u_n \Big[ u_{tj}  - (1 + u_n ^2 )^{ - \frac{1}{2}} \sum\limits_{k = 1}^{n - 1} {u_{kkj} }  + (1 + u_n ^2 )^{ - 1} u_t u_n u_{nj}  + 2(1 + u_n ^2 )^{ - \frac{5}{2}} u_n u_{nj} u_{nn} \Big]\\
=& - 3(1 + u_n ^2 )^{ - 1} u_n u_{tj}  +2 (1 + u_n ^2 )^{ - \frac{3}{2}} u_n\sum\limits_{k \in G}{u_{kkj} }  -3 (1 + u_n ^2 )^{ - 2} u_t u_n^2 u_{nj} \\
& -6 (1 + u_n ^2 )^{ - \frac{7}{2}} u_n ^2 u_{nj} u_{nn} + O(\mathcal H_{\phi}).
\end{align*}

So
\begin{align*}
Q_j=&   4(1 + u_n ^2 )^{ - \frac{3}{2}} u_n\sum\limits_{k \in G} {u_{kkj} } \cdot u_{jn}u_n^2 +2(1 + u_n ^2 )^{ - \frac{1}{2}}\sum_{k \in G}\frac{1}{u_{kk}}[u_nu_{kkj}-2u_{kk}u_{jn}]^2  \\
&+2\Big[- 3(1 + u_n ^2 )^{ - 1} u_n u_{tj}  -3 (1 + u_n ^2 )^{ - 2} u_t u_n^2 u_{nj}  -6 (1 + u_n ^2 )^{ - \frac{7}{2}} u_n ^2 u_{nj} u_{nn} \Big] u_{jn}u_n^2  \\
&+ \Big[ - (1 + u_n ^2 )^{ - 1}u_t  + 3(1 + u_n ^2 )^{ - 2} u_n ^2 u_t + \Big(- 2(1 + u_n ^2 )^{ - \frac{5}{2}}  + 12(1 + u_n ^2 )^{ - \frac{7}{2}} u_n ^2\Big) u_{nn}\Big]u_{jn}^2 u_n^2    \\
&  + 2 \Big[ - (1 + u_n ^2 )^{ - 1} u_t u_n - 2(1 + u_n ^2 )^{ - \frac{5}{2}} u_n  u_{nn} \Big]u_n u_{jn}^2  \\
&+6u_nu_{jn} \Big[ u_{tj}+ (1 + u_n ^2 )^{ - 1} u_t u_n u_{nj}   + 2(1 + u_n ^2 )^{ - \frac{5}{2}} u_n u_{nj} u_{nn} \Big]- 6u_{jn}^2u_{t}  \\
& +2\sum_{i\in G}\sum_{\alpha=1, \a \ne i}^{n}
\frac{1}{u_{ii}}m_{\alpha\a}[u_nu_{ij\alpha}-2u_{i\a}u_{jn}]^2 + O(\mathcal H_{\phi})\\
=&4(1 + u_n ^2 )^{ - \frac{3}{2}} u_n\sum\limits_{k \in G} {u_{kkj} } \cdot u_{jn}u_n^2 +2(1 + u_n ^2 )^{ - \frac{1}{2}}\sum_{k \in G}\frac{1}{u_{kk}}[u_nu_{kkj}-2u_{kk}u_{jn}]^2  \\
&- 6(1 + u_n ^2 )^{ - 1} u_n^2 \cdot u_n u_{tj} u_{jn}  -3 (1 + u_n ^2 )^{ - 2}  u_n^4 \cdot u_t u_{nj}^2 +3 (1 + u_n ^2 )^{ - 1}  u_n^2 \cdot u_t u_{nj}^2 \\
&+ 6(1 + u_n ^2 )^{ - \frac{5}{2}} u_n^2 \cdot u_{nn} u_{jn}^2 + 6u_nu_{jn}u_{tj} - 6u_{jn}^2u_{t}  \\
& +2\sum_{i\in G}\sum_{\alpha=1, \a \ne i}^{n}
\frac{1}{u_{ii}}m_{\alpha\a}[u_nu_{ij\alpha}-2u_{i\a}u_{jn}]^2 + O(\mathcal H_{\phi})
\end{align*}
where
\begin{align*}
&4(1 + u_n ^2 )^{ - \frac{3}{2}} u_n\sum\limits_{k \in G} {u_{kkj} } \cdot u_{jn}u_n^2 +2(1 + u_n ^2 )^{ - \frac{1}{2}}\sum_{k \in G}\frac{1}{u_{kk}}[u_nu_{kkj}-2u_{kk}u_{jn}]^2  \\
=&2(1 + u_n ^2 )^{ - \frac{1}{2}}\Big( \sum_{k \in G}\frac{1}{u_{kk}}[u_nu_{kkj}-2u_{kk}u_{jn}]^2  + 2(1 + u_n ^2 )^{ - 1}u_n^2 u_{jn} \sum\limits_{k \in G} [u_nu_{kkj}-2u_{kk}u_{jn} + 2u_{kk}u_{jn}] \Big)  \\
=&2(1 + u_n ^2 )^{ - \frac{1}{2}}\Big( \sum_{k \in G}\frac{1}{u_{kk}}[u_nu_{kkj}-2u_{kk}u_{jn}  +(1 + u_n ^2 )^{ - 1}u_n^2 u_{jn} u_{kk} ]^2  - \sum_{k \in G}\frac{1}{u_{kk}}[(1 + u_n ^2 )^{ - 1}u_n^2 u_{jn} u_{kk} ]^2 \Big)  \\
&+ 8(1 + u_n ^2 )^{ - \frac{3}{2}}u_n^2 u_{jn}^2 \sum\limits_{k \in G} u_{kk}  \\
=&2(1 + u_n ^2 )^{ - \frac{1}{2}}\sum_{k \in G}\frac{1}{u_{kk}}[u_nu_{kkj}-2u_{kk}u_{jn}  +(1 + u_n ^2 )^{ - 1}u_n^2 u_{jn} u_{kk} ]^2    \\
&- 2(1 + u_n ^2 )^{ - \frac{5}{2}} u_n^4 u_{jn}^2  \sum_{k \in G} u_{kk}+ 8(1 + u_n ^2 )^{ - \frac{3}{2}}u_n^2 u_{jn}^2 \sum\limits_{k \in G} u_{kk}\\
=&2(1 + u_n ^2 )^{ - \frac{1}{2}}\sum_{k \in G}\frac{1}{u_{kk}}[u_nu_{kkj}-2u_{kk}u_{jn}  +(1 + u_n ^2 )^{ - 1}u_n^2 u_{jn} u_{kk} ]^2    \\
&- 2(1 + u_n ^2 )^{ - 2} u_n^4 u_{jn}^2 [u_t  - (1 + u_n ^2 )^{ - \frac{3}{2}} u_{nn}]+ 8(1 + u_n ^2 )^{ - 1}u_n^2 u_{jn}^2 [u_t  - (1 + u_n ^2 )^{ - \frac{3}{2}} u_{nn}] + O(\mathcal H_{\phi}) \\
=&2(1 + u_n ^2 )^{ - \frac{1}{2}}\sum_{k \in G}\frac{1}{u_{kk}}[u_nu_{kkj}-2u_{kk}u_{jn}  +(1 + u_n ^2 )^{ - 1}u_n^2 u_{jn} u_{kk} ]^2 + O(\mathcal H_{\phi})   \\
&- 2(1 + u_n ^2 )^{ - 2} u_n^4 \cdot u_tu_{jn}^2 + 8(1 + u_n ^2 )^{ - 1}u_n^2 \cdot u_tu_{jn}^2 +[2 (1 + u_n ^2 )^{ - \frac{7}{2}} u_n^4- 8(1 + u_n ^2 )^{ - \frac{5}{2}}u_n^2  ]u_{nn}u_{jn}^2
\end{align*}
So we can get
\begin{align*}
Q_j -2u_n u_{jn}u_{tj} =& 2(1 + u_n ^2 )^{ - \frac{1}{2}}\sum_{k \in G}\frac{1}{u_{kk}}[u_nu_{kkj}-2u_{kk}u_{jn}  +(1 + u_n ^2 )^{ - 1}u_n^2 u_{jn} u_{kk} ]^2    \\
&- 6(1 + u_n ^2 )^{ - 1} u_n^2 \cdot u_n u_{tj} u_{jn}  -5 (1 + u_n ^2 )^{ - 2}  u_n^4 \cdot u_t u_{nj}^2 + 11 (1 + u_n ^2 )^{ - 1}  u_n^2 \cdot u_t u_{nj}^2 \\
&+[2 (1 + u_n ^2 )^{ - \frac{7}{2}} u_n^4- 2(1 + u_n ^2 )^{ - \frac{5}{2}}u_n^2  ]u_{nn}u_{jn}^2 + 4u_nu_{jn}u_{tj} - 6u_{jn}^2u_{t}  \\
& +2\sum_{i\in G}\sum_{\alpha=1, \a \ne i}^{n}
\frac{1}{u_{ii}}m_{\alpha\a}[u_nu_{ij\alpha}-2u_{i\a}u_{jn}]^2 + O(\mathcal H_{\phi}) \\
\end{align*}
Hence
\begin{align*}
Q_j -2u_n u_{jn}u_{tj} \leq&- 6(1 + u_n ^2 )^{ - 1} u_n^2 \cdot u_n u_{tj} u_{jn}  -5 (1 + u_n ^2 )^{ - 2}  u_n^4 \cdot u_t u_{nj}^2 + 11 (1 + u_n ^2 )^{ - 1}  u_n^2 \cdot u_t u_{nj}^2 \\
&+[2 (1 + u_n ^2 )^{ -2} u_n^4- 2(1 + u_n ^2 )^{ - 1}u_n^2  ] \cdot u_{t}u_{jn}^2 + 4u_nu_{jn}u_{tj} - 6u_{jn}^2u_{t} + O(\mathcal H_{\phi}) \\
=& - 6(1 + u_n ^2 )^{ - 1} u_n^2 \cdot u_{jn} [u_t u_{jn} - \frac{\hat h_{jn}}{u_t}]  -3 (1 + u_n ^2 )^{ - 2}  u_n^4 \cdot u_t u_{nj}^2  \\
&+ 9 (1 + u_n ^2 )^{ - 1}  u_n^2 \cdot u_t u_{nj}^2 + 4u_{jn} [u_t u_{jn} - \frac{\hat h_{jn}}{u_t}] - 6u_{jn}^2u_{t} + O(\mathcal H_{\phi}) \\
=& \frac{\hat h_{jn}}{u_t}  [6(1 + u_n ^2 )^{ - 1} u_n^2 -4 ] \cdot u_{jn}   \\
&+[ - 2+ 3(1 + u_n ^2 )^{ - 1}  u_n^2 -3 (1 + u_n ^2 )^{ - 2}  u_n^4 ]\cdot u_t u_{nj}^2 + O(\mathcal H_{\phi})\\
=& \frac{\hat h_{jn}}{u_t}  [6(1 + u_n ^2 )^{ - 1} u_n^2 -4 ] \cdot u_{jn} -  \frac{5}{4}u_t u_{nj}^2 \\
&-3 [ (1 + u_n ^2 )^{ - 1}  u_n^2 - \frac{1}{2}]^2\cdot u_t u_{nj}^2 + O(\mathcal H_{\phi}) \\
=& -u_t \Big(\frac{\hat h_{jn}}{u_t^2}  [3(1 + u_n ^2 )^{ - 1} u_n^2 -2 ] - u_{jn} \Big)^2 + u_t  [3(1 + u_n ^2 )^{ - 1} u_n^2 -2 ]^2\frac{\hat h_{jn}^2}{u_t^4}   -  \frac{1}{4}u_t u_{nj}^2 \\
&-3 [ (1 + u_n ^2 )^{ - 1}  u_n^2 - \frac{1}{2}]^2\cdot u_t u_{nj}^2 + O(\mathcal H_{\phi}) \\
\leq& u_t  [3(1 + u_n ^2 )^{ - 1} u_n^2 -2 ]^2\frac{\hat h_{jn}^2}{u_t^4} + O(\mathcal H_{\phi}) = O(\mathcal H_{\phi})
\end{align*}

So we can get
\begin{eqnarray}\label{3.50}
m_{\alpha\beta}\phi_{\alpha\beta}- \phi_t &\le&  C(\phi+\sum_{i,j\in B}|\nabla
a_{ij}|) \nonumber\\
&&-\frac{1}{{\sigma}^3_1(B)}\sum_{\alpha,\beta=1}^n\sum_{i\in
B} m_{\alpha\beta}[{\sigma}_1(B)a_{ii,\alpha}-a_{ii}\sum_{j\in
B}a_{jj,\alpha}][{\sigma}_1(B)a_{ii,\beta}-a_{ii}\sum_{j\in
B}a_{jj,\beta}]\nonumber\\
&&-\frac{1}{{\sigma}_1(B)}\sum_{\alpha,\beta=1}^n\sum_{i\neq j,
i,j\in B}m_{\alpha\beta}a_{ij,\alpha}a_{ij,\beta},
\end{eqnarray}
Following the proof of Theorem \ref{th3.1}, we get \eqref{3.45}.
\qed

\begin{remark}\label{rem3.6}
The constant rank properties ( that is, Corollary \ref{cor3.2}) still holds for the equation \eqref{1.13}.
\end{remark}

\subsection{Constant rank theorem of the spatial fundamental form for the equation \eqref{1.14}}

In this subsection, we consider the mean curvature parabolic equation, that is
\begin{align}\label{3.51}
 u_t  = (1 + |\nabla u|^2)^{\frac{3}{2}}\texttt{div}(\frac{{\nabla u}}{{\sqrt {1 + |\nabla u|^2 } }}) = M_{\alpha \beta} (\nabla u)u_{\alpha \beta}, \text{ in } \Omega \times (0, T],
 \end{align}
where
\begin{align}\label{3.52}
M_{\alpha \beta} (\nabla u) = (1 + |\nabla u|^2 ) \delta _{\alpha \beta}  -  u_\alpha  u_\beta.
\end{align}
We will establish the constant rank theorem
for the spatial second fundamental form $a$ as follows.
\begin{theorem}\label{th3.7}
Suppose $u \in C^{3,1}(\Omega \times (0,T])$ is a spacetime quasiconcave to the parabolic equation \eqref{3.51} and satisfies \eqref{1.4}.
Then the second fundamental form of spatial level sets $\Sigma^c = \{x \in \Omega | u(x,t) =
c\}$ has the same constant rank in $\Omega$ for any fixed $ t \in (0, T]$. Moreover, let $l(t)$ be the minimal
rank of the second fundamental form in $\Omega$, then $l(s)
\leqslant l(t)$ for all $0 < s \leqslant t \leqslant T$.
\end{theorem}

PROOF. The proof is similar to the the proof of Theorem \ref{th3.1}, with some modifications.

Suppose $a(x,t)$ attains minimal rank $l$ at some point $(x_0,t_0) \in \Omega \times (0, T]$. We may
assume $l\leqslant n-2$, otherwise there is nothing to prove. So there is a small neighborhood $\mathcal {O}\times
(t_0-\delta, t_0]$ of $(x_0, t_0)$, such that there are $l$
"good" eigenvalues of $(a_{ij})$ which are bounded below by a
positive constant, and the other $n-1-l$ "bad" eigenvalues of
$(a_{ij})$ are very small. Denote $G$ be the index set of these
"good" eigenvalues and $B$ be the index set of "bad" eigenvalues.
We will prove the differential inequality
\begin{equation} \label{3.53}
\sum_{\alpha,\beta=1}^{n}M_{\alpha\beta}\phi_{\alpha\beta}(x,t)-\phi_t\leq
C(\phi+|\nabla \phi|),~~\forall ~(x,t)\in \mathcal {O}
\times(t_0-\delta,t_0],
\end{equation}
where $\phi$ is defined in \eqref{3.4} and $C$ is a positive constant  independent of $\phi$.
Then by the strong maximum principle and the method of continuity, Theorem \ref{th3.7} holds.

For any fixed point  $(x,t) \in \mathcal {O}\times (t_0-\delta,
t_0]$, we may express $(a_{ij})$ in a form of
\eqref{2.3}, by choosing $e_1,\cdots, e_{n-1},e_n$ such that
\begin{equation}\label{3.54}
 |\n u(x,t)|=u_n(x,t)>0\ \mbox{and}\
\Big(u_{ij} \Big)_{1 \leq i,j \leq n-1} \mbox{is diagonal at}\ (x,t).
\end{equation}
Following the proof of Theorem \ref{th3.1}, we get from \eqref{3.21}
\begin{eqnarray}\label{3.55}
&&M_{\alpha\beta}\phi_{\alpha\beta}- \phi_t\nonumber\\
&=& u_n^{-3}\sum_{j\in
B}\left[\sigma_l(G)+\frac{{\sigma}^2_1(B|j)-{\sigma}_2(B|j)}{{\sigma}^2_1(B)}\right] \Big(  Q_j -2 u_n u_{jn} u_{jt} \Big) \nonumber\\
&&-\frac{1}{{\sigma}^3_1(B)}\sum_{\alpha,\beta=1}^n\sum_{i\in
B} M_{\alpha\beta}[{\sigma}_1(B)a_{ii,\alpha}-a_{ii}\sum_{j\in
B}a_{jj,\alpha}][{\sigma}_1(B)a_{ii,\beta}-a_{ii}\sum_{j\in
B}a_{jj,\beta}]\nonumber\\
&&-\frac{1}{{\sigma}_1(B)}\sum_{\alpha,\beta=1}^n\sum_{i\neq j,
i,j\in B}M_{\alpha\beta}a_{ij,\alpha}a_{ij,\beta}+O(\mathcal H_{\phi}),
\end{eqnarray}
where
\begin{align*}
Q_j=&  2\sum\limits_{kl = 1}^n \frac{{\partial M_{kl} }}{{\partial u_n }} u_{klj} u_{jn}u_n^2 +  \sum\limits_{kl = 1}^n \frac{{\partial ^2 M_{kl} }}{{\partial u_n ^2 }} u_{kl} u_{jn}^2 u_n^2    \\
&  + 2 \sum\limits_{kl = 1}^n \frac{{\partial M_{kl} }}{{\partial u_n }} u_{kl}u_n u_{jn}^2 +6u_nu_{jn}\sum_{kl=1}^n
M_{kl}u_{jkl} - 6u_{jn}^2\sum_{kl=1}^n M_{kl}u_{kl}  \\
&+2\sum_{i\in G}\sum_{\alpha,\beta=1}^{n}
\frac{1}{u_{ii}}M_{\alpha\beta}[u_nu_{ij\alpha}-2u_{i\a}u_{jn}][u_nu_{ij\b}-2u_{i\b}u_{jn}].
\end{align*}

Under the coordinate \eqref{3.54}, we get
\begin{align*}
M_{kl}  = 0, k \ne l; \quad M_{kk}  = 1 + u_n ^2 , \quad k < n; \quad M_{nn}  = 1;
\end{align*}

\begin{align*}
\frac{{\partial M_{kl} }}{{\partial u_n }} = 0,~~ k \ne l;  \quad \frac{{\partial M_{kk} }}{{\partial u_n }} =  2 u_n ,~~k < n; \quad \frac{{\partial M_{nn} }}{{\partial u_n }} = 0;
\end{align*}
and
\begin{align*}
\frac{{\partial ^2 M_{kl} }}{{\partial u_n ^2 }} = 0,~~k \ne l; \quad \frac{{\partial ^2 M_{kk} }}{{\partial u_n ^2 }} =  2,~~k < n; \quad \frac{{\partial ^2 M_{nn} }}{{\partial u_n ^2 }} = 0.
\end{align*}

From \eqref{3.9}, \eqref{3.11}
\begin{equation}\label{3.56}
u_{kk} = O(\mathcal H_{\phi}), \quad u_{kkj} = O(\mathcal H_{\phi}), ~~\forall k \in B, j \in B.
\end{equation}
From \eqref{3.13}, we get
\begin{eqnarray}\label{3.57}
\hat h_{jn}^2= O(\mathcal H_{\phi}) ,\forall j \in B.
\end{eqnarray}

From the equation \eqref{3.51}, we know
\begin{align*}
 u_t  = \sum\limits_{k = 1}^n {M_{kk} u_{kk} }  = (1 + u_n ^2 ) \sum\limits_{k = 1}^{n - 1} {u_{kk} }  + u_{nn},
\end{align*}
so we get
\begin{align*}
 (1 + u_n ^2 )\sum\limits_{k = 1}^{n - 1} {u_{kk} }=  u_t  - u_{nn} ,
\end{align*}
and by $u_{kk} \leq 0$ for $k <n$, it yields
\begin{align*}
u_{nn}  \ge u_t.
\end{align*}
Hence we can get
\begin{align*}
\sum\limits_{k l= 1}^n {\frac{{\partial M_{kl} }}{{\partial u_n }} u_{kl} }  =&  \sum\limits_{k= 1}^n {\frac{{\partial M_{kk} }}{{\partial u_n }} u_{kk} }    = \sum\limits_{k= 1}^{n-1} {\frac{{\partial M_{kk} }}{{\partial u_n }} u_{kk} }+ \frac{{\partial M_{nn} }}{{\partial u_n }} u_{nn} \\
=& 2 u_n\sum\limits_{k= 1}^{n-1} { u_{kk} }  \\
=& 2u_n  (1 + u_n ^2 )^{ - 1} [ u_t  -  u_{nn} ]\\
=& 2(1 + u_n ^2 )^{ - 1} u_t u_n - 2(1 + u_n ^2 )^{ - 1} u_n  u_{nn}.
\end{align*}
and
\begin{align*}
\sum\limits_{kl = 1}^n \frac{{\partial ^2 M_{kl} }}{{\partial u_n ^2 }} u_{kl}=&\sum\limits_{k = 1}^n \frac{{\partial ^2 M_{kk} }}{{\partial u_n ^2 }} u_{kk}=\sum\limits_{k = 1}^{n-1} \frac{{\partial ^2 M_{kk} }}{{\partial u_n ^2 }} u_{kk} + \frac{{\partial ^2 M_{nn} }}{{\partial u_n ^2 }} u_{nn} \\
=& 2 \sum\limits_{k = 1}^{n-1} u_{kk}  = 2 (1 + u_n ^2 )^{ - 1}  [ u_t  -  u_{nn}  ] \\
=&  2(1 + u_n ^2 )^{ - 1}u_t  -2(1 + u_n ^2 )^{ - 1} u_{nn}
\end{align*}

For $j \in B$, differentiating the equation once in $x_j$, we get
\begin{align*}
u_{tj} =& \sum\limits_{k = 1}^n {M_{kk} u_{kkj} }  + \sum\limits_{k l= 1}^n {\frac{{\partial M_{kl} }}{{\partial u_p }}u_{pj} u_{kl} }  \\
=& \sum\limits_{k = 1}^n {M_{kk} u_{kkj} }  + \sum\limits_{k l= 1}^n {\frac{{\partial M_{kl} }}{{\partial u_j}}u_{jj} u_{kl} }  + \sum\limits_{k l= 1}^n {\frac{{\partial M_{kl} }}{{\partial u_n }}u_{nj} u_{kl} }\\
=& \sum\limits_{k = 1}^n {M_{kk} u_{kkj} }  + \sum\limits_{k = 1}^n {\frac{{\partial M_{kk} }}{{\partial u_n }}u_{nj} u_{kk} } + O(\mathcal H_{\phi}),
\end{align*}
so
\begin{align*}
\sum\limits_{kl = 1}^n {M_{kl} u_{klj} }   =&  \sum\limits_{k = 1}^n {M_{kk} u_{kkj} }   \\
=& u_{tj} - \sum\limits_{k = 1}^n {\frac{{\partial M_{kk} }}{{\partial u_n }}u_{nj} u_{kk} } + O(\mathcal H_{\phi})\\
=& u_{tj}-2 u_n u_{nj} \sum\limits_{k = 1}^{n - 1} {u_{kk} }  + O(\mathcal H_{\phi})\\
=& u_{tj}-2 u_n u_{nj} (1 + u_n ^2 )^{ - 1} [u_t  -  u_{nn} ]  + O(\mathcal H_{\phi})\\
=& u_{tj}-2 (1 + u_n ^2 )^{ - 1} u_t u_n u_{nj}   + 2(1 + u_n ^2 )^{ - 1} u_n u_{nj} u_{nn}+ O(\mathcal H_{\phi})
\end{align*}
Hence from \eqref{3.56},
\begin{align*}
\sum\limits_{k l= 1}^n {\frac{{\partial M_{kl} }}{{\partial u_n }} u_{klj} }  =&  \sum\limits_{k= 1}^n {\frac{{\partial M_{kk} }}{{\partial u_n }} u_{kkj} }  = \sum\limits_{k= 1}^{n-1} {\frac{{\partial M_{kk} }}{{\partial u_n }} u_{kkj} }+ \frac{{\partial M_{nn} }}{{\partial u_n }} u_{nnj} \\
=& 2 u_n\sum\limits_{k= 1}^{n-1} { u_{kkj} }   \\
=& 2 u_n \sum\limits_{k \in G} { u_{kkj} }+ O(\mathcal H_{\phi})
\end{align*}

So
\begin{align*}
Q_j=&   4 u_n\sum\limits_{k \in G} {u_{kkj} } \cdot u_{jn}u_n^2  +2(1 + u_n ^2 )^{ - 1} [u_t - u_{nn} ] u_n^2 u_{jn}^2  + 4 u_n (1 + u_n ^2 )^{ - 1} [u_t - u_{nn} ] u_n u_{jn}^2  \\
&+6u_nu_{jn} \Big[ u_{tj}-2 (1 + u_n ^2 )^{ - 1} u_t u_n u_{nj}   + 2(1 + u_n ^2 )^{ - 1} u_n u_{nj} u_{nn} \Big]- 6u_{jn}^2u_{t}  \\
& +2(1 + u_n ^2 )\sum_{k \in G}\frac{1}{u_{kk}}[u_nu_{kkj}-2u_{kk}u_{jn}]^2 +2\sum_{i\in G}\sum_{\alpha=1, \a \ne i}^{n}
\frac{1}{u_{ii}}M_{\alpha\a}[u_nu_{ij\alpha}-2u_{i\a}u_{jn}]^2 + O(\mathcal H_{\phi})\\
=&   4 u_n\sum\limits_{k \in G} {u_{kkj} } \cdot u_{jn}u_n^2   +2(1 + u_n ^2 )\sum_{k \in G}\frac{1}{u_{kk}}[u_nu_{kkj}-2u_{kk}u_{jn}]^2   \\
& -6 (1 + u_n ^2 )^{ - 1}u_n^2 \cdot   u_{nj}^2[ u_t - u_{nn}] + 6u_nu_{jn} u_{tj} - 6u_{jn}^2u_{t}  \\
&+2\sum_{i\in G}\sum_{\alpha=1, \a \ne i}^{n} \frac{1}{u_{ii}}M_{\alpha\a}[u_nu_{ij\alpha}-2u_{i\a}u_{jn}]^2 + O(\mathcal H_{\phi})\\
\end{align*}
where
\begin{align*}
& 4 u_n\sum\limits_{k \in G} {u_{kkj} } \cdot u_{jn}u_n^2   +2(1 + u_n ^2 )\sum_{k \in G}\frac{1}{u_{kk}}[u_nu_{kkj}-2u_{kk}u_{jn}]^2   \\
=&2(1 + u_n ^2 )\Big( \sum_{k \in G}\frac{1}{u_{kk}}[u_nu_{kkj}-2u_{kk}u_{jn}]^2  + 2(1 + u_n ^2 )^{ - 1}u_n^2 u_{jn} \sum\limits_{k \in G} [u_nu_{kkj}-2u_{kk}u_{jn} + 2u_{kk}u_{jn}] \Big)  \\
=&2(1 + u_n ^2 )\Big( \sum_{k \in G}\frac{1}{u_{kk}}[u_nu_{kkj}-2u_{kk}u_{jn}  +(1 + u_n ^2 )^{ - 1}u_n^2 u_{jn} u_{kk} ]^2  - \sum_{k \in G}\frac{1}{u_{kk}}[(1 + u_n ^2 )^{ - 1}u_n^2 u_{jn} u_{kk} ]^2 \Big)  \\
&+ 8u_n^2 u_{jn}^2 \sum\limits_{k \in G} u_{kk}  \\
=&2(1 + u_n ^2 )\sum_{k \in G}\frac{1}{u_{kk}}[u_nu_{kkj}-2u_{kk}u_{jn}  +(1 + u_n ^2 )^{ - 1}u_n^2 u_{jn} u_{kk} ]^2    \\
&- 2(1 + u_n ^2 )^{ - 1} u_n^4 u_{jn}^2  \sum_{k \in G} u_{kk}+ 8u_n^2 u_{jn}^2 \sum\limits_{k \in G} u_{kk}\\
=&2(1 + u_n ^2 )\sum_{k \in G}\frac{1}{u_{kk}}[u_nu_{kkj}-2u_{kk}u_{jn}  +(1 + u_n ^2 )^{ - 1}u_n^2 u_{jn} u_{kk} ]^2    \\
&- 2(1 + u_n ^2 )^{ - 2} u_n^4 \cdot u_{jn}^2 [u_t  -  u_{nn}]+ 8(1 + u_n ^2 )^{ - 1}u_n^2 \cdot u_{jn}^2 [u_t  -  u_{nn}] + O(\mathcal H_{\phi})
\end{align*}
So we can get
\begin{align*}
Q_j -2u_n u_{jn}u_{tj} =& 2(1 + u_n ^2 )\sum_{k \in G}\frac{1}{u_{kk}}[u_nu_{kkj}-2u_{kk}u_{jn}  +(1 + u_n ^2 )^{ - 1}u_n^2 u_{jn} u_{kk} ]^2    \\
&+[2(1 + u_n ^2 )^{ - 1}u_n^2  - 2(1 + u_n ^2 )^{ - 2} u_n^4 ]u_{jn}^2 [u_t- u_{nn}]+ 4u_nu_{jn}u_{tj} - 6u_{jn}^2u_{t}  \\
& +2\sum_{i\in G}\sum_{\alpha=1, \a \ne i}^{n}
\frac{1}{u_{ii}}M_{\alpha\a}[u_nu_{ij\alpha}-2u_{i\a}u_{jn}]^2 + O(\mathcal H_{\phi}) \\
\leq& 4u_nu_{jn}u_{tj} - 6u_{jn}^2u_{t} + O(\mathcal H_{\phi}) \\
=&4 u_{jn}[u_t u_{jn} - \frac{\hat h_{jn}}{u_t}]^2 - 6u_{jn}^2u_{t} + O(\mathcal H_{\phi}) \\
=&-2 u_t [u_{jn} +\frac{\hat h_{jn}}{u_t^2} ]^2 + 2\frac{\hat h_{jn}^2}{u_t^3} + O(\mathcal H_{\phi}) \\
\leq& O(\mathcal H_{\phi})
\end{align*}

So we can get
\begin{eqnarray}\label{3.58}
M_{\alpha\beta}\phi_{\alpha\beta}- \phi_t &\le&  C(\phi+\sum_{i,j\in B}|\nabla
a_{ij}|) \nonumber\\
&&-\frac{1}{{\sigma}^3_1(B)}\sum_{\alpha,\beta=1}^n\sum_{i\in
B}M_{\alpha\beta}[{\sigma}_1(B)a_{ii,\alpha}-a_{ii}\sum_{j\in
B}a_{jj,\alpha}][{\sigma}_1(B)a_{ii,\beta}-a_{ii}\sum_{j\in
B}a_{jj,\beta}]\nonumber\\
&&-\frac{1}{{\sigma}_1(B)}\sum_{\alpha,\beta=1}^n\sum_{i\neq j,
i,j\in B}M_{\alpha\beta}a_{ij,\alpha}a_{ij,\beta},
\end{eqnarray}
Following the proof of Theorem \ref{th3.1}, we get \eqref{3.53}.
\qed

\begin{remark}\label{rem3.8}
The constant rank properties ( that is, Corollary \ref{cor3.2}) still holds for the equation \eqref{1.14}.
\end{remark}

\section{Constant rank theorem of the spacetime second fundamental form}
\setcounter{equation}{0} \setcounter{theorem}{0}

In this section, we start to consider the spacetime level sets $\hat \Sigma ^{c} = \{(x,t) \in
\Omega \times (0, T)|u(x,t) = c\}$, and as in Section 2, the Weingarten tensor is
\begin{equation}\label{4.1}
\hat a_{\a \b} =-\frac{|u_t|}{|D u|{u_t}^3} \hat A_{\a \b}, \quad  1 \leq \a, \b \leq n,
\end{equation}
where
\begin{equation}\label{4.2}
\hat A_{\a \b} = \hat h_{\a \b}
-\frac{u_\a u_\gamma \hat h_{\b \gamma}}{\hat W(1+\hat W)u_t^2} -\frac{u_\b u_\gamma \hat h_{\a \gamma}}{\hat W(1+ \hat W)u_t^2}
+ \frac{u_\a u_\b u_\gamma u_\eta \hat h_{\gamma \eta}}{\hat W^2(1+\hat W)^2u_t^4}, \quad \hat W = \frac{|D u|}{|u_t|}.
\end{equation}

Suppose $\hat{a}(x,t)=(\hat a_{\a \b})_{n \times n}$ attains the minimal
rank $l$ at some point $(x_0, t_0) \in \Omega \times (0, T]$. We may assume
$l\leqslant n-1$, otherwise there is nothing to prove. At $(x_0,t_0)$, we may choose $e_1,\cdots, e_{n-1},e_n$ such that
\begin{equation}\label{4.3}
 |\n u(x_0,t_0)|=u_n(x_0,t_0)>0\ \mbox{ and } \
\Big(u_{ij} \Big)_{1 \leq i,j\leq n-1} \mbox{ is diagonal at}\ (x_0,t_0).
\end{equation}
Without loss of generality we assume $ u_{11} \leq u_{22}\leq \cdots
\leq u_{n-1n-1} $. So, at $(x_0,t_0)$, from \eqref{4.1}, we have the matrix
$\Big(\hat a_{ij} \Big)_{1 \leq i,j \leq n-1}$ is also diagonal, and $\hat a_{11} \geq \hat a_{22} \geq \cdots \geq \hat a_{n-1
n-1}$. From lemma \ref{lem2.5}, there is a positive constant $
C_0$ such that at $(x_0,t_0)$

CASE 1:
\begin{eqnarray*}
&&\hat a_{11}  \geq \cdots \geq \hat a_{l-1l-1} \geq
C_0 , \quad \hat a_{ll} = \cdots = \hat a_{n-1n-1} =0 , \\
&&\hat a_{nn} -\sum\limits_{i = 1}^{l-1} {\frac{{\hat a_{in} ^2 }} {{\hat a_{ii}
}}} \geq C_0 ,  \quad \hat a_{in} = 0, \quad l \leqslant i \leqslant n-1.
\end{eqnarray*}

CASE 2:
\begin{eqnarray*}
&&\hat a_{11}  \geq \cdots \geq \hat a_{ll} \geq C_0 , \quad \hat a_{l+1l+1} =
\cdots = \hat a_{n-1n-1} =0, \\
&& \hat a_{tt}  = \sum\limits_{i = 1}^{l}{\frac{{\hat a_{in} ^2 }} {{\hat a_{ii} }}}
,  \quad \hat a_{in} = 0, \quad  l+1 \leqslant i \leqslant n-1.
\end{eqnarray*}

\subsection{ CASE 1 }

In this subsection, we consider CASE 1, that is, at $(x_0,t_0)$, we have
\begin{eqnarray*}
&&\hat a_{11}  \geq \cdots \geq \hat a_{l-1l-1} \geq
C_0 , \quad \hat a_{ll} = \cdots = \hat a_{n-1n-1} =0 , \\
&&\hat a_{nn} -\sum\limits_{i = 1}^{l-1} {\frac{{\hat a_{in} ^2 }} {{\hat a_{ii}
}}} \geq C_0 ,  \quad \hat a_{in} = 0, \quad l \leqslant i \leqslant n-1.
\end{eqnarray*}
Then there is a neighborhood $\mathcal {O}\times
(t_0-\delta, t_0]$ of $(x_0, t_0)$, such that for any fixed point  $(x,t) \in \mathcal {O}\times (t_0-\delta,
t_0]$, we may choose $e_1,\cdots, e_{n-1},e_n$ such that
\begin{equation}\label{4.4}
 |\n u(x,t)|=u_n(x,t)>0\ \mbox{ and } \
\Big(u_{ij} \Big)_{1 \leq i,j \leq n-1} \mbox{is diagonal at}\ (x,t).
\end{equation}
Similarly we assume $ u_{11} \leq u_{22}\leq \cdots
\leq u_{n-1n-1} $. So, at $(x,t) \in \mathcal {O}\times (t_0-\delta,
t_0]$, from \eqref{4.1}, we have the matrix
$\Big(\hat a_{ij} \Big)_{1 \leq i,j \leq n-1}$ is also diagonal, and $\hat a_{11} \geq \hat a_{22} \geq \cdots \geq \hat a_{n-1
n-1}$. There is a positive constant $C_0>0$ depending only on
$\| u \|_{C^{4}}$ and $\mathcal {O}\times (t_0-\delta, t_0]$,
such that
\begin{eqnarray*}
&&\hat a_{11}  \geq \cdots \geq \hat a_{l-1l-1} \geq
C_0 ,  \\
&&\hat a_{nn} -\sum\limits_{i = 1}^{l-1} {\frac{{\hat a_{in} ^2 }} {{\hat a_{ii}
}}} \geq C_0.
\end{eqnarray*}
for $(x,t) \in \mathcal {O}\times (t_0-\delta, t_0]$. For convenience we
denote $ G = \{ 1, \cdots ,l-1 \} $ and $ B = \{
l, \cdots, n-1 \} $ be the "good" and "bad" sets of indices
respectively.

Since
\begin{equation}\label{4.5}
\hat a_{ij}= \frac{|\n u|}{|Du|} a_{ij}, 1 \leq i, j \leq n-1,
\end{equation}
there is a positive constant $C>0$ depending only on
$\|u\|_{C^{4}}$ and $\mathcal {O}\times (t_0-\delta, t_0]$,
such that
\begin{equation}\label{4.6}
a_{11}  \geq \cdots \geq  a_{l-1l-1} \geq C,\quad (x,t) \in \mathcal {O}\times (t_0-\delta, t_0],
\end{equation}
and
\begin{equation}\label{4.7}
a_{ll}(x_0,t_0) = \cdots = a_{n-1n-1}(x_0,t_0) = 0.
\end{equation}
So the spatial
second fundamental form $a=(a_{ij})_{n-1 \times n-1}$ attains the minimal rank $l-1$ at $(x_0,t_0)$.
From Theorem $\ref{th3.1}$, the constant rank theorem holds for the spatial
second fundamental form $a$. So we can get $a_{ii}=0, \forall i \in B$ for any $(x,t) \in \mathcal {O} \times (t_0-\delta, t_0]$.
Furthermore,
\begin{equation}\label{4.8}
\hat a_{ii}=0, \forall i \in B.
\end{equation}
We denote $M=(\hat a_{ij})_{1 \leq i, j \leq n-1}$,
so
\begin{equation} \label{4.9}
\sigma_{l+1}(M)=\sigma_{l}(M) \equiv 0, \quad \text{for
every}\quad (x,t) \in \mathcal {O} \times (t_0-\delta, t_0].
\end{equation}
Then we have
\begin{equation} \label{4.10}
0 \leq \sigma_{l+1}(\hat a)\leq \sigma_{l+1}(M)+\hat a_{nn}\sigma_l(M) = 0.
\end{equation}
So \begin{equation} \label{4.11}
\sigma_{l+1}(\hat a)\equiv 0, \quad \text{for every}\quad
(x,t) \in \mathcal {O} \times (t_0-\delta, t_0].
\end{equation}

By the continuity method, Theorem \ref{th1.2} holds under the CASE 1.

\subsection{ CASE 2 }

In this subsection, we consider CASE 2. From Lemma \ref{lem2.5}, $\hat{a}(x,t)=(\hat a_{\a \b})_{n \times n}$ attains the minimal
rank $l$ at some point $(x_0, t_0) \in \Omega \times (0, T]$ and at $(x_0,t_0)$, we may choose $e_1,\cdots, e_{n-1},e_n$ such that
\begin{equation}
 |\n u|=u_n>0\ \mbox{ and } \
(u_{ij})_{1 \leq i, j \leq n-1} \mbox{ is diagonal at}\ (x_0,t_0). \notag
\end{equation}
Then we have
\begin{eqnarray*}
&&\hat a_{11}  \geq \cdots \geq \hat a_{ll} \geq C_0 , \quad \hat a_{l+1l+1} =
\cdots = \hat a_{n-1n-1} =0, \\
&& \hat a_{nn}  = \sum\limits_{i = 1}^{l}{\frac{{\hat a_{in} ^2 }} {{\hat a_{ii} }}}
,  \quad \hat a_{in} = 0, \quad  l+1 \leqslant i \leqslant n-1.
\end{eqnarray*}
Then there is a small enough neighborhood $\mathcal {O}\times
(t_0-\delta, t_0]$ of $(x_0, t_0)$, such that for any fixed point  $(x,t) \in \mathcal {O}\times (t_0-\delta,
t_0]$, we may choose $e_1,\cdots, e_{n-1},e_n$ such that
\begin{equation}\label{4.12}
 |\n u(x,t)|=u_n(x,t)>0\ \mbox{ and }\
(u_{ij})_{1 \leq i, j \leq n-1} \mbox{ is diagonal at}\ (x,t).
\end{equation}
Similarly we assume $ u_{11} \leq u_{22}\leq \cdots
\leq u_{n-1n-1} $. So, at $(x,t) \in \mathcal {O}\times (t_0-\delta,
t_0]$, from \eqref{4.1}, we have the matrix
$\Big(\hat a_{ij}\Big)_{1 \leq i,j\leq n-1}$ is also diagonal, and $\hat a_{11} \geq \hat a_{22} \geq \cdots \geq \hat a_{n-1
n-1}$. There is a positive constant $C>0$ depending only on
$\|u\|_{C^{4}}$ and $\mathcal {O}\times (t_0-\delta, t_0]$,
such that
\begin{align*}
\hat a_{11}  \geq \cdots \geq \hat a_{ll} \geq C,
\end{align*}
for all $(x,t) \in \mathcal {O}\times (t_0-\delta, t_0]$.
 For convenience we
denote $ G = \{ 1, \cdots ,l \} $ and $ B = \{
l+1, \cdots, n-1 \} $ be the "good" and "bad" sets of indices
respectively. Since
\begin{equation*}
\hat a_{ij}= \frac{|\n u|}{|Du|} a_{ij}, 1 \leq i, j \leq n-1,
\end{equation*}
there is a positive constant $C>0$ depending only on
$\|u\|_{C^{4}}$ and $\mathcal {O}\times (t_0-\delta, t_0]$,
such that
\begin{equation}\label{4.13}
a_{11}  \geq \cdots \geq  a_{l-1l-1} \geq C,\quad (x,t) \in \mathcal {O}\times (t_0-\delta, t_0],
\end{equation}
and
\begin{equation*}
a_{ll}(x_0,t_0) = \cdots = a_{n-1n-1}(x_0,t_0) = 0.
\end{equation*}
So the spatial
second fundamental form $a=(a_{ij})_{n-1 \times n-1}$ attains the minimal rank $l$ at $(x_0,t_0)$.
From Theorem $\ref{th3.1}$, the constant rank theorem holds for the spatial
second fundamental form $a$. So we can get $a_{ii}=0, \forall i \in B$ for any $(x,t) \in \mathcal {O} \times (t_0-\delta, t_0]$.
Furthermore, $u_{x_i x_i}=0, \forall i \in B$.

In order to simplify the calculations, we need a new spacetime coordinate system.
For any fixed point  $(x,t) \in \mathcal {O}\times (t_0-\delta,
t_0]$, $\{e_1,\cdots, e_{n-1},e_n\}$ is the coordinate satisfying \eqref{4.12} and \eqref{4.13}, and $e_{n+1}$ is the time coordinate.
First, by translating $\{e_n, e_{n+1} \}$, we  get the coordinate $\{e_1,\cdots, e_{n-1}, \hat e_n, \hat e_{n+1}\}$ with
\begin{align}\label{4.14}
z = (z_1 , \cdots ,z_n ,z_{n + 1} ) = (x,t)O,
\end{align}
where
\begin{align}\label{4.15}
O = \left( {O_{ab} } \right)_{n + 1 \times n + 1}  = \left( {\begin{array}{*{20}c}
   1 & {} & {} & {} & {}  \\
   {} &  \ddots  & {} & {} & {}  \\
   {} & {} & 1 & {} & {}  \\
   {} & {} & {} & {\cos \theta } & {\sin \theta }  \\
   {} & {} & {} & { - \sin \theta } & {\cos \theta }  \\
\end{array}} \right),
\end{align}
such that
\begin{equation}\label{4.16}
 u_{z_{n+1}} = |D u| > 0,\quad  u_{z_1} = \cdots =  u_{z_n} = 0, \quad \text{ at } (x,t).
\end{equation}
So from \eqref{1.4}, we have $ \theta \in (0, \frac{\pi}{2})$. Now we fix the coordinate $\{e_{l+1},\cdots, e_{n-1}, \hat e_{n+1}\}$ and translate
$\{e_1,\cdots, e_{l}, \hat e_n\}$, and we get the coordinates $\{\bar e_1, \cdots, \bar e_{l}, e_{l+1},\cdots, e_{n-1}, \bar e_n, \hat e_{n+1}\}$ with
\begin{align}\label{4.17}
 y = (y_1 , \cdots ,y_n ,y_{n + 1} ) = (z_1 , \cdots ,z_n ,z_{n + 1} )T,
 \end{align}
where
\begin{align}\label{4.18}
 T = \left( {T_{\alpha \beta } } \right)_{n + 1 \times n + 1}  = \left( {\begin{array}{*{20}c}
   {T_{11} } &  \cdots  & {T_{1l} } & {} & {} & {} & {T_{1n} } & {}  \\
    \vdots  &  \ddots  &  \vdots  & {} & {} & {} &  \vdots  & {}  \\
   {T_{l1} } &  \cdots  & {T_{ll} } & {} & {} & {} & {T_{ln } } & {}  \\
   {} & {} & {} & 1 & {} & {} & 0 & {}  \\
   {} & {} & {} & {} &  \ddots  & {} &  \vdots  & {}  \\
   {} & {} & {} & {} & {} & 1 & 0 & {}  \\
   {T_{n1} } &  \cdots  & {T_{nl} } & 0 &  \cdots  & 0 & 1 & {}  \\
   {} & {} & {} & {} & {} & {} & {} & 1  \\
\end{array}} \right)
\end{align}
such that
\begin{equation}\label{4.19}
\bigg( u_{y_{\a} y_{\b}} \bigg)_{1 \leq \a, \b\leq n} \mbox{is diagonal at}\ (x,t).
\end{equation}
Finally, we get a new spacetime coordinate $\{\bar e_1, \cdots, \bar e_{l}, e_{l+1},\cdots, e_{n-1}, \bar e_n, \hat e_{n+1}\}$ with
\begin{align}\label{4.20}
 y = (y_1 , \cdots ,y_n ,y_{n + 1} ) = (x,t)P, \quad P = OT
\end{align}
such that
\begin{align}
\label{4.21}&u_{y_{n+1}} = |D u| > 0,\quad  u_{y_1} = \cdots =  u_{y_n} = 0, \quad \text{ at } (x,t), \\
\label{4.22}& \bigg( u_{y_{\a} y_{\b}} \bigg)_{1 \leq \a, \b\leq n} \mbox{is diagonal at}\ (x,t).
\end{align}
From \eqref{2.13}-\eqref{2.15}, we get
\begin{equation}\label{4.23}
\bar a_{\a \b} =-\frac{1}{{u_{y_{n+1}}}} u_{y_{\a} y_{\b}}, \quad  1 \leq \a, \b \leq n,
\end{equation}
Without loss of generality, we can assume $ \frac{\partial ^2 u} {\partial y_1 \partial y_1} \leqslant \cdots \leqslant
\frac{\partial ^2 u} {\partial y_l \partial y_l} \leqslant -C < 0 $, where the positive
constant $C > 0$ depending only on $\left\| u \right\|_{C^{3,1} }$. Then we have
\begin{align}\label{4.24}
\bar a_{11}  \geq \cdots \geq \bar a_{ll} \geq C.
\end{align}
For convenience we denote $ G = \{1, \cdots , l\} $ and $ B =
\{ l+1, \cdots ,n-1\} $ which mean good terms and bad ones of
indices respectively. Without confusion we will also simply denote $ G = \{ \bar a_{11}, \cdots , \bar a_{ll} \} $ and $
B = \{ \bar a_{l+1 l+1}  , \cdots ,\bar a_{n-1 n-1} \} $.
We set
\begin{eqnarray}\label{4.25}
\phi = \sigma_{l+1}(\bar a).
\end{eqnarray}
In the following, we will prove a differential inequality
\begin{align}\label{4.26}
\sum_{ij=1}^n F^{ij} \phi_{ij}  - \phi _t \le  C(\phi  + |\nabla _x \phi |) \quad \text{ in } \mathcal {O}\times (t_0-\delta,
t_0].
\end{align}
In fact, if $t_0 = T$ and $(x, t) \in \mathcal {O}\times \{ t_0\}$, we only have \eqref{2.30} instead of \eqref{2.29} from Lemma \ref{lem2.9} ( see Remark \ref{rem2.10}). So in order to utilizing \eqref{2.29}, we just prove \eqref{4.26} holds for any $(x, t) \in \mathcal {O}\times (t_0-\delta,
t_0)$, with a constant $C$ independent of $ dist(\mathcal {O}\times (t_0-\delta,
t_0], \partial(\Omega \times (0, T]) )$ and then by a approximation, \eqref{4.26} holds for $t = t_0$.
Then by the strong maximum principle and the method of continuity, we can prove Theorem \ref{th1.2} under CASE 2.

For convenience, we will use $i,j,k,l =1, \cdots, n$ to represent the $x$ coordinates, $t$ still the time coordinate, and
 $\alpha, \beta, \gamma, \eta =1, \cdots, n+1$ the $y$ coordinates. And we have
\begin{align}
\label{4.27}&\frac{{\partial y_\alpha  }}{{\partial x_i }} = P_{i\alpha }  \\
\label{4.28}&\frac{{\partial y_\alpha  }}{{\partial t}} = P_{n + 1\alpha }
\end{align}
In the following, we always denote
\begin{align*}
& u_{i}=\frac{{\partial u}}{{\partial x_i }}, u_{t}=\frac{{\partial u}}{{\partial t }} , u_{\alpha}=\frac{{ \partial u}}{{\partial y_\alpha}},  u_{n+1}=\frac{{ \partial u}}{{\partial y_{n+1}}},\\
& u_{ij}=\frac{{\partial^2 u}}{{\partial x_i \partial x_j}}, u_{it}=\frac{{\partial^2 u}}{{\partial x_i \partial t }}, u_{tt}=\frac{{\partial^2 u}}{{ \partial t ^2 }}, u_{i\alpha}=\frac{{\partial^2 u}}{{\partial x_i \partial y_\alpha}}, \\
&u_{\alpha t}=\frac{{\partial^2 u}}{{\partial y_\alpha \partial t}},u_{\alpha \beta}=\frac{{ \partial^2 u}}{{\partial y_\alpha \partial y_\beta}}, \text{ etc. }
\end{align*}
Also, we will use notion $h=O(\phi)$ if $|h(x,t)| \le C (\phi) $ for $(x,t) \in \mathcal {O}\times (t_0-\delta, t_0)$ with positive constant $C$ under
control, and $h=O(\phi + | \nabla \phi| )$ has a similar meaning.

From the above discussions, for any $(x,t) \in \mathcal {O} \times (t_0-\delta, t_0)$ with the coordinate \eqref{4.21} and \eqref{4.22}, we get
\begin{align}\label{4.29}
 u_{\alpha \alpha }= \frac{\partial ^2 u} {\partial y_{\alpha} \partial y_{\alpha}} = \frac{\partial ^2 u} {\partial x_{\alpha} \partial x_{\alpha}}= 0,  \quad \forall \alpha  \in B.
\end{align}

Under above assumptions, we can get
\begin{proposition}\label{pro4.9}
For any $(x,t) \in \mathcal {O} \times (t_0-\delta, t_0)$ with the coordinate \eqref{4.21} and \eqref{4.22}, we can get
\begin{align}\label{4.30}
\bar a_{\alpha \alpha}(x,t) =0, \quad  \alpha \in B.
\end{align}
Furthermore, we have from the semipositive definite of $\bar a$,
\begin{eqnarray}
\label{4.31}&&\bar a_{\alpha \beta}(x,t) = 0, \quad  \alpha \in B, \beta  \in B\cup G, \\
\label{4.32}&& \bar a_{\alpha n}(x,t)=\hat a_{n \alpha}(x,t) = 0, \quad  \alpha \in B, \\
\label{4.33} && {D \bar a_{\alpha \beta} }(x,t) = 0, \quad  \alpha \in B,  \beta \in B \cup G,\\
\label{4.34}&& {D \bar a_{\alpha n} }(x,t) = 0, \quad  \alpha \in B.
\end{eqnarray}
\end{proposition}

PROOF. The proof is directly from the constant rank theorem of $a$ and Lemma \ref{lem2.9}.

For any $(x,t) \in \mathcal {O} \times (t_0-\delta, t_0)$ with the coordinate \eqref{4.21} and \eqref{4.22}, we can get from \eqref{4.29}
\begin{equation*}
\bar a_{\a \a}(x,t)=- \frac{1}{|D u|} u_{\a \a}=0, \forall \a \in B.
\end{equation*}
From the positive definite of $\bar a$ at $(x,t)$, we get
\begin{align*}
&\bar a_{\alpha \beta}(x,t) = 0, \quad  \alpha \in B, \beta  \in B\cup G, \\
&\bar a_{\a n}(x,t)=0, \forall \a \in B.
\end{align*}
And from Lemma \ref{lem2.9} (i.e. Remark \ref{rem2.10}), we get
\begin{align*}
&|D \bar a_{\alpha \beta}|(x,t) = 0, \quad  \alpha \in B, \beta  \in B\cup G, \\
&|D \bar a_{\a n}|(x,t)=0, \forall \a \in B.
\end{align*}
So the lemma holds.

\qed

\begin{lemma}\label{lem4.10}
\begin{align}
\label{4.35}& D u_{\alpha}=0, \quad \alpha \in B,  \\
\label{4.36}& D u_{\alpha \beta}=0, \quad \alpha \in B, \beta \in B, \\
\label{4.37}& D u_{\alpha \beta } =0, \quad \alpha \in B, \beta \in G \cup \{ n \}.
\end{align}
\end{lemma}

PROOF. By the constant rank properties Corollary \ref{cor3.2}, \eqref{4.35} holds since the $y_\a$ coordinate is
the $x_\a$ coordinate for $\a \in B$.

By \eqref{2.13}, \eqref{2.14}, \eqref{4.33} and \eqref{4.34}, we get for $\alpha \in B$ and $\beta = 1, \cdots, n$,
\begin{align*}
0={D \bar a_{\alpha \beta} } =-\frac{|u_{n+1}|}{|D u|{u_{n+1}}^3} {D \bar A_{\alpha \beta} }=-\frac{|u_{n+1}|}{|D u|{u_{n+1}}^3}
{D \bar h_{\alpha \beta} },
\end{align*}
so from \eqref{2.15}, we get
\begin{align*}
0=D \bar h_{\alpha \beta}=& u_{n+1}^2 D u_{\alpha \beta} + 2u_{n+1} D u_{n+1} u_{\alpha \beta}- u_{n+1} u_{\alpha  n+1} D u_{\beta}- u_{n+1} u_{\b n+1} D u_{\a} \\
=& u_{n+1}^2 D u_{\alpha \beta}, \qquad \qquad \alpha \in B, \beta =1, \cdots, n.
\end{align*}
Hence the lemma holds. \qed

\begin{lemma}\label{lem4.11}
\begin{align}
\label{4.38}& u_{y_ny_n}=O( \phi), \\
\label{4.39}& u_{x_i y_n} = O( \phi ), i <n; \\
\label{4.40}& u_{y_n y_n x_i}=O( \phi + |\nabla \phi |) , \quad i= 1, \cdots, n-1, \\
\label{4.41}& u_{y_n y_n x_n}=2 \frac{1}{u_{y_{n+1}}} u_{y_n x_n} u_{y_ny_{n+1}} + O( \phi + |\nabla \phi |)
\end{align}
\end{lemma}
Proof: In the $y$ coordinates, we have from \eqref{4.30}
\begin{align}
\phi  = \sigma _{l + 1} (\overline a ) = \sigma _l (G)\overline a _{nn}  \ge 0, \notag
\end{align}
so we have
\begin{align}\label{4.42}
\overline a _{nn}  = O(\phi ).
\end{align}
By \eqref{2.13}-\eqref{2.15}, we have
\begin{align*}
\bar a_{nn}  =-\frac{|u_{y_{n+1}}|}{|D u|{u_{y_{n+1}}}^3} {\bar A_{nn} }=-\frac{|u_{y_{n+1}}|}{|D u|{u_{y_{n+1}}}^3}
{ \bar h_{nn} }=-\frac{|u_{y_{n+1}}|}{|D u|{u_{y_{n+1}}}^3}{ {u_{n+1}}^2 u_{y_n y_n} },
\end{align*}
so
\begin{equation}\label{4.43}
\bar A_{nn} = O(\phi ), \quad \bar h_{nn} = O(\phi ),  \quad u_{y_{n} y_{n}}= O(\phi ).
\end{equation}

Taking the first derivatives of $\phi$ in $x$, we have
\begin{align}\label{4.44}
\phi _i  =& \frac{{\partial \phi }}{{\partial x_i }} = \sum\limits_{\alpha  = 1}^n {\sigma _l (\overline a |\alpha )\overline a _{\alpha \alpha ,i} }  \notag  \\
=& \sum\limits_{\alpha  \in G} {\sigma _l (\overline a |\alpha )\bar a_{\alpha \alpha ,i} }  + \sum\limits_{\alpha  \in B} {\sigma _l (\overline a |\alpha )\overline a _{\alpha \alpha ,i} }  + \sum\limits_{\alpha  = n} {\sigma _l (\overline a |\alpha )\overline a _{\alpha \alpha ,i} }  \notag  \\
=& \sigma _l (G)\overline a _{nn,i}  + O(\phi ),
\end{align}
and from \eqref{2.13}-\eqref{2.15}
\begin{align}\label{4.45}
\bar a_{nn, i}  =&(-\frac{|u_{n+1}|}{|D u|{u_{n+1}}^3} )_i { \bar A_{nn} }-\frac{|u_{n+1}|}{|D u|{u_{n+1}}^3} { \bar A_{nn, i} }  = O(\phi )-\frac{|u_{n+1}|}{|D u|{u_{n+1}}^3} { \bar A_{nn, i} }  \notag \\
=&-\frac{|u_{n+1}|}{|D u|{u_{n+1}}^3}{ \bar h_{nn,i} }+ O(\phi )\notag \\
=&-\frac{|u_{n+1}|}{|D u|{u_{n+1}}^3}[ u_{y_{n+1}}^2  u_{y_ny_{n} x_i} - 2 u_{y_{n+1}} u_{y_n x_i} u_{y_ny_{n+1}}  ]+ O(\phi ),
\end{align}
so
\begin{equation}\label{4.46}
\bar a_{nn, i} = O( \phi + |\nabla \phi |), \quad  \bar A_{nn,i} = O( \phi + |\nabla \phi |), \quad \bar h_{nn,i} = O( \phi + |\nabla \phi |),
\end{equation}
and
\begin{align}\label{4.47}
 u_{y_n y_n x_i}=2 \frac{1}{u_{y_{n+1}}} u_{y_n x_i} u_{y_ny_{n+1}} + O( \phi + |\nabla \phi |).
\end{align}
It is easy to know for $i= 1, \cdots, n-1$,
\begin{align}\label{4.48}
u_{y_n x_i} =& u_{y_ny_{\a}} \frac{\partial y_\a }{ \partial x_i } = u_{y_ny_{n}} \frac{\partial y_n }{ \partial x_i }+ u_{y_ny_{n+1}} \frac{\partial y_{n+1} }{ \partial x_i }= O( \phi )+ u_{y_ny_{n+1}} \frac{\partial y_{n+1} }{ \partial x_i }  \notag  \\
=& u_{y_ny_{n+1}} \frac{\partial z_{n+1} }{ \partial x_i }+ O( \phi )= 0+ O( \phi ).
\end{align}
Hence the lemma holds from \eqref{4.47} and \eqref{4.48}.
\qed

\begin{lemma}\label{lem4.12}
\begin{align}
\label{4.49}&u_{y_{n+1}} u_{x_n y_n}= u_{x_n} u_{y_ny_{n+1}}+ O( \phi ),  \\
\label{4.50}&u_{y_{n+1}} u_{y_n t} = u_t u_{y_ny_{n+1}} + O( \phi ).
\end{align}
\end{lemma}

PROOF. By the chain rule, we get
\begin{align*}
u_{x_n} u_{y_ny_{n+1}}= u_{y_\a} \frac{\partial y_\a }{ \partial x_n } u_{y_n y_{n+1}}= u_{y_{n+1}} \frac{\partial y_{n+1} }{ \partial x_n } u_{y_n y_{n+1}},
\end{align*}
and
\begin{align*}
u_{y_{n+1}} u_{x_n y_n}= u_{y_{n+1}} u_{y_\a y_n } \frac{\partial y_\a }{ \partial x_n }= u_{y_{n+1}} u_{y_n y_n } \frac{\partial y_n }{ \partial x_n }+ u_{y_{n+1}} u_{y_{n+1} y_n } \frac{\partial y_{n+1} }{ \partial x_n } = O( \phi )+ u_{y_{n+1}} u_{y_{n+1} y_n } \frac{\partial y_{n+1} }{ \partial x_n }.
\end{align*}
so \eqref{4.49}holds.

Similarly, we have
\begin{align*}
u_{t} u_{y_n y_{n+1}} =u_{y_{\a}} \frac{\partial y_\a }{ \partial t }  u_{y_ny_{n+1}} = u_{y_{n+1}} \frac{\partial y_{n+1} }{ \partial t } u_{y_ny_{n+1}} ,
\end{align*}
and
\begin{align*}
u_{y_{n+1}} u_{y_n t} =u_{y_{n+1}}  u_{y_ny_{\a}} \frac{\partial y_\a }{ \partial t } = u_{y_ny_{n}} \frac{\partial y_n }{ \partial t }+ u_{y_{n+1}} u_{y_ny_{n+1}} \frac{\partial y_{n+1} }{ \partial t }= O( \phi )+ u_{y_{n+1}} u_{y_ny_{n+1}} \frac{\partial y_{n+1} }{ \partial t }.
\end{align*}
so \eqref{4.50} holds.
\qed

\begin{lemma}\label{lem4.13}
\begin{align}\label{4.51}
\sum\limits_{kl = 1}^n {F^{kl} } u_{kl y_\gamma y_\gamma}= 0,  \quad  \text{ for } \gamma \in B.
\end{align}
\end{lemma}

PROOF. From \eqref{3.31} (i.e. the constant rank properties Corollary \ref{cor3.2}) and \eqref{4.36}, we have for $\gamma \in B$
\begin{align}\label{4.52}
\sum\limits_{kl = 1}^n {F^{kl} } u_{kl y_\gamma y_\gamma}=& = \sum\limits_{kl = 1}^n {F^{kl} } u_{kl y_\gamma y_\gamma} -  u_{y_\gamma y_\gamma t} \notag\\
=& \sum\limits_{kl = 1}^n {F^{kl} } u_{kl x_\gamma x_\gamma} -  u_{x_\gamma x_\gamma t}= 2\sum_{i\in G}\sum_{k,l=1}^{n}
F^{kl}\frac{u_{x_n}^2 u_{ik x_\gamma} u_{il x_\gamma} }{u_{ii}}.
\end{align}
In fact, for $i \in G, \gamma \in B$, we have from \eqref{4.36} and \eqref{4.37},
\begin{align}\label{4.53}
u_{ik x_\gamma} =u_{ik y_\gamma} =u_{ y_\a x_k y_\gamma} \frac{\partial y_\a}{\partial x_i} =\sum_{\alpha \leq n} \frac{\partial u_{ y_\a  y_\gamma}}{\partial x_k}  \frac{\partial y_\a}{\partial x_i} + u_{ y_{n+1} x_k y_\gamma} \frac{\partial y_{n+1}}{\partial x_i}
=0 + u_{ y_{n+1} k y_\gamma} \frac{\partial z_{n+1}}{\partial x_i}=0.
\end{align}
So \eqref{4.51} holds from \eqref{4.52} and \eqref{4.53}.

\qed

\begin{lemma}\label{lem4.14}

\begin{align}\label{4.54}
\sum\limits_{ij = 1}^n {F^{ij} \phi _{ij} }  =\sigma _l (G)\sum\limits_{ij = 1}^n {F^{ij} \overline a _{nn,ij} }  - 2\sigma _l (G)\sum\limits_{ij = 1}^n {F^{ij} \sum\limits_{\eta  \in G} {\frac{{\overline a _{n\eta ,i} \overline a _{\eta n,j} }}{{\overline a _{\eta \eta } }}} } + O( \phi + |\nabla_x \phi |).
\end{align}

\end{lemma}

PROOF. Taking the second derivatives of $\phi$ in $y$ coordinates, we have
\begin{align}\label{4.55}
\phi _{\alpha \beta }  =& \frac{{\partial ^2 \phi }}{{\partial y_\alpha  \partial y_\beta  }} = \sum\limits_{\gamma  = 1}^n {\frac{{\partial \sigma _{l + 1} }}{{\partial \overline a _{\gamma \gamma } }}\overline a _{\gamma \gamma ,\alpha \beta } }  + \sum\limits_{\gamma  \ne \eta } {\frac{{\partial ^2 \sigma _{l + 1} }}{{\partial \overline a _{\gamma \gamma } \partial \overline a _{\eta \eta } }}\overline a _{\gamma \gamma ,\alpha } \overline a _{\eta \eta ,\beta } }  + \sum\limits_{\gamma  \ne \eta } {\frac{{\partial ^2 \sigma _{l + 1} }}{{\partial \overline a _{\gamma \eta } \partial \overline a _{\eta \gamma } }}\overline a _{\gamma \eta ,\alpha } \overline a _{\eta \gamma ,\beta } }   \notag \\
=& \sum\limits_{\gamma  = 1}^n {\sigma _l (\overline a |\gamma )\overline a _{\gamma \gamma ,\alpha \beta } }  + \sum\limits_{\gamma  \ne \eta } {\sigma _{l - 1} (\overline a |\gamma \eta )\overline a _{\gamma \gamma ,\alpha } \overline a _{\eta \eta ,\beta } }  - \sum\limits_{\gamma  \ne \eta } {\sigma _{l - 1} (\overline a |\gamma \eta )\overline a _{\gamma \eta ,\alpha } \overline a _{\eta \gamma ,\beta } },
\end{align}
where
\begin{align}\label{4.56}
\sum\limits_{\gamma  = 1}^n {\sigma _l (\overline a |\gamma )\overline a _{\gamma \gamma ,\alpha \beta } }  =& \sum\limits_{\gamma  \in G} {\sigma _l (\overline a |\gamma )\overline a _{\gamma \gamma ,\alpha \beta } }  + \sum\limits_{\gamma  \in B} {\sigma _l (\overline a |\gamma )\overline a _{\gamma \gamma ,\alpha \beta } }  + \sum\limits_{\gamma  = n} {\sigma _l (\overline a |\gamma )\overline a _{\gamma \gamma ,\alpha \beta } }   \notag \\
=& \sigma _l (G)\sum\limits_{\gamma  \in B} {\overline a _{\gamma \gamma ,\alpha \beta } }  + \sigma _l (G)\overline a _{nn,\alpha \beta }  + O(\phi )
\end{align}

\begin{align}\label{4.57}
\sum\limits_{\gamma  \ne \eta } {\sigma _{l - 1} (\overline a |\gamma \eta )\overline a _{\gamma \gamma ,\alpha } \overline a _{\eta \eta ,\beta } }  =& \sum\limits_{\scriptstyle \gamma \eta  \in G \hfill \atop \scriptstyle \gamma  \ne \eta  \hfill} {\sigma _{l - 1} (\overline a |\gamma \eta )\overline a _{\gamma \gamma ,\alpha } \overline a _{\eta \eta ,\beta } }  + \sum\limits_{\scriptstyle \gamma  = n \hfill \atop
\scriptstyle \eta  \in G \hfill} {\sigma _{l - 1} (\overline a |\gamma \eta )\overline a _{\gamma \gamma ,\alpha } \overline a _{\eta \eta ,\beta } }  \notag  \\
&+ \sum\limits_{\scriptstyle \gamma  \in G \hfill \atop \scriptstyle \eta  = n \hfill} {\sigma _{l - 1} (\overline a |\gamma \eta )\overline a _{\gamma \gamma ,\alpha } \overline a _{\eta \eta ,\beta } }   \notag \\
=& O(\phi ) + \sum\limits_{\eta  \in G} {\sigma _{l - 1} (G|\eta )\overline a _{\eta \eta ,\beta } } \overline a _{nn,\alpha }  + \sum\limits_{\gamma  \in G} {\sigma _{l - 1} (G|\gamma )\overline a _{\gamma \gamma ,\alpha } } \overline a _{nn,\beta }  \notag  \\
=& \sigma _l (G)[\sum\limits_{\eta  \in G} {\frac{{\overline a _{\eta \eta ,\beta } }}{{\overline a _{\eta \eta } }}} \overline a _{nn,\alpha }  + \sum\limits_{\gamma  \in G} {\frac{{\overline a _{\gamma \gamma ,\alpha } }}{{\overline a _{\gamma \gamma } }}} \overline a _{nn,\beta } ] + O(\phi )
\end{align}
and
\begin{align}\label{4.58}
\sum\limits_{\gamma  \ne \eta } {\sigma _{l - 1} (\overline a |\gamma \eta )\overline a _{\gamma \eta ,\alpha } \overline a _{\eta \gamma ,\beta } }  =& \sum\limits_{\scriptstyle \gamma \eta  \in G \hfill \atop \scriptstyle \gamma  \ne \eta  \hfill} {\sigma _{l - 1} (\overline a |\gamma \eta )\overline a _{\gamma \eta ,\alpha } \overline a _{\eta \gamma ,\beta } }  + \sum\limits_{\scriptstyle \gamma  = n \hfill \atop
\scriptstyle \eta  \in G \hfill} {\sigma _{l - 1} (\overline a |\gamma \eta )\overline a _{\gamma \eta ,\alpha } \overline a _{\eta \gamma ,\beta } }   \notag \\
&+ \sum\limits_{\scriptstyle \gamma  \in G \hfill \atop\scriptstyle \eta  = n \hfill} {\sigma _{l - 1} (\overline a |\gamma \eta )\overline a _{\gamma \eta ,\alpha } \overline a _{\eta \gamma ,\beta } }   \notag \\
=& O(\phi ) + \sum\limits_{\eta  \in G} {\sigma _{l - 1} (G|\eta )\overline a _{n\eta ,\alpha } \overline a _{\eta n,\beta } }  + \sum\limits_{\gamma  \in G} {\sigma _{l - 1} (G|\gamma )\overline a _{\gamma n,\alpha } \overline a _{n\gamma ,\beta } }   \notag \\
=& 2\sigma _l (G)\sum\limits_{\eta  \in G} {\frac{{\overline a _{n\eta ,\alpha } \overline a _{\eta n,\beta } }}{{\overline a _{\eta \eta } }}}  + O(\phi ) \end{align}
So we have
\begin{align}\label{4.59}
\phi _{\alpha \beta }  =& \sigma _l (G)\sum\limits_{\gamma  \in B} {\overline a _{\gamma \gamma ,\alpha \beta } }  + \sigma _l (G)\overline a _{nn,\alpha \beta }  - 2\sigma _l (G)\sum\limits_{\eta  \in G} {\frac{{\overline a _{n\eta ,\alpha } \overline a _{\eta n,\beta } }}{{\overline a _{\eta \eta } }}}  \notag \\
&+ \sigma _l (G)[\sum\limits_{\eta  \in G} {\frac{{\overline a _{\eta \eta ,\beta } }}{{\overline a _{\eta \eta } }}} \overline a _{nn,\alpha }  + \sum\limits_{\gamma  \in G} {\frac{{\overline a _{\gamma \gamma ,\alpha } }}{{\overline a _{\gamma \gamma } }}} \overline a _{nn,\beta } ] + O(\phi )
\end{align}
So we have
\begin{align}\label{4.60}
\sum\limits_{ij = 1}^n {F^{ij} \phi _{ij} }  =& \sum\limits_{ij = 1}^n {F^{ij} \sum\limits_{\alpha \beta  = 1}^{n + 1} {P_{i\alpha } P_{j\beta } \phi _{\alpha \beta } } }\notag  \\
=& \sigma _l (G)\sum\limits_{\gamma  \in B} {\sum\limits_{ij = 1}^n {F^{ij} \sum\limits_{\alpha \beta  = 1}^{n + 1} {P_{i\alpha } P_{j\beta } \overline a _{\gamma \gamma ,\alpha \beta } } } }  + \sigma _l (G)\sum\limits_{ij = 1}^n {F^{ij} \sum\limits_{\alpha \beta  = 1}^{n + 1} {P_{i\alpha } P_{j\beta } } \overline a _{nn,\alpha \beta } } \notag \\
&- 2\sigma _l (G)\sum\limits_{ij = 1}^n {F^{ij} \sum\limits_{\eta  \in G} {\frac{{[\sum\limits_{\alpha  = 1}^{n + 1} {P_{i\alpha } \overline a _{n\eta ,\alpha } } ][\sum\limits_{\beta  = 1}^{n + 1} {P_{j\beta } \overline a _{\eta n,\beta } } ]}}{{\overline a _{\eta \eta } }}} } \notag \\
&+ \sigma _l (G)\sum\limits_{ij = 1}^n {F^{ij} [\sum\limits_{\eta  \in G} {\frac{{\sum\limits_{\beta  = 1}^{n + 1} {P_{j\beta } } \overline a _{\eta \eta ,\beta } }}{{\overline a _{\eta \eta } }}} \sum\limits_{\alpha  = 1}^{n + 1} {P_{i\alpha } \overline a _{nn,\alpha } }  + \sum\limits_{\gamma  \in G} {\frac{{\sum\limits_{\alpha  = 1}^{n + 1} {P_{i\alpha } } \overline a _{\gamma \gamma ,\alpha } }}{{\overline a _{\gamma \gamma } }}} \sum\limits_{\beta  = 1}^{n + 1} {P_{j\beta } \overline a _{nn,\beta } } ]}  + O(\phi ) \notag\\
=& \sigma _l (G)\sum\limits_{\gamma  \in B} {\sum\limits_{ij = 1}^n {F^{ij} \overline a _{\gamma \gamma ,ij} } }  + \sigma _l (G)\sum\limits_{ij = 1}^n {F^{ij} \overline a _{nn,ij} }  - 2\sigma _l (G)\sum\limits_{ij = 1}^n {F^{ij} \sum\limits_{\eta  \in G} {\frac{{\overline a _{n\eta ,i} \overline a _{\eta n,j} }}{{\overline a _{\eta \eta } }}} } \notag \\
&+ \sigma _l (G)\sum\limits_{ij = 1}^n {F^{ij} [\sum\limits_{\eta  \in G} {\frac{{\overline a _{\eta \eta ,j} }}{{\overline a _{\eta \eta } }}} \overline a _{nn,i}  + \sum\limits_{\gamma  \in G} {\frac{{\overline a _{\gamma \gamma ,i} }}{{\overline a _{\gamma \gamma } }}} \overline a _{nn,j} ]}  + O(\phi ) \notag\\
=& \sigma _l (G)\sum\limits_{\gamma  \in B} {\sum\limits_{ij = 1}^n {F^{ij} \overline a _{\gamma \gamma ,ij} } }  + \sigma _l (G)\sum\limits_{ij = 1}^n {F^{ij} \overline a _{nn,ij} }  - 2\sigma _l (G)\sum\limits_{ij = 1}^n {F^{ij} \sum\limits_{\eta  \in G} {\frac{{\overline a _{n\eta ,i} \overline a _{\eta n,j} }}{{\overline a _{\eta \eta } }}} }+ O(\phi  + |\nabla _x \phi |).
\end{align}

For $\gamma \in B$, we have
\begin{align}
\bar a_{\gamma \gamma, ij}  =&-\frac{|u_{n+1}|}{|D u|{u_{n+1}}^3} { \bar A_{\gamma \gamma, ij} }  =-\frac{|u_{n+1}|}{|D u|{u_{n+1}}^3}{ \bar h_{\gamma \gamma,ij} }= -\frac{|u_{n+1}|}{|D u|{u_{n+1}}^3}[ u_{y_{n+1}}^2  u_{y_\gamma y_{\gamma} x_i x_j} ], \notag
\end{align}
so
\begin{align}\label{4.61}
 \sigma _l (G)\sum\limits_{\gamma  \in B} {\sum\limits_{ij = 1}^n {F^{ij} \overline a _{\gamma \gamma ,ij} } } = -\frac{|u_{n+1}|}{|D u|{u_{n+1}}^3}\sigma _l (G)\sum\limits_{\gamma  \in B} {u_{y_{n+1}}^2 \sum\limits_{ij = 1}^n {F^{ij} u_{y_\gamma y_{\gamma} x_i x_j} } } =0.
\end{align}

From \eqref{4.60} and \eqref{4.61}, \eqref{4.54} holds.

\qed

\begin{lemma}\label{lem4.15}
\begin{align}\label{4.62}
\phi _t  = -u_{y_{n+1}}^{-3}\sigma _l (G)[u_{y_{n+1}}^2 u_{y_n y_n t} - 2u_{y_{n+1}} u_{y_n y_{n+1}} u_{y_n t} ] + O(\phi ),
\end{align}
and
\begin{align}\label{4.63}
\sum_{i,j=1}^{n}F^{ij}\phi_{ij}
=&u_{y_{n+1}}^{-3}\sigma_l(G)[-u_{y_{n+1}}^2\sum_{i,j=1}^n F^{ij}u_{ij y_n y_n} + 6u_{y_{n+1}} u_{y_n y_{n+1}}\sum_{i,j=1}^nF^{ij }u_{ij y_n}
-6 u_{y_n y_{n+1}}^2\sum_{i,j=1}^{n}F^{ij}u_{ij}]  \notag  \\
&+2 u_{y_{n+1}}^{-3} \sigma_l(G)\sum_{\a \in G}\sum_{i,j=1}^{n}
F^{ij}\frac{1}{u_{\a \a}}[u_{y_{n+1}} u_{i y_n y_\alpha} - 2 u_{i y_\a} u_{y_n y_{n+1}}][u_{y_{n+1}} u_{j y_n y_\alpha} - 2 u_{j y_\a} u_{y_n y_{n+1}}]  \\
&+ O( \phi + |\nabla_x \phi |).\notag
\end{align}

\end{lemma}

PROOF. Similarly with \eqref{4.44}, taking the first derivatives of $\phi$ in $t$, we have
\begin{align}\label{4.64}
\phi _t  =& \frac{{\partial \phi }}{{\partial t}} = \sum\limits_{\alpha  = 1}^n {\sigma _l (\overline a |\alpha )\overline a _{\alpha \alpha ,t} }  = \sigma _l (G)\overline a _{nn ,t}  + O(\phi ) \notag \\
=&-\frac{|u_{n+1}|}{|D u|{u_{n+1}}^3}\sigma _l (G){ \bar h_{nn,t} }+ O(\phi )\notag \\
=&-\frac{1}{{u_{n+1}}^3}\sigma _l (G)[ u_{y_{n+1}}^2  u_{y_ny_{n} t} - 2 u_{y_{n+1}} u_{y_n t} u_{y_ny_{n+1}}  ]+ O(\phi ).
\end{align}

In the following, we prove \eqref{4.63}. In fact, the calculation is similar as in \cite{BGMX} and \cite{GX}.

It is easy to know
\begin{align}\label{4.65}
\sum\limits_{ij = 1}^n {F^{ij} \overline a _{nn,ij} }  - 2\sum\limits_{ij = 1}^n {F^{ij} \sum\limits_{\eta  \in G} {\frac{{\overline a _{n\eta ,i} \overline a _{\eta n,j} }}{{\overline a _{\eta \eta } }}} } = \sum\limits_{\a \b = 1}^{n+1} {G^{\a \b}  \overline a _{nn, \a \b} }  - 2\sum\limits_{\a \b = 1}^{n+1} {G^{\a \b} \sum\limits_{\eta  \in G} {\frac{{\overline a _{n\eta ,\a} \overline a _{\eta n,\b} }}{{\overline a _{\eta \eta } }}} },
\end{align}
where
\begin{align}\label{4.66}
G^{\a \b} = \sum\limits_{ij = 1}^n F^{ij} P_{i \a} P_{j \b}.
\end{align}
By \eqref{2.13}-\eqref{2.15}, we have
\begin{align}\label{4.67}
\bar a_{nn, \a }  =&(-\frac{|u_{n+1}|}{|D u|{u_{n+1}}^3} )_\a  { \bar A_{nn} }-\frac{|u_{n+1}|}{|D u|{u_{n+1}}^3} { \bar A_{nn, \a } }  = O(\phi )-\frac{|u_{n+1}|}{|D u|{u_{n+1}}^3} { \bar A_{nn, \a } }  \notag \\
=&-\frac{|u_{n+1}|}{|D u|{u_{n+1}}^3}{ \bar h_{nn,\a } }+ O(\phi )\notag \\
=&-\frac{1}{{u_{n+1}}^3}[ u_{y_{n+1}}^2  u_{y_ny_{n} y_\a } - 2 u_{y_{n+1}} u_{y_n y_\a } u_{y_ny_{n+1}}  ]+ O(\phi ),
\end{align}
and
\begin{eqnarray}\label{4.68}
\bar a_{nn,\alpha \beta}&=&(-\frac{|u_{n+1}|}{|D u|{u_{n+1}}^3} )_{\a \b} { \bar A_{nn} } + (-\frac{|u_{n+1}|}{|D u|{u_{n+1}}^3} )_{\a} { \bar A_{nn, \b} } + (-\frac{|u_{n+1}|}{|D u|{u_{n+1}}^3} )_{\b} { \bar A_{nn, \a} }-\frac{|u_{n+1}|}{|D u|{u_{n+1}}^3} { \bar A_{nn, \a \b} } \notag  \\
&=& O(\phi )+ (-\frac{|u_{n+1}|}{|D u|{u_{n+1}}^3} )_{\a} { \bar A_{nn, \b} } + (-\frac{|u_{n+1}|}{|D u|{u_{n+1}}^3} )_{\b} { \bar A_{nn, \a} }-\frac{|u_{n+1}|}{|D u|{u_{n+1}}^3} { \bar A_{nn, \a \b} }  \notag \\
&=& (-\frac{|u_{n+1}|}{|D u|{u_{n+1}}^3} )_{\a} { \bar A_{nn, \b} } + (-\frac{|u_{n+1}|}{|D u|{u_{n+1}}^3} )_{\b} { \bar A_{nn, \a} }-\frac{1}{{u_{n+1}}^3} { \bar h_{nn, \a \b} }+ O(\phi ).
\end{eqnarray}
and
\begin{align}\label{4.69}
\bar h_{nn, \alpha \b} =& {u_{n+1}}^2 u_{y_n y_n \alpha \b} +2 u_{n+1}u_{y_{n+1} y_\a} u_{y_n y_n y_\b}+2 u_{n+1}u_{y_{n+1} y_\b} u_{y_n y_n y_\a}+2 u_{y_{n+1}y_{n+1} }u_{y_n  y_\a} u_{y_n y_\b}  \notag \\
&-2 u_{y_{n+1} }u_{y_n  y_\a y_\b} u_{y_n y_{n+1}}- 2 u_{y_{n+1} }u_{y_n  y_\a} u_{y_n y_{n+1}y_\b}- 2 u_{y_{n+1} }u_{y_n y_\b } u_{y_n y_{n+1}y_\a} \notag \\
& - 2 u_{y_{n+1} y_\b}u_{y_n  y_\a} u_{y_n y_{n+1}}- 2 u_{y_{n+1}y_\a }u_{y_n y_\b } u_{y_n y_{n+1}}\notag \\
=& {u_{n+1}}^2 u_{y_n y_n \alpha \b} + 2 u_{y_{n+1}y_\a }u_{y_n y_\b } u_{y_n y_{n+1}}+2 u_{y_{n+1} y_\b}u_{y_n  y_\a} u_{y_n y_{n+1}}+2 u_{y_{n+1}y_{n+1} }u_{y_n  y_\a} u_{y_n y_\b}  \notag \\
&-2 u_{y_{n+1} }u_{y_n  y_\a y_\b} u_{y_n y_{n+1}}- 2 u_{y_{n+1} }u_{y_n  y_\a} u_{y_n y_{n+1}y_\b}- 2 u_{y_{n+1} }u_{y_n y_\b } u_{y_n y_{n+1}y_\a} \notag \\
&+2 u_{y_{n+1} y_\a} [- u_{y_{n+1}}^2 \bar a_{nn,\beta}] + 2 u_{y_{n+1} y_\b}[ - u_{y_{n+1}}^2 \bar a_{nn,\alpha }]+ O(\phi )
\end{align}
Hence,
\begin{align}\label{4.70}
\sum\limits_{\alpha \b = 1}^{n+1} {G^{\alpha \b}  \overline a _{nn, \alpha \b} } =&\sum\limits_{\alpha \b = 1}^{n+1} G^{\alpha \b}  [-\frac{1}{{u_{n+1}}^3} { \bar h_{nn, \alpha \b} }]+ O(\phi +|\nabla_x \phi| )  \notag\\
=&-\frac{1}{{u_{n+1}}^3} \sum\limits_{\alpha \b = 1}^{n+1} G^{\alpha \b}  [{ u_{n+1}}^2 u_{y_n y_n \alpha \b} + 4 u_{y_{n+1}y_\a }u_{y_n y_\b } u_{y_n y_{n+1}}+2 u_{y_{n+1}y_{n+1} }u_{y_n  y_\a} u_{y_n y_\b} \notag \\
& \qquad \qquad\qquad \quad -2 u_{y_{n+1} }u_{y_n  y_\a y_\b} u_{y_n y_{n+1}}- 4 u_{y_{n+1} }u_{y_n  y_\a} u_{y_n y_{n+1}y_\b}]+O(\phi +|\nabla_x \phi| ).
\end{align}
where
\begin{align}\label{4.71}
\sum_{\alpha \b=1}^{n+1}G^{\alpha\beta}u_{y_{n+1}y_\a }u_{y_n y_\b } u_{y_n y_{n+1}}=& u_{y_n y_{n+1}}^2\sum_{\alpha =1}^{n+1}
G^{\alpha n+1}u_{y_{n+1}y_\a }+O(\phi )  \notag \\
=&u_{y_n y_{n+1}}^2(\sum_{\alpha,\beta=1}^{n+1}- \sum_{\b
=1}^{n}\sum_{\a=1}^{n+1})G^{\alpha\beta}u_{\a\b} +O(\phi ),
\end{align}
\begin{align}\label{4.72}
\sum_{\alpha,\beta=1}^{n+1} G^{\alpha\beta} u_{y_{n+1} }u_{y_n  y_\a} u_{y_n y_{n+1}y_\b}
=&u_{y_{n+1} }u_{y_n y_{n+1}}\sum_{\beta=1}^{n+1} G^{n+1 \beta}  u_{y_n y_{n+1}y_\b}+O(\phi ) \notag \\
=&u_{y_{n+1} }u_{y_n y_{n+1}}(\sum_{\alpha,\beta=1}^{n+1}
-\sum_{\a=1}^{n}\sum_{\b=1}^{n+1} )G^{\alpha\beta}u_{y_n \a \b}+O(\phi ),
\end{align}
and
\begin{eqnarray}\label{4.73}
&&\sum_{\alpha,\beta=1}^{n+1} G^{\alpha\beta} u_{y_{n+1}y_{n+1} }u_{y_n  y_\a} u_{y_n y_\b} = u_{y_{n+1}y_{n+1} } G^{n+1 n+1} u_{y_n y_{n+1}}^2+O(\phi ) \notag \\
&=&u_{y_n y_{n+1}}^2 \sum_{\alpha,\beta=1}^{n+1}
G^{\alpha\beta}u_{y_\a y_\b} - 2u_{y_n y_{n+1}}^2 \sum_{\a=1}^{n} G^{\a
n+1}u_{y_{n+1} y_\a} - u_{y_n y_{n+1}}^2 \sum_{\alpha,\beta=1}^{n}
G^{\alpha\beta}u_{y_\a  y_\b}+O(\phi ).
\end{eqnarray}
So
\begin{eqnarray}\label{4.74}
{u_{n+1}}^3\sum\limits_{\alpha \b = 1}^{n+1} {G^{\alpha \b}  \overline a _{nn, \alpha \b} }
&=&-{ u_{n+1}}^2\sum\limits_{\alpha \b = 1}^{n+1} G^{\alpha \b}   u_{y_n y_n \alpha \b}  + 6 u_{y_{n+1} }u_{y_n y_{n+1}}\sum_{\alpha,\beta=1}^{n+1}
G^{\alpha\beta}u_{y_n \a \b}\nonumber\\
&&-6u_{y_n y_{n+1}}^2 \sum_{\alpha,\beta=1}^{n+1} G^{\alpha\beta}u_{\a\b} -4 u_{y_{n+1} }u_{y_n y_{n+1}} \sum_{\a=1}^{n}\sum_{\b=1}^{n+1} G^{\alpha\beta}u_{y_n \a \b}   \nonumber\\
&&+8 u_{y_n y_{n+1}}^2 \sum_{\a=1}^{n} G^{\a
n+1}u_{y_{n+1} y_\a}  +6 u_{y_n y_{n+1}}^2 \sum_{\alpha,\beta=1}^{n}
G^{\alpha\beta}u_{y_\a  y_\b} +O(\phi ).
\end{eqnarray}
and
\begin{eqnarray}\label{4.75}
&& u_{y_{n+1} }u_{y_n y_{n+1}} \sum_{\a=1}^{n}\sum_{\b=1}^{n+1} G^{\alpha\beta}u_{y_n \a \b} \nonumber\\
&=&u_{y_{n+1} }u_{y_n y_{n+1}} \sum_{\b=1}^{n+1} \bigg(\sum_{\a \in
B} G^{\alpha\beta}u_{y_n \a \b} +\sum_{\a \in G} G^{\alpha\beta}u_{y_n \a \b} \bigg) \nonumber\\
&=& u_{y_{n+1} }u_{y_n y_{n+1}} \sum_{\b=1}^{n+1} \sum_{\a \in G} G^{\alpha\beta}u_{y_n \a \b}\nonumber\\
&=& u_{y_n y_{n+1}} \sum_{\b=1}^{n+1} \sum_{\a \in G} G^{\alpha\beta}[-u_{y_{n+1} }^2 \bar a_{n \a,\b} + u_{y_\a y_\b}u_{y_n y_{n+1}} + u_{y_n y_\b} u_{y_\a y_{n+1}}]+O(\phi )\nonumber\\
&=& - u_{y_{n+1} }^2u_{y_n y_{n+1}} \sum_{\b=1}^{n+1} \sum_{\a \in G} G^{\alpha\beta} \bar a_{n \a,\b}+ u_{y_n y_{n+1}}^2 \sum_{\a \in G} G^{\alpha\a} u_{y_\a y_\a}+ 2 u_{y_n y_{n+1}}^2 \sum_{\a \in G} G^{\alpha n+1} u_{y_\a y_{n+1}}+O(\phi ).
\end{eqnarray}
\eqref{4.74} and \eqref{4.75} yield
\begin{eqnarray}\label{4.76}
&&{u_{n+1}}^3\sum\limits_{\alpha \b = 1}^{n+1} {G^{\alpha \b}  \overline a _{nn, \alpha \b} } \nonumber\\
&=&-{ u_{n+1}}^2\sum\limits_{\alpha \b = 1}^{n+1} G^{\alpha \b}   u_{y_n y_n \alpha \b}  + 6 u_{y_{n+1} }u_{y_n y_{n+1}}\sum_{\alpha,\beta=1}^{n+1}
G^{\alpha\beta}u_{y_n \a \b}-6u_{y_n y_{n+1}}^2 \sum_{\alpha,\beta=1}^{n+1} G^{\alpha\beta}u_{\a\b}   \nonumber\\
&& +4 u_{y_{n+1} }^2u_{y_n y_{n+1}} \sum_{\b=1}^{n+1} \sum_{\a \in G} G^{\alpha\beta} \bar a_{n \a,\b} +2 u_{y_n y_{n+1}}^2 \sum_{\alpha \in G}
G^{\alpha\a}u_{y_\a  y_\a} +O(\phi ).
\end{eqnarray}
So,
\begin{eqnarray}\label{4.77}
&&\sum\limits_{\a \b = 1}^{n+1} {G^{\a \b}  \overline a _{nn, \a \b} }  - 2\sum\limits_{\a \b = 1}^{n+1} {G^{\a \b} \sum\limits_{\eta  \in G} {\frac{{\overline a _{n\eta ,\a} \overline a _{\eta n,\b} }}{{\overline a _{\eta \eta } }}} } \notag  \\
&=&\frac{1}{ u_{y_{n+1} }^3} \left[-{ u_{n+1}}^2\sum\limits_{\alpha \b = 1}^{n+1} G^{\alpha \b}   u_{y_n y_n \alpha \b}  + 6 u_{y_{n+1} }u_{y_n y_{n+1}}\sum_{\alpha,\beta=1}^{n+1}G^{\alpha\beta}u_{y_n \a \b}-6u_{y_n y_{n+1}}^2 \sum_{\alpha,\beta=1}^{n+1} G^{\alpha\beta}u_{\a\b} \right]  \notag  \\
&&-2\sum_{\eta \in G}\left[ \sum\limits_{\a \b = 1}^{n+1} {G^{\a \b}  {\frac{{\overline a _{n\eta ,\a} \overline a _{\eta n,\b} }}{{\overline a _{\eta \eta } }}} } -2 \frac{ u_{y_n y_{n+1}} }{ u_{y_{n+1} }} \sum_{\b=1}^{n+1}  G^{\eta \beta} \bar a_{n \eta,\b} -  \frac{ u_{y_n y_{n+1}}^2 }{ u_{y_{n+1} }^3} G^{\eta \eta} u_{y_\eta  y_\eta}  \right] +O(\phi ).
\end{eqnarray}

In fact, for any $\eta \in G$,
\begin{eqnarray}
&&\sum\limits_{\a \b = 1}^{n+1} {G^{\a \b}  {\frac{{\overline a _{n\eta ,\a} \overline a _{\eta n,\b} }}{{\overline a _{\eta \eta } }}} } -2 \frac{ u_{y_n y_{n+1}} }{ u_{y_{n+1} }} \sum_{\b=1}^{n+1}  G^{\eta \beta} \bar a_{n \eta,\b} -  \frac{ u_{y_n y_{n+1}}^2 }{ u_{y_{n+1} }^3} G^{\eta \eta} u_{y_\eta  y_\eta}  \notag \\
&=&-\frac{ 1 }{ u_{y_{n+1} }^3}[\sum\limits_{\a \b = 1}^{n+1} {G^{\a \b}  {\frac{{\overline h _{n\eta ,\a} \overline h _{\eta n,\b} }}{{\overline h _{\eta \eta } }}} } -2 \frac{ u_{y_n y_{n+1}} }{ u_{y_{n+1} }} \sum_{\b=1}^{n+1}  G^{\eta \beta} \bar h_{n \eta,\b} -  u_{y_n y_{n+1}}^2 G^{\eta \eta} u_{y_\eta  y_\eta} ]  \notag \\
&=&-\frac{ 1 }{ u_{y_{n+1} }^3}\left\{\sum_{\alpha,\beta=1}^{n}
G^{\alpha\beta}\frac{1}{u_{y_{n+1} }^2u_{\eta \eta}}[ u_{y_{n+1} }^2u_{n\eta\alpha}- u_{y_{n+1} } u_{\eta\a}u_{y_n y_{n+1}}][ u_{y_{n+1} }^2u_{n\eta\b}- u_{y_{n+1} }u_{\eta\b}u_{y_n y_{n+1}}]\right. +O(\phi )\notag  \\
&& \qquad \qquad  +2\sum_{\alpha=1}^{n}
G^{\alpha n+1}\frac{1}{u_{y_{n+1} }^2u_{\eta \eta}}[ u_{y_{n+1} }^2u_{n\eta\alpha}- u_{y_{n+1} }u_{\eta\a}u_{y_n y_{n+1}}][ u_{y_{n+1} }^2u_{n\eta n+1}-2 u_{y_{n+1} }u_{\eta n+1}u_{y_n y_{n+1}}]  +O(\phi )\notag  \\
&& \qquad \qquad + G^{n+1 n+1}\frac{1}{u_{y_{n+1} }^2u_{\eta \eta}} [ u_{y_{n+1} }^2u_{n\eta n+1}-2 u_{y_{n+1} }u_{\eta n+1}u_{y_n y_{n+1}}][ u_{y_{n+1} }^2u_{n\eta n+1}-2 u_{y_{n+1} }u_{\eta n+1}u_{y_n y_{n+1}}]  \notag  \\
&& \qquad \qquad -2\sum_{\a=1}^{n} G^{\alpha \eta}\frac{1}{u_{y_{n+1} }^2u_{\eta \eta}}[ u_{y_{n+1} }^2u_{n\eta\a}- u_{y_{n+1} }u_{\eta \a}u_{y_n y_{n+1}}][ u_{y_{n+1} }u_{\eta \eta}u_{y_n y_{n+1}}]  +O(\phi )\notag \\
&& \qquad \qquad -2 G^{n+1 \eta}\frac{1}{u_{y_{n+1} }^2u_{\eta \eta}}[ u_{y_{n+1} }^2u_{n\eta n+1}-2 u_{y_{n+1} }u_{\eta n+1}u_{y_n y_{n+1}}][ u_{y_{n+1} }u_{\eta \eta}u_{y_n y_{n+1}}]  \notag \\
&& \qquad \qquad
\left.+ G^{\eta \eta}\frac{1}{u_{y_{n+1} }^2u_{\eta \eta}} ( u_{y_{n+1} }u_{\eta \eta}u_{y_n y_{n+1}})^2\right\}
\notag \\
&=&-\frac{ 1 }{ u_{y_{n+1} }^3}\sum_{\alpha,\beta=1}^{n+1}
G^{\alpha\beta}\frac{1}{u_{\eta \eta}}[ u_{y_{n+1} }u_{n\eta\alpha}- 2 u_{\eta\a}u_{y_n y_{n+1}}][ u_{y_{n+1} }u_{n\eta\b}- 2 u_{\eta\b}u_{y_n y_{n+1}}]+O(\phi ). \notag
\end{eqnarray}
Then we get
\begin{align}\label{4.78}
&\sum\limits_{\a \b = 1}^{n+1} {G^{\a \b}  \overline a _{nn, \a \b} }  - 2\sum\limits_{\a \b = 1}^{n+1} {G^{\a \b} \sum\limits_{\eta  \in G} {\frac{{\overline a _{n\eta ,\a} \overline a _{\eta n,\b} }}{{\overline a _{\eta \eta } }}} }  \notag  \\
=&u_{y_{n+1}}^{-3}[-u_{y_{n+1}}^2\sum_{i,j=1}^n F^{ij}u_{ij y_n y_n} + 6u_{y_{n+1}} u_{y_n y_{n+1}}\sum_{i,j=1}^nF^{ij }u_{ij y_n}
-6 u_{y_n y_{n+1}}^2\sum_{i,j=1}^{n}F^{ij}u_{ij}]  \notag  \\
&+2 u_{y_{n+1}}^{-3} \sum_{\a \in G}\sum_{i,j=1}^{n}
F^{ij}\frac{1}{u_{\a \a}}[u_{y_{n+1}} u_{i y_n y_\alpha} - 2 u_{i y_\a} u_{y_n y_{n+1}}][u_{y_{n+1}} u_{j y_n y_\alpha} - 2 u_{j y_\a} u_{y_n y_{n+1}}] + O( \phi + |\nabla_x \phi |).
\end{align}
Hence \eqref{4.63} holds from \eqref{4.65} and \eqref{4.78}.
\qed

\begin{lemma}\label{lem4.16}

\begin{align}\label{4.79}
&\sum_{\a \in G}\sum_{i,j=1}^{n}
F^{ij}\frac{1}{u_{y_\a y_\a}}[u_{y_{n+1}} u_{i y_n y_\alpha} - 2 u_{i y_\a} u_{y_n y_{n+1}}][u_{y_{n+1}} u_{j y_n y_\alpha} - 2 u_{j y_\a} u_{y_n y_{n+1}}] \notag \\
\leq& \sum_{k \in G}\sum_{i,j=1}^{n}
F^{ij}\frac{1}{u_{x_k x_k}}[u_{y_{n+1}} u_{i y_n x_k} - 2 u_{i x_k} u_{y_n y_{n+1}}][u_{y_{n+1}} u_{j y_n x_k} - 2 u_{j x_k} u_{y_n y_{n+1}}]+ O(\phi  + |\nabla _x \phi |).
\end{align}

\end{lemma}

PROOF. First, we consider a special case: $ F^{ij}= \delta_{ij}$.
That is, we need to prove
\begin{align}\label{4.80}
\sum_{\a \in G}\sum_{i=1}^{n}\frac{1}{u_{y_\a y_\a}}[u_{y_{n+1}} u_{i y_n y_\alpha} - 2 u_{i y_\a} u_{y_n y_{n+1}}]^2
\leq \sum_{k \in G}\sum_{i=1}^{n}\frac{1}{u_{x_k x_k}}[u_{y_{n+1}} u_{i y_n x_k} - 2 u_{i x_k} u_{y_n y_{n+1}}]^2+ O(\phi  + |\nabla _x \phi |).
\end{align}
Form \eqref{4.35}-\eqref{4.37}, \eqref{4.40} and \eqref{4.41}, we have
\begin{align}
\label{4.81}& u_{y_{n+1}} u_{i y_n y_\alpha} - 2 u_{i y_\a} u_{y_n y_{n+1}} = 0, \quad \a \in G, i \in B;  \\
\label{4.82}& u_{y_{n+1}} u_{i y_n y_\alpha} - 2 u_{i y_\a} u_{y_n y_{n+1}} = 0,  \quad \a \in B, i \in G\cup \{ n \}; \\
\label{4.83}& u_{y_{n+1}} u_{i y_n y_n} - 2 u_{i y_n} u_{y_n y_{n+1}} = O(\phi  + |\nabla _x \phi |), \quad i \in G\cup \{ n \}.
\end{align}
Since $\nabla_y^2 u = \bigg(u_{y_i y_j}\bigg)_{1 \leq i,j \leq n} \leq 0$ is diagonal, by the approximation, we have for $i \in G \cup \{ n \}$
\begin{align}\label{4.84}
&\sum_{\a \in G}\frac{1}{u_{y_\a y_\a}}[u_{y_{n+1}} u_{i y_n y_\alpha} - 2 u_{i y_\a} u_{y_n y_{n+1}}]^2  \notag \\
=& \mathop {\lim }\limits_{\varepsilon  \to 0 + } (u_{y_{n+1}} \nabla_y u_{i y_n} - 2 \nabla_y u_{i} u_{y_n y_{n+1}})\left( {\nabla_y^2 u - \varepsilon I_n} \right)^{ - 1} (u_{y_{n+1}} \nabla_y u_{i y_n} - 2 \nabla_y u_{i} u_{y_n y_{n+1}})^T  + O(\phi  + |\nabla _x \phi |),
 \end{align}
where $\varepsilon >0$ small, and
\begin{align}\label{4.85}
&(u_{y_{n+1}} \nabla_y u_{i y_n} - 2 \nabla_y u_{i} u_{y_n y_{n+1}})\left( {\nabla_y^2 u - \varepsilon I_n} \right)^{ - 1} (u_{y_{n+1}} \nabla_y u_{i y_n} - 2 \nabla_y u_{i} u_{y_n y_{n+1}})^T  \notag \\
=& (u_{y_{n+1}} \nabla_z u_{i y_n} - 2 \nabla_z u_{i} u_{y_n y_{n+1}})T^T\left( {\nabla_y^2 u - \varepsilon I_n} \right)^{ - 1} T (u_{y_{n+1}} \nabla_z u_{i y_n} - 2 \nabla_z u_{i} u_{y_n y_{n+1}})^T   \notag \\
=&  (u_{y_{n+1}} \nabla_z u_{i y_n} - 2 \nabla_z u_{i} u_{y_n y_{n+1}})\left( {\nabla_z^2 u - \varepsilon I_n} \right)^{ - 1}  (u_{y_{n+1}} \nabla_z u_{i y_n} - 2 \nabla_z u_{i} u_{y_n y_{n+1}})^T .
\end{align}

Denote
\begin{align}\label{4.86}
 C := u_{z_n z_n}  - \varepsilon  - \sum\limits_{i = 1}^l {\frac{{u_{z_i z_n} ^2 }}{{u_{z_i z_i }  - \varepsilon }}}  < 0,
 \end{align}
then
\begin{align}
\left( {\nabla_{z}^2 u - \varepsilon I_n} \right)^{ - 1}  =& diag(\frac{1}{{u_{z_1 z_1 }  - \varepsilon }}, \cdots ,\frac{1}{{u_{z_l z_l }  - \varepsilon }}, - \frac{1}{\varepsilon }, \cdots , -\frac{1}{\varepsilon },0)  \notag \\
&+ \frac{1}{C}( - \frac{{u_{z_1 z_n} }}{{u_{z_1 z_1 }  - \varepsilon }}, \cdots , - \frac{{u_{z_l z_n} }}{{u_{z_l z_l }  - \varepsilon }},0, \cdots ,0,1)^T ( - \frac{{u_{z_1 z_n} }}{{u_{z_1 z_1 }  - \varepsilon }}, \cdots , - \frac{{u_{z_l z_n} }}{{u_{z_l z_l }  - \varepsilon }},0, \cdots ,0,1)  \notag \\
\le& diag(\frac{1}{{u_{z_1 z_1 }  - \varepsilon }}, \cdots ,\frac{1}{{u_{z_l z_l }  - \varepsilon }}, 0, \cdots ,0,0). \notag
\end{align}
So
\begin{align}\label{4.87}
&(u_{y_{n+1}} \nabla_z u_{i y_n} - 2 \nabla_z u_{i} u_{y_n y_{n+1}})\left( {\nabla_z^2 u - \varepsilon I_n} \right)^{ - 1}  (u_{y_{n+1}} \nabla_z u_{i y_n} - 2 \nabla_z u_{i} u_{y_n y_{n+1}})^T  \notag \\
\le& \sum_{k \in G}\frac{1}{u_{z_k z_k} - \varepsilon }[u_{y_{n+1}} u_{i y_n z_k} - 2 u_{i z_k} u_{y_n y_{n+1}}]^2
=\sum_{k \in G}\frac{1}{u_{x_k x_k} - \varepsilon }[u_{y_{n+1}} u_{i y_n x_k} - 2 u_{i x_k} u_{y_n y_{n+1}}]^2.
\end{align}
Then we have for $i \in G\cup \{ n \}$
\begin{align}\label{4.88}
&\sum_{\a \in G}\frac{1}{u_{y_\a y_\a}}[u_{y_{n+1}} u_{i y_n y_\alpha} - 2 u_{i y_\a} u_{y_n y_{n+1}}]^2  \notag \\
\le& \mathop {\lim }\limits_{\varepsilon  \to 0 + } \sum_{k \in G}\frac{1}{u_{x_k x_k} - \varepsilon }[u_{y_{n+1}} u_{i y_n x_k} - 2 u_{i x_k} u_{y_n y_{n+1}}]^2 + O(\phi  + |\nabla _x \phi |)  \notag  \\
=& \sum_{k \in G}\frac{1}{u_{x_k x_k} }[u_{y_{n+1}} u_{i y_n x_k} - 2 u_{i x_k} u_{y_n y_{n+1}}]^2 + O(\phi  + |\nabla _x \phi |).
\end{align}
Hence, \eqref{4.80} holds from \eqref{4.81} and \eqref{4.88}.

For the general case, the CLAIM also holds following the above proof.
\qed

\begin{theorem} \label{th4.17}
Under the assumptions of Theorem 1.2 and the above notations, we
have for any fixed point  $(x,t) \in \mathcal {O}\times (t_0-\delta, t_0)$,
\begin{equation}\label{4.89}
\sum_{ ij=1 } ^n F^{ij} \phi_{ ij} -\phi_t  \le c_0( \phi+ |\nabla_x \phi| )
\end{equation}

\end{theorem}

PROOF. From \eqref{4.62}, \eqref{4.63} and \eqref{4.79},
\begin{align}\label{4.90}
\sum_{i,j=1}^{n}F^{ij}\phi_{ij} -\phi _t
\leq& u_{y_{n+1}}^{-3}\sigma_l(G)\Big[-u_{y_{n+1}}^2(\sum_{i,j=1}^n F^{ij}u_{ij y_n y_n} - u_{y_n y_n t}) -2 u_{y_{n+1}} u_{y_n y_{n+1}} u_{y_n t} \notag \\
&\qquad \qquad\quad + 6u_{y_{n+1}} u_{y_n y_{n+1}}\sum_{i,j=1}^nF^{ij }u_{ij y_n}-6 u_{y_n y_{n+1}}^2\sum_{i,j=1}^{n}F^{ij}u_{ij}\Big]  \notag  \\
&+2 u_{y_{n+1}}^{-3} \sigma_l(G)\sum_{k \in G}\sum_{i,j=1}^{n}
F^{ij}\frac{1}{u_{x_k x_k}}[u_{y_{n+1}} u_{i y_n x_k} - 2 u_{i x_k} u_{y_n y_{n+1}}][u_{y_{n+1}} u_{j y_n x_k} - 2 u_{j x_k} u_{y_n y_{n+1}}]  \notag  \\
&+ O( \phi + |\nabla_x \phi |).
\end{align}

From the equation \eqref{1.1}, we get
\begin{align}\label{4.91}
u_{y_n y_n t}  =& \sum\limits_{ij = 1}^n {F^{ij} } u_{y_n y_n ij}  + \sum\limits_{k = 1}^n {F^{u_k } } u_{ky_n y_n }  + F^{u,u} u_{y_n y_n }  \notag \\
&+ \sum\limits_{ijkl = 1}^n {F^{ij,kl} u_{ijy_n } u_{kly_n } }  + 2\sum\limits_{ijk = 1}^n {F^{ij,u_k } } u_{ijy_n } u_{ky_n }  + 2\sum\limits_{ij = 1}^n {F^{ij,u} } u_{ijy_n } u_{y_n }  \notag  \\
& + 2\sum\limits_{ijk = 1}^n {F^{ij,x_k } } u_{ijy_n } \frac{{\partial x_k }}{{\partial y_n }}+ 2\sum\limits_{ij = 1}^n {F^{ij,t} } u_{ijy_n } \frac{{\partial t}}{{\partial y_n }} + \sum\limits_{kl = 1}^n {F^{u_k ,u_l } } u_{ky_n } u_{ly_n }  + 2\sum\limits_{k = 1}^n {F^{u_k ,u} } u_{ky_n } u_{y_n }   \notag \\
&+ 2\sum\limits_{kl = 1}^n {F^{u_k ,x_l } } u_{ky_n } \frac{{\partial x_l }}{{\partial y_n }} + 2\sum\limits_{k = 1}^n {F^{u_k ,t} } u_{ky_n } \frac{{\partial t}}{{\partial y_n }} + F^{u,u} u_{y_n } ^2  + 2\sum\limits_{k = 1}^n {F^{u,x_k } } u_{y_n } \frac{{\partial x_k }}{{\partial y_n }}\notag \\
& + 2F^{u,t} u_{y_n } \frac{{\partial t}}{{\partial y_n }}  + \sum\limits_{ik = 1}^n {F^{x_i ,x_k } } \frac{{\partial x_i }}{{\partial y_n }}\frac{{\partial x_k }}{{\partial y_n }} + 2\sum\limits_{i = 1}^n {F^{x_i ,t} } \frac{{\partial x_i }}{{\partial y_n }}\frac{{\partial t}}{{\partial y_n }} + F^{t,t} \left( {\frac{{\partial t}}{{\partial y_n }}} \right)^2   \notag \\
=& \sum\limits_{ij = 1}^n {F^{ij} } u_{y_n y_n ij}  + F^{u_n } u_{x_n y_n y_n }+ O( \phi + |\nabla_x \phi |)  \notag  \\
&+ \sum\limits_{ijkl = 1}^n {F^{ij,kl} u_{ijy_n } u_{kly_n } }  + 2\sum\limits_{ij = 1}^n {F^{ij,u_n } } u_{ijy_n } u_{x_n y_n }  + 2\sum\limits_{ijk = 1}^n {F^{ij,x_k } } u_{ijy_n } \frac{{\partial x_k }}{{\partial y_n }}  \notag \\
&+ 2\sum\limits_{ij = 1}^n {F^{ij,t} } u_{ijy_n } \frac{{\partial t}}{{\partial y_n }} + F^{u_n ,u_n } u_{x_n y_n } ^2  + 2\sum\limits_{l = 1}^n {F^{u_n ,x_l } } u_{x_n y_n } \frac{{\partial x_l }}{{\partial y_n }} \\
&+ 2F^{u_n ,t} u_{x_n y_n } \frac{{\partial t}}{{\partial y_n }} + \sum\limits_{ik = 1}^n {F^{x_i ,x_k } } \frac{{\partial x_i }}{{\partial y_n }}\frac{{\partial x_k }}{{\partial y_n }} + 2\sum\limits_{i = 1}^n {F^{x_i ,t} } \frac{{\partial x_i }}{{\partial y_n }}\frac{{\partial t}}{{\partial y_n }} + F^{t,t} \left( {\frac{{\partial t}}{{\partial y_n }}} \right)^2, \notag
\end{align}
And from \eqref{4.49},
\begin{align}\label{4.92}
&u_{x_n} u_{y_ny_{n+1}}= u_{y_{n+1}} u_{x_n y_n}+ O( \phi ),
\end{align}
so
\begin{align}\label{4.93}
u_{x_n y_n y_n }  = 2\frac{1}{{u_{y_{n + 1} } }}u_{x_n y_n } u_{y_n y_{n + 1} } + O( \phi ) = 2\frac{{u_{y_n y_{n + 1} } }}{{u_{y_{n + 1} } }}u_{x_n y_n } + O( \phi ) = 2\frac{{u_{x_n y_n } }}{{u_{x_n } }}u_{x_n y_n }+ O( \phi ).
\end{align}
Hence
\begin{align}\label{4.94}
\sum_{i,j=1}^{n}F^{ij}\phi_{ij} -\phi _t
\leq& u_{y_{n+1}}^{-3}\sigma_l(G)\Big(Q  -2 u_{y_{n+1}} u_{y_n y_{n+1}} u_{y_n t} \Big)+ O( \phi + |\nabla_x \phi |).
\end{align}
where
\begin{align}\label{4.95}
Q =& \sum\limits_{ijkl = 1}^n {F^{ij,kl} u_{ijy_n } u_{kly_n }u_{y_{n+1}}^2 }  + 2\sum\limits_{ij = 1}^n {F^{ij,u_n } } u_{ijy_n } u_{x_n y_n } u_{y_{n+1}}^2 + 2\sum\limits_{ijk = 1}^n {F^{ij,x_k } } u_{ijy_n } \frac{{\partial x_k }}{{\partial y_n }} u_{y_{n+1}}^2 \notag \\
&+ 2\sum\limits_{ij = 1}^n {F^{ij,t} } u_{ijy_n } \frac{{\partial t}}{{\partial y_n }}u_{y_{n+1}}^2 + F^{u_n ,u_n } u_{x_n y_n } ^2u_{y_{n+1}}^2  + 2\sum\limits_{l = 1}^n {F^{u_n ,x_l } } u_{x_n y_n } \frac{{\partial x_l }}{{\partial y_n }} u_{y_{n+1}}^2\notag  \\
&+ 2F^{u_n ,t} u_{x_n y_n } \frac{{\partial t}}{{\partial y_n }} u_{y_{n+1}}^2+ \sum\limits_{ik = 1}^n {F^{x_i ,x_k } } \frac{{\partial x_i }}{{\partial y_n }}\frac{{\partial x_k }}{{\partial y_n }} u_{y_{n+1}}^2+ 2\sum\limits_{i = 1}^n {F^{x_i ,t} } \frac{{\partial x_i }}{{\partial y_n }}\frac{{\partial t}}{{\partial y_n }}u_{y_{n+1}}^2 + F^{t,t} \left( {\frac{{\partial t}}{{\partial y_n }}} \right)^2 u_{y_{n+1}}^2  \notag \\
&+ 2F^{u_n } \frac{1}{{u_{x_n } }}u_{x_n y_n } ^2u_{y_{n+1}}^2  + 6 \frac{ u_{y_{n+1}}^2 } {u_{x_{n} }} u_{y_n x_{n} } \sum\limits_{ij = 1}^n {F^{ij} u_{ijy_n } }  - 6u_{y_{n+1}}^2 \frac{ u_{y_n x_{n } } ^2}{ u_{x_n}^2 } \sum\limits_{ij = 1}^n {F^{ij} u_{ij} }  \\
&+2\sum_{k \in G}\sum_{i,j=1}^{n}
F^{ij}\frac{1}{u_{x_k x_k}}[u_{y_{n+1}} u_{i y_n x_k} - 2 u_{i x_k} u_{y_n y_{n+1}}][u_{y_{n+1}} u_{j y_n x_k} - 2 u_{j x_k} u_{y_n y_{n+1}}], \notag
\end{align}
From \eqref{4.50}, we have
\begin{align}\label{4.96}
u_{y_{n+1}} u_{y_n y_{n+1}} u_{y_n t} =u_{t}  u_{y_ny_{n+1}}^2 + O(\phi ) \geq  O(\phi ),
\end{align}
Set
\begin{align*}
&s=\frac{1}{u_{x_n}},  A_{ij} = s u_{ij} = \frac{u_{ij}}{u_{x_n}}, \theta = (0, \cdots, 0, 1), \\
&X_{\a \b} =u_{x_{\a} x_{\b} y_n }u_{x_{n}},  \\
&Y = u_{x_n y_n }u_{x_{n}},  \\
&Z_k= \frac{{\partial x_k}}{{\partial y_n }} u_{x_{n}},  \\
&D = \frac{{\partial t}}{{\partial y_n }} u_{x_{n}}, \\
&V = ((X_{\alpha\beta} ),Y,(Z_i), D)\in \mathcal{S}^n \times \mathbb{R} \times \mathbb{R}^{n}\times
\mathbb{R} ;
\end{align*}
then we get
\begin{align*}
X_{i \alpha}=0,\quad A_{i \a} =0, \quad X_{i \alpha} -2 A_{i \a} Y =0,  \quad i \in B.
\end{align*}
So it yields
\begin{align} \label{3.20}
Q = \frac{u_{y_{n+1}}^2}{u_{x_n}^2} Q^*(V, V),
\end{align}
where $Q^*(V, V)$ is defined in \eqref{2.23}.

From the structural condition \eqref{1.3} (i.e. Lemma \ref{lem2.7}),
\begin{equation*}
Q^*( V, V)\le 0.
\end{equation*}
so we get
\begin{align}\label{4.98}
Q=\frac{u_{y_{n+1}}^2}{u_{x_n}^2}  Q^*( V,  V) \leq  0.
\end{align}

Hence \eqref{4.89} holds from \eqref{4.94}, \eqref{4.96} and \eqref{4.98}.
\qed

\subsection{Constant rank theorem of the spacetime fundamental form for the equations \eqref{1.12}-\eqref{1.14}}

Following the proof of Theorem \ref{th1.2}, we establish the constant rank theorem for the spacetime fundamental form for the equations \eqref{1.12}-\eqref{1.14} in this subsection as follows.
\begin{theorem}\label{th4.18}
Suppose $u \in C^{3,1}(\Omega \times [0,T])$ is a spacetime quasiconcave to the parabolic equation \eqref{1.12}-\eqref{1.14}, and satisfies the condition \eqref{1.4}.
Then the second fundamental form of spatial level sets $\hat \Sigma^c =\{(x, t) \in \Omega \times (0, T)| u(x,t) =
c\}$ has the same constant rank in $\Omega$ for any fixed $t \in (0, T]$. Moreover, let $l(t)$ be the minimal
rank of the second fundamental form in $\Omega$, then $l(s)
\leqslant l(t)$ for all $0< s \leqslant t \leqslant T$.
\end{theorem}

PROOF. The proof is following the proof of Theorem \ref{th1.2}.

Suppose $\hat{a}(x,t)=(\hat a_{\a \b})_{n \times n}$ attains the minimal
rank $l$ at some point $(x_0, t_0) \in \Omega \times (0, T]$. We may assume
$l\leqslant n-1$, otherwise there is nothing to prove. At $(x_0,t_0)$, we may choose $e_1,\cdots, e_{n-1},e_n$ such that
\begin{equation}\label{4.3}
 |\n u(x_0,t_0)|=u_n(x_0,t_0)>0\ \mbox{ and } \
\Big(u_{ij}\Big)_{1 \leq i,j\leq n-1} \mbox{is diagonal at}\ (x_0,t_0).
\end{equation}
Without loss of generality we assume $ u_{11} \leq u_{22}\leq \cdots
\leq u_{n-1n-1} $. So, at $(x_0,t_0)$, from \eqref{4.1}, we have the matrix
$\Big(\hat a_{ij} \Big)_{1 \leq i,j \leq n-1}$ is also diagonal, and $\hat a_{11} \geq \hat a_{22} \geq \cdots \geq \hat a_{n-1
n-1}$. From lemma \ref{lem2.5}, there is a positive constant $
C_0$ such that at $(x_0,t_0)$

CASE 1:
\begin{eqnarray*}
&&\hat a_{11}  \geq \cdots \geq \hat a_{l-1l-1} \geq
C_0 , \quad \hat a_{ll} = \cdots = \hat a_{n-1n-1} =0 , \\
&&\hat a_{nn} -\sum\limits_{i = 1}^{l-1} {\frac{{\hat a_{in} ^2 }} {{\hat a_{ii}
}}} \geq C_0 ,  \quad \hat a_{in} = 0, \quad l \leqslant i \leqslant n-1.
\end{eqnarray*}

CASE 2:
\begin{eqnarray*}
&&\hat a_{11}  \geq \cdots \geq \hat a_{ll} \geq C_0 , \quad \hat a_{l+1l+1} =
\cdots = \hat a_{n-1n-1} =0, \\
&& \hat a_{tt}  = \sum\limits_{i = 1}^{l}{\frac{{\hat a_{in} ^2 }} {{\hat a_{ii} }}}
,  \quad \hat a_{in} = 0, \quad  l+1 \leqslant i \leqslant n-1.
\end{eqnarray*}

For the CASE 1, the theorem holds from Subsection 4.1 and Theorem 3.3, Theorem 3.5, Theorem 3.7.

For the CASE 2, we need to prove the differential inequality \eqref{4.89}, which is similar to the proof of Theorem 3.3, Theorem 3.5, and Theorem 3.7,
with some modifications. In the following, we just prove \eqref{4.89} for the equation \eqref{1.12}. And for the equation \eqref{1.13} and \eqref{1.14},
the proofs follow from the proofs of Theorem 3.5 and Theorem 3.7 with the same modifications.

For the equation \eqref{1.12}, following the assumptions and notations, we need to prove
\begin{equation} \label{4.100}
\sum_{\alpha,\beta=1}^{n}L_{\alpha\beta}\phi_{\alpha\beta}(x,t)-\phi_t\leq
C(\phi+|\nabla_x \phi|),~~\forall ~(x,t)\in \mathcal {O}
\times(t_0-\delta_0,t_0),
\end{equation}
where $\phi$ is defined in \eqref{4.25} and $C$ is a positive constant  independent of $\phi$.
Then by a approximation, \eqref{4.100} holds for $t = t_0$.
Then by the strong maximum principle and the method of continuity, we can prove Theorem \ref{th4.18} under CASE 2.

For any fixed point  $(x,t) \in \mathcal {O}\times (t_0-\delta,
t_0]$, following the Subsection 4.2, we first choose the space coordinates $e_1,\cdots, e_{n-1},e_n$, $e_{n+1}$ still the time coordinate with
\begin{align*}
 y = (y_1 , \cdots ,y_n ,y_{n + 1} ) = (x,t)P, \quad P = OT,
\end{align*}
such that
\begin{equation}\label{4.101}
 |\n u(x,t)|=u_n(x,t)>0\ \mbox{and}
\Big(u_{ij} \Big)_{1 \leq i,j \leq n-1} \mbox{ is diagonal at}\ (x,t).
\end{equation}
Finally, we get a new spacetime coordinate $\{\bar e_1, \cdots, \bar e_{l}, e_{l+1},\cdots, e_{n-1}, \bar e_n, \hat e_{n+1}\}$
such that
\begin{align}
\label{4.102}&u_{y_{n+1}} = |D u| > 0,\quad  u_{y_1} = \cdots =  u_{y_n} = 0, \quad \text{ at } (x,t), \\
\label{4.103}& \bigg( u_{y_{\a} y_{\b}} \bigg)_{1 \leq \a, \b\leq n} \mbox{is diagonal at}\ (x,t).
\end{align}

Also we will use $i,j,k,l =1, \cdots, n$ to represent the $x$ coordinates, $t$ still the time coordinate, and
$\alpha, \beta, \gamma, \eta =1, \cdots, n+1$ the $y$ coordinates.

Following the proof of Theorem \ref{th4.17}, we get from \eqref{4.94}
\begin{align}\label{4.104}
\sum_{i,j=1}^{n}L_{ij}\phi_{ij} -\phi _t
\leq& u_{y_{n+1}}^{-3}\sigma_l(G)\Big(Q  -2 u_{y_{n+1}} u_{y_n y_{n+1}} u_{y_n t} \Big)+ O( \phi + |\nabla_x \phi |) \notag \\
=& u_{y_{n+1}}^{-3}\sigma_l(G) \frac{u_{y_{n+1}}^2}{u_{x_n}^2} \Big(Q ^*(V,V) -2 u_{x_{n}} u_{x_n y_n } u_{y_n t} \Big)+ O( \phi + |\nabla_x \phi |).
\end{align}
where
\begin{align*}
Q ^*(V,V) =&  2\sum\limits_{kl = 1}^n \frac{{\partial L_{kl} }}{{\partial u_n }} u_{kly_n} u_{n y_n}u_{n}^2 +  \sum\limits_{kl = 1}^n \frac{{\partial ^2 L_{kl} }}{{\partial u_n ^2 }} u_{kl} u_{n y_n}^2 u_n^2    \\
&  + 2 \sum\limits_{kl = 1}^n \frac{{\partial L_{kl} }}{{\partial u_n }} u_{kl}u_n u_{n y_n}^2 +6u_nu_{n y_n}\sum_{kl=1}^n
L_{kl}u_{kl y_n} - 6u_{n y_n}^2\sum_{kl=1}^n L_{kl}u_{kl}  \\
&+2\sum_{i\in G}\sum_{\alpha,\beta=1}^{n}
\frac{1}{u_{ii}}L_{\alpha\beta}[u_nu_{i\alpha y_n}-2u_{i\a}u_{n y_n}][u_nu_{i\b y_n}-2u_{i\b}u_{n y_n}].
\end{align*}

Under the coordinate \eqref{4.101} ( i.e. \eqref{3.38}), we still have \eqref{a1} - \eqref{a7}.
So from the equation \eqref{1.12}, we know
\begin{align*}
 u_t  = L_{kk} u_{kk},
\end{align*}
and differentiating the equation in $y_n$,
\begin{align}
u_{ty_n}  =& L_{kk} u_{kk y_n}  + \frac{{\partial L_{kl} }}{{\partial u_i }}u_{i y_n} u_{kl} =L_{kk} u_{kk y_n}  + \frac{{\partial L_{kl} }}{{\partial u_n }}u_{n y_n} u_{kl}+O(\phi) \notag \\
=& L_{kk} u_{kk y_n}  + (p - 2)\frac{{L_{kk} }}{{u_n }}u_{n y_n} u_{kk} +O(\phi)\notag \\
=& L_{kk} u_{kk y_n}  + (p - 2)\frac{{u_t }}{{u_n }}u_{n y_n} +O(\phi).
\end{align}

So
\begin{align*}
Q ^*(V,V) =&  2\sum\limits_{k = 1}^n (p - 2)\frac{{L_{kk} }}{{u_n }} u_{kk y_n} u_{n y_n}u_n^2 +  \sum\limits_{k = 1}^n (p - 2)(p - 3)\frac{{L_{kk} }}{{u_n ^2 }} u_{kk} u_{n y_n}^2 u_n^2    \\
&  + 2 \sum\limits_{k = 1}^n (p - 2)\frac{{L_{kk} }}{{u_n }} u_{kk}u_n u_{n y_n}^2 +6u_nu_{n y_n}\sum_{k=1}^n
L_{kk}u_{kk y_n} - 6u_{n y_n}^2 u_t  \\
&+2\sum_{i\in G}\sum_{\alpha=1}^{n}
\frac{1}{u_{ii}}L_{\alpha\a}[u_nu_{i\alpha y_n}-2u_{i\a}u_{n y_n}]^2  +O(\phi)\\
=&  2(p - 2)[u_{t y_n}  - (p - 2)\frac{{u_t }}{{u_n }}u_{n y_n}] u_{n y_n}u_n +  (p - 2)(p - 3) u_{t} u_{n y_n}^2   \\
&  + 2 (p - 2)u_t u_{n y_n}^2 + 6 u_n u_{n y_n}[u_{ t y_n}  - (p - 2)\frac{{u_t }}{{u_n }}u_{n y_n}]- 6u_{n y_n}^2 u_t  \\
&+2\sum_{i\in G}\sum_{\alpha=1}^{n}
\frac{1}{u_{ii}}L_{\alpha\a}[u_nu_{i\alpha y_n}-2u_{i\a}u_{n y_n}]^2  +O(\phi)\\
=&  ( 2p+2) u_n u_{n y_n}u_{t y_n} - (p^2+p)u_t u_{n y_n}^2   \\
&+2\sum_{i\in G}\sum_{\alpha=1}^{n}
\frac{1}{u_{ii}}L_{\alpha\a}[u_nu_{i\alpha y_n}-2u_{i\a}u_{n y_n}]^2 +O(\phi).
\end{align*}
Hence from \eqref{4.49} and \eqref{4.50},
\begin{align*}
Q ^*(V,V) -2u_n u_{n y_n}u_{t y_n} =& 2p u_n u_{ny_n }u_{t y_n} - (p^2+p)u_t u_{n y_n}^2   \\
&+2\sum_{i\in G}\sum_{\alpha=1}^{n}
\frac{1}{u_{ii}}L_{\alpha\a}[u_nu_{i\alpha y_n}-2u_{i\a}u_{n y_n}]^2  +O(\phi) \\
=& 2p  u_{ny_n}[ u_t u_{n y_n}+O(\phi) ]- (p^2+p)u_t u_{jn}^2   \\
&+2\sum_{i\in G}\sum_{\alpha=1}^{n}
\frac{1}{u_{ii}}L_{\alpha\a}[u_nu_{i\alpha y_n}-2u_{i\a}u_{n y_n}]^2   +O(\phi)\\
=& - (p^2-p)u_t u_{jn}^2 +2\sum_{i\in G}\sum_{\alpha=1}^{n}
\frac{1}{u_{ii}}L_{\alpha\a}[u_nu_{i\alpha y_n}-2u_{i\a}u_{n y_n}]^2   +O(\phi)\\
\leq& C\phi.
\end{align*}
So we get \eqref{4.100}.
\qed

\textbf{Acknowledgement}.  The author would like to thank Prof. Xi-Nan Ma and Prof. Paolo Salani for the joint work \cite{CMS} and the enlightment in this paper. Also, the author would like to express sincere gratitude to Prof. Pengfei Guan for the advice on the choice of coordinates in Dec. 2012.

\end{document}